\documentstyle[cmp]{svjour}

\def\be#1{\begin{equation} \label{#1}}

\def\bi{\begin{itemize}}
\def\R{{\hbox{\bf R}}}

\def\dist{{\hbox{\rm dist}}}
\def\R{{\hbox{\bf R}}}

\def\sgn{{\hbox{\rm sgn}}}

\def\Z{{\hbox{\bf Z}}}

\def\eps{\varepsilon}

\def\lesssim{\mathrel{\hbox{\rlap{\hbox{\lower4pt\hbox{$\sim$}}}\hbox{$<$}}}}
\def\gtrsim{\mathrel{\hbox{\rlap{\hbox{\lower4pt\hbox{$\sim$}}}\hbox{$>$}}}}
\def\eqref#1{({\ref{#1}})}
\def \endprf{\hfill  {\vrule height6pt width6pt depth0pt}\medskip}
\def\emph#1{{\it #1}}
\def\textbf#1{{\bf #1}}
\def\divider#1{$\bullet\quad${\bf #1}}

\begin{document}

\title{Global regularity of wave maps II.  Small energy in two dimensions }
\author{Terence Tao}
\institute{Department of Mathematics, UCLA, Los Angeles CA 90095-1555
\email{tao@math.ucla.edu}}
\date{Nov 26, 2000}

\begin{abstract}
We show that wave maps from Minkowski space $\R^{1+n}$ to a sphere $S^{m-1}$ are globally smooth if the initial data is smooth and has small norm in the critical Sobolev space $\dot H^{n/2}$, in all dimensions $n \geq 2$.  This generalizes the results in the prequel \cite{tao:wavemap1} of this paper, which addressed the high-dimensional case $n \geq 5$.  In particular, in two dimensions we have global regularity whenever the energy is small, and global regularity for large data is thus reduced to demonstrating non-concentration of energy.
\end{abstract}

\maketitle

\section{Introduction}

Throughout this paper $m \geq 2$, $n \geq 1$ will be fixed integers, and all constants may depend on $m$ and $n$.

Let $\R^{1+n}$ be $n+1$ dimensional  Minkowski space  
with flat metric $\eta := \hbox{diag}(-1, 1, \ldots, 1)$, and let $S^{m-1} \subset \R^m$ denote the unit sphere in the Euclidean space $\R^m$.  Elements $\phi$ of $\R^m$ will be viewed as column vectors, while their adjoints $\phi^\dagger$ are row vectors.  We let $\partial_\alpha$ and $\partial^\alpha$ for $\alpha = 0,\ldots, n$ be the usual derivatives with respect to the Minkowski metric $\eta$, subject to the usual summation conventions.  We let $\Box := \partial_\alpha \partial^\alpha = \Delta - \partial_t^2$ denote the D'Lambertian.  We shall write $\phi_{,\alpha}$ and $\phi^{,\alpha}$ for $\partial_\alpha \phi$ and $\partial^\alpha \phi$ respectively.

Define a \emph{classical wave map} to be any function $\phi$ defined on an open set in $\R^{1+n}$ taking values on the sphere $S^{m-1}$ which is smooth, equal to a constant outside of a finite union of light cones, and obeys the equation
\be{wavemap-eq}
\Box \phi = - \phi \phi_{,\alpha}^\dagger \phi^{,\alpha}.
\end{equation}

For any time $t$, we use $\phi[t] := (\phi(t), \partial_t \phi(t))$ to denote the position and velocity of $\phi$ at time $t$.  We refer to $\phi[0]$ as the \emph{initial data} of $\phi$.  Note that in order for $\phi[0]$ to be the initial data for a classical wave map, $\phi[0]$ must be smooth, equal to a constant outside of a compact set, and satisfy the consistency conditions
\be{consist}
\phi^\dagger(0) \phi(0) = 1; \quad \phi^\dagger(0) \partial_t \phi(0) = 0.
\end{equation}
We shall refer to data $\phi[0]$ which satisfy these properties as \emph{classical initial data}.

The purpose of this paper is to prove the following regularity result for classical wave maps.  Let $\dot H^s := (\sqrt{-\Delta})^{-s} L^2(\R^n)$ denote the usual homogeneous Sobolev spaces.

\begin{theorem}\label{main2}
Let $n \geq 2$, and suppose that $\phi[0]$ is classical initial data which has sufficiently small $\dot H^{n/2} \times \dot H^{n/2-1}$ norm.  Then $\phi$ can be extended to a classical wave map globally in time.  Furthermore, if $s$ is sufficiently close to $n/2$, we have the global bounds
\be{bound}
\| \phi[t] \|_{L^\infty_t (\dot H^s_x \times \dot H^{s-1}_x)} \lesssim \| \phi[0] \|_{\dot H^s_x \times \dot H^{s-1}_x}.
\end{equation}
\end{theorem}

In particular, in the energy-critical two-dimensional case one has global regularity for wave maps with small energy.  From this and standard arguments based on finite speed of propagation (see e.g. \cite{christ.spherical.wave}), we see that the problem of global regularity for general smooth data is thus reduced to demonstrating the non-concentration of energy. This non-concentration is known if one assumes some symmetry on the data and some curvature assumptions on the target manifold (\cite{christ.spherical.wave}, \cite{shatah.shadi.noblow}, \cite{struwe:sphere}, \cite{struwe:equivariant}), but is not known in general.  For further discussion on these problems see, e.g. \cite{kman.barrett}, \cite{kman.selberg:survey}, \cite{shatah-struwe}, \cite{struwe.barrett}.  A similar result, but with the Sobolev norm $\dot H^{n/2}$ replaced a slightly smaller Besov counterpart $\dot B^{n/2,1}_2$, was obtained in \cite{tataru:wave2}.  Indeed, our paper shall be largely be a (self-contained) combination of \cite{tataru:wave2} and the prequel \cite{tao:wavemap1} to this paper, although there are some additional technical issues arising here which do not occur in the two papers just mentioned.

Theorem \ref{main2} was proven in \cite{tao:wavemap1} in the 
high-dimensional case $n \geq 5$.  In that paper the main techniques were Littlewood-Paley  decomposition, Strichartz estimates\footnote{Readers familiar with the literature may be surprised that Strichartz estimates are able to handle the critical problem for wave maps.  The reason is that the renormalization almost reduces the strength of the non-linearity to the level of a pure power.  A more precise statement 
is that the renormalization ensures that at in the event of high-low frequency interactions, at least one of the two derivatives in the cubic non-linearity will land on the lowest frequency term.  To compare this with the pure power problem, observe that if we could somehow ensure that \emph{both} derivatives in the non-linearity landed on the two lowest frequency terms, then we could differentiate the equation to obtain something like $\Box \nabla_{x,t} \phi = -\nabla_{x,t} \phi \phi_{,\alpha}^\dagger \phi^{,\alpha}$, which is a cubic semi-linear equation in $\nabla_{x,t} \phi$ with an additional 
null structure.}
and an adapted co-ordinate frame constructed by parallel transport.  To cover the low dimensional cases $n=2,3,4$ we shall keep the Littlewood-Paley decomposition and adapted co-ordinate frame construction (with only minimal changes from \cite{tao:wavemap1}), but we shall abandon the use of Strichartz estimates as the range of these estimates becomes far too restrictive to be of much use, especially when $n=2$.  Instead, we shall adapt the more intricate spaces (including $\dot X^{s,b}$ spaces) and estimates developed in \cite{tataru:wave2}, as substitutes for the Strichartz estimates.  This will make the argument much lengthier and involved, although the overall strategy is little changed\footnote{Indeed, the basic renormalization argument only covers about a third of the paper, from Section \ref{xsb-sec} to Section \ref{linear-sec}.  The bulk of the paper is concerned instead with constructing rather complicated function spaces as substitutes for the Strichartz spaces, and proving the relevant estimates for those spaces.} from that in \cite{tao:wavemap1}.  One major new difficulty is that multiplication by $L^\infty_t L^\infty_x$ functions is not well-behaved on $\dot X^{s,b}$ spaces, and so we will need to replace $L^\infty_t L^\infty_x$ with a more complicated Banach algebra.

In \cite{tao:wavemap1} the non-linearity was placed (after localizing in frequency and switching to the adapted co-ordinate frame) in the familiar space $L^1_t \dot H^{n/2-1}_x$.  When $n \geq 5$ this was relatively easy to achieve, since one had access to $L^2_t L^4_x$ and $L^2_t L^\infty_x$ Strichartz estimates.  For $n=4$ one loses the $L^2_t L^4_x$ estimate, but one could probably use $\dot X^{s,b}$ spaces and null form estimates (which would place $\phi^\dagger_{,\alpha} \phi^{,\alpha}$ in spaces such as $L^2_t L^2_x$) as a substitute, so that one could continue to place the non-linearity in such good spaces as $L^1_t \dot H^{n/2-1}_x$.  When $n=3$ one also (barely) loses the $L^2_t L^\infty_x$ estimate, although in principle this could be compensated for by the $L^p$ null form estimates in \cite{wolff:cone}, \cite{tao:cone} for certain $p < 2$.  However in the energy-critical case $n=2$, the best Strichartz estimate available is only $L^4_t L^\infty_x$, and it appears that even the best possible $L^p$ null form estimates\footnote{An examination of the known counter-examples suggests that it may just be possible to place $\phi_{,\alpha} \psi^{,\alpha}$ in $L^{4/3}_t L^2_x$, which in principle is just barely enough to obtain $L^1_t L^2_x$ control on the non-linearity thanks to the $L^4_t L^\infty_x$ Strichartz estimate.  However this would require (among other things) a reworking of the endpoint argument of \cite{tao:cone} and would therefore not be a simplification to this paper. } are not strong enough to place the non-linearity in a space such as $L^1_t L^2_x$, even after using the adapted co-ordinate frame and introducing $\dot X^{s,b}$ type spaces.

Because of this, we can only place a small portion of the non-linearity in $L^1_t L^2_x$.  Following \cite{tataru:wave2}, we shall place the other portions of the non-linearity either in an $\dot X^{s,b}$ type space, or in $L^1_t L^2_x$ spaces corresponding to null frames.  To obtain these types of control on the non-linearity, we shall use null-form estimates, as well as the decomposition, introduced in \cite{tataru:wave2}, of free solutions as a superposition of travelling waves, each of which is in $L^2_t L^\infty_x$ with respect to a certain null frame.  This decomposition, combined with the $L^2_t L^2_x$ control coming from $\dot X^{s,b}$ estimates, is crucial to recover the $L^1_t L^2_x$ type control of the non-linearity which we need to close the argument.

The high-dimensional argument in \cite{tao:wavemap1} did not need to exploit the null structure in \eqref{wavemap-eq}.  However, one does not have this luxury in the low-dimensional cases, and we shall need in particular to rely on the identity
\be{null-form}
2\phi_{,\alpha} \psi^{,\alpha} = \Box(\phi \psi) - \phi \Box \psi - \Box \phi \psi
\end{equation}
heavily (cf. \cite{tataru:wave1}, \cite{tataru:wave2}, and elsewhere).  This identity is useful when $\phi$, $\psi$, $\phi \psi$ are relatively close to the light cone in frequency space, although when one is far away from the light cone this identity can become counter-productive.

It is quite possible that Theorem \ref{main2} can be extended to other manifolds than the sphere\footnote{In dimensions $n \geq 5$ we have recently learnt (S. Klainerman, personal communication) that this has been achieved by Klainerman and Rodnianski in the case when the target manifold is a group.}, and to scattering and well-posedness results.  We refer the reader to \cite{tao:wavemap1} as we have nothing of interest to add to that discussion here (other than a large increase in complexity).  Indeed, we would strongly recommend to anyone interested in these problems for small data that they first study the high-dimensional case before attempting the low-dimensional one.  (For large data one of course has blowup in dimensions greater than two due to the supercritical nature of the energy conservation law; see \cite{shatah.shadi.blow}.)

This work was conducted at UCLA, Tohoku University, UNSW, and the French Alps.  The author thanks Daniel Tataru, Mark Keel, and Sergiu Klainerman for very helpful discussions, insights, and encouragement, and to Joachim Krieger for pointing out an error in the original manuscript. The author is a Clay Prize Fellow and is supported by grants from the Sloan and Packard foundations.
 
\section{Notation and preliminary reduction}\label{xsb-sec}

We shall need some small exponents
$$0 < \delta_0 \ll \delta_1 \ll \delta_2 \ll \delta_3 \ll \delta_4 \ll 1.$$
The exact choice of these exponents is not important, but for concreteness we shall choose them as follows.  We first choosing $0 < \delta_4 \ll 1$ to be a small absolute constant depending on $n$ ($\delta_4 = 1/100n$ shall do), and then set $\delta_i := \delta_{i+1}^{10}$ for $i=3,2,1,0$. We shall implicitly be inserting the disclaimer ``assuming $\delta_4$ is sufficiently small depending on $n$'' in all the arguments which follow. Thus any exponential term involving $\delta_4$ shall dominate a corresponding term involving $\delta_3$, and so forth down to $\delta_0$, which is dominated by everything.

Let $j, k$ be integers and $i=0,1,2,3,4$.  We use $\chi^{(i)}_{j \leq k}$ or $\chi^{(i)}_{k \geq j}$ to denote a quantity of the form $\min(1, 2^{-\delta(j-k)})$, where $\delta > C^{-1} \delta_i^2$ for some absolute constant $C > 0$ depending only on $n$.  We also use $\chi^{(i)}_{j = k}$ to denote a quantity of the form $2^{-\delta|j-k|}$ with the same assumptions on $\delta$.  Thus $\chi^{(i)}_{j \leq k}$ is small unless $j \leq k + O(1)$, and $\chi^{(i)}_{j = k}$ is small unless $j = k + O(1)$.  

Similarly, we $\chi^{-(i)}_{j \leq k} = \chi^{-(i)}_{k \geq j}$ to denote a quantity of the form $\max(1, 2^{\delta(j-k)})$, where $\delta < C \delta_i^{1/2}$ for some absolute constant $C > 0$ depending only on $n$, and also use $\chi^{-(i)}_{j=k}$ to denote quantity of the form $2^{\delta|j-k|}$ with the same assumptions on $\delta$.  The $\chi^{(i)}$ thus various exponential gains in our estimates, while $\chi^{-(i)}$ represent various exponential losses.  Note that a $\chi^{(i)}$ gain will dominate a corresponding $\chi^{-(j)}$ loss whenever $i > j$. 

As usual, we use $A \lesssim B$ or $A = O(B)$ to denote the estimate $A \leq CB$, where $C$ is some quantity depending only on $n$, $m$, and the $\delta_i$.  All sums will be over the integers $\Z$ unless otherwise specified.

We fix $0 < \eps \ll 1$ to be a small constant depending only on $n$, $m$, and the $\delta_i$ ($\eps := \delta_0^{nm}$ will 
suffice\footnote{Of course, by our above construction this value of $\eps$ is absurdly small, as we have wildly exaggerated the separation of scales in exponents that we shall actually need.  One can improve the value of $\eps$ substantially, but we shall not attempt to do so here.}).  We shall implicitly insert the disclaimer ``assuming $\eps$ is sufficiently small depending on $n$, $m$ and the $\delta_i$'' in all the arguments which follow.  Eventually we shall assume that the initial data has a $\dot H^{n/2} \times \dot H^{n/2-1}$ norm of $\ll \eps^2$. 

We shall parameterize spacetime $\R^{1+n}$ in the standard Euclidean frame $\{(t,x): t \in \R, x \in \R^n\}$ with the Euclidean inner product $(t,x) \cdot (t',x') := tt' + x \cdot x'$; we not use the Minkowski metric $\eta$ much (except in Case 4(e) of Section \ref{o-lemma-sec}).  In the proof of our estimates in the second half of the paper we shall also introduce null frames for $\R^{1+n}$, but we shall not need them for quite a while.

We fix $T>0$ to be a given time.  It will be important that our implicit constants do not depend on $T$.   For the first half of this paper, which is concerned with the iteration scheme and the renormalization, our functions shall only be supported on the slab $[-T,T] \times \S(1)$, but in the second half, which is concerned with the function spaces and the estimates, we shall mainly work on all of $\R^{1+n}$ (as one then gains access to the spacetime Fourier transform) and then apply standard restriction arguments to return to $[-T,T] \times \R^n$. 

We define the Lebesgue spaces $L^q_t L^r_x$ by the norm
$$ \| \phi \|_{L^q_t L^r_x} := (\int (\int |\phi(t,x)|^r\ dx)^{q/r}\ dt)^{1/q}$$
with the usual modifications when $r=\infty$ or $q=\infty$.

If $\phi(t,x)$ is a function on $[-T,T] \times \R^n$ or $\R^{1+n}$, we define the spatial Fourier transform $\hat \phi(t,\xi)$ by
$$ \hat \phi(t,\xi) := \int_{\R^n} e^{-2\pi i x \cdot \xi} \phi(t,x)\ dx$$
with the inverse transform given by
$$ \check F(t,x) := \int_{\R^n} e^{2\pi i x \cdot \xi} F(t,\xi)\ d\xi.$$
The spatial Fourier transform is distinct from the \emph{spacetime} Fourier transform ${\cal F} \phi(\tau,\xi)$, which we shall need in the second half of the paper.  We define the \emph{spatial Fourier support} or \emph{$\xi$-Fourier support} of $\phi$ to be the set $\{ \xi: \hat \phi(t,\xi) \neq 0 \hbox{ for some } t \}$.  

We shall write $D_0$ for $|\xi|$, so that $D_0$ measures the strength of the operator $\nabla_x$.  Thus, for instance, the set $\{ D_0 \sim 2^k\}$ denotes the frequency region $\{ \xi: |\xi| \sim 2^k\}$.  

We now set up some Littlewood-Paley operators, which shall play a central role in our arguments.  Fix $m_0(\xi)$ to be a non-negative radial bump function supported on $D_0 \leq 2$ which equals 1 on the ball $D_0 \leq 1$.  For each integer $k$, define the operators $P_{\leq k} = P_{<k+1}$ to the ball $D_0 \lesssim 2^k$ by the formula
$$ \widehat {\phi_{\leq k}}(t,\xi) := m_0(2^{-k}\xi) \hat \phi(t,\xi),$$
and the projection operators $P_k$ to the frequency annulus $D_0 \sim 2^k$ by the formula
$$ P_k := P_{\leq k} - P_{<k}.$$
Thus $P_k$ has symbol $m(2^{-k} \cdot)$, where $m(\xi) := m_0(\xi) - m_0(2\xi)$.  Note that $m$ is supported on the annulus $2^{-2} \leq D_0 \leq 2^2$.

We also define more general projections $P_{k_1 \leq \cdot \leq k_2}$ and $P_{\geq k}$ by
$$ P_{k_1 \leq \cdot \leq k_2} := P_{\leq k_2} - P_{< k_1}; \quad P_{\geq k} := 1 - P_{<k}.$$
Similarly define $P_{k_1 < \cdot \leq k_2}$, $P_{>k}$, etc.  Note that these operators all have kernels which are measures of finite mass, and so these operators are bounded on all translation-invariant Banach spaces.  Also observe that these operators commute with differentiation operators (or other Fourier multipliers), as well as any time cutoffs.

We shall write $\phi_k$ for $P_k \phi$.  Similarly define $\phi_{\leq k}$, $\phi_{\geq k}$, $\phi_{k_1 < \cdot < k_2}$, etc.  We shall occasionally extend this notation to $\psi$ (e.g. $\psi_{\leq k} := P_{\leq k} \psi$), but other subscripted functions ($U_k$, $w_k$, $A_{\leq k;\alpha}$, etc.) will not be subject to this convention, although the subscripts will often be suggestive.  For instance, in our main argument $A_{\leq k; \alpha}$ and $U_{\leq k}$ will have Fourier support on $D_0 \lesssim 2^k$, and $w_k$ has Fourier support on $D_0 \sim 2^k$.

Differentiation operators and other Fourier multipliers take precedence over arithmetic operations, thus $\Box \phi \psi$ denotes $(\Box \phi) \psi$ rather than $\Box(\phi \psi)$, and $P_k \phi \psi$ denotes $(P_k \phi) \psi$ rather than $P_k(\phi \psi)$.

We shall often exploit the \emph{Littlewood-Paley product trichotomy}, which asserts if $\phi$ and $\psi$ have Fourier support on $D_0 \sim 2^{k_1}$ and $D_0 \sim 2^{k_2}$ respectively, then $P_k(\phi \psi)$ vanishes unless one of the following three statements\footnote{One could also distinguish the ``low-low interaction'' case $k = k_2 + O(1)$; $k_1 = k_2 + O(1)$ from the three listed above.  This is sometimes useful in Schrodinger or KdV contexts as it isolates the parallel interaction case, but is not particularly worthwhile for wave equations because there is essentially no change in the geometry of the light cone in the transition between this case and the other three.} is true:
\begin{itemize}
\item (Low-high interaction) $k = k_2 + O(1)$; $k_1 \leq k_2 + O(1)$.
\item (High-low interaction) $k = k_1 + O(1)$; $k_2 \leq k_1 + O(1)$.
\item (High-high interaction) $k \leq k_2 + O(1)$; $k_1 = k_2 + O(1)$.
\end{itemize}
The same is true if $\phi \psi$ is replaced by similar expressions such as $\phi^{,\alpha} \psi_{,\alpha}$ or $L(\phi, \psi)$ (see below).  In the wave map problem with $n \geq 2$ the low-high (and by symmetry, the high-low) interactions are dominant; the high-high interactions are weaker because in two and higher dimensions it is unlikely that two high frequencies will be so opposed as to cancel and form a low frequency\footnote{Another way to see this is to consider an expression such as $P_0(\phi_{k,\alpha} \psi_k^{,\alpha})$ for $k > 0$.  If $\phi, \psi \in \dot H^{n/2}_x$, then $\phi_{k,\alpha}$ and $\psi_k^{,\alpha}$ are in $L^2_x$ with a norm of $O(2^{-(\frac{n}{2}-1)k})$. For $n > 2$ this gives a decay in $k$.  For $n=2$ one can exploit the null form structure via \eqref{null-form} to obtain a similar heuristic decay.  For $n=1$ no decay is available, even with the null structure.}. This qualitative statement will be made concrete by the presence of factors such as $2^{-|k-k_2|/10}$ in several of our estimates.  This fact implies that the ``high to low'' frequency cascade is heuristically negligible\footnote{More precisely, if we only retained the high-high interaction component of the non-linearity and suppressed the other two components, then one could obtain global regularity, well-posedness, and scattering for small $\dot H^{n/2}$ data for $n \geq 2$ directly from the iteration arguments in \cite{tataru:wave2} (or in \cite{tataru:wave1} for $n \geq 4$).  Indeed, for $n \geq 5$ one could even do this just by iterating in Strichartz spaces as one can see by inspecting the arguments in \cite{tao:wavemap1}.  Note that once the high-high interaction is neglected, the wave map equation essentially becomes an upper-triangular system in the sense that low frequencies affect the high, but not vice versa.  In principle, this allows us to solve the system recursively, and our construction of the gauge $U_{\leq k}$ shall reflect this philosophy.}, and so we only need to concentrate
 on stemming a possible ``low to high'' frequency cascade coming from a divergence of the low-high and high-low interactions.  Our main weapon in achieving this is a ``gauge transformation'' or ``renormalization'', which effectively moves a derivative from the high frequency term to the low one.  We remark that in the one-dimensional case that the ``high to low'' cascade is \emph{not} negligible and is not damped by the renormalization.  Indeed, this cascade causes ill-posedness at the critical regularity \cite{tao:ill}.

Note that if $\phi$ is a smooth function which is equal to a constant $e$ outside of a compact set, then the $\phi_k$ are rapidly decreasing in space and we have the Littlewood-Paley decomposition
$\phi = e + \sum_k \phi_k$.  In particular, we have
$$ \phi_{,\alpha} = \sum_k \phi_{k,\alpha}$$
for all indices $\alpha$.

We shall rely frequently on \emph{Bernstein's inequality}, which asserts (among other things) that
\be{bernstein}
\| f\|_{L^\infty_x} \lesssim 2^{nk/2} \|f\|_{L^2_x}
\end{equation}
whenever $f$ has Fourier support on $D_0 \lesssim 2^k$, and (by duality) that
\be{bernstein-dual}
\| f\|_{L^2_x} \lesssim 2^{nk/2} \|f\|_{L^1_x}
\end{equation}
under the same assumption.  More generally, we have
\be{bernstein-gen}
\| f\|_{L^\infty_x} \lesssim |R|^{1/2} \|f\|_{L^2_x}
\end{equation}
whenever $f$ has Fourier support on a box $R$.
  When estimating a product expression such as $P_k(\phi_{k_1} \psi_{k_2})$ we shall usually apply Bernstein's inequality at the lowest frequency $\min(k,k_1,k_2)$, in order to minimize the factor $2^{nk/2}$ which appears.

We shall now introduce a convenient notation for describing multi-linear expressions of product type.  More precisely, if $\phi^{(1)}(t,x), \ldots, \phi^{(s)}(t,x)$ are scalar functions, we use $L(\phi^{(1)}, \ldots, \phi^{(s)})(t,x)$ to denote any multi-linear expression of the form
$$ L(\phi^{(1)}, \ldots, \phi^{(s)})(t,x) := \int K(y_1, \ldots, y_{S(c)}) \phi^{(1)}(t,x-y_1) \ldots \phi^{(s)}(t, x-y_{S(c)})\ dy_1 \ldots dy_{S(c)}$$
where the kernel $K$ is a measure with bounded mass (the exact value of $K$ might change from line to line).  The kernel of $L$ shall never depend on the index $\alpha$.
Also we extend the notation to the case when $\phi^{(1)}, \ldots \phi^{(s)}$ take values as $m$-dimensional vectors or $m \times m$ matrices, by making $K$ into an appropriate tensor.  For instance, if $\phi^{(1)} = (\phi^{(1), i})_{i=1}^m$ is an $m$-dimensional vector field and $\phi^{(2)} = (\phi^{(2),jk})_{j,k=1}^m$ is an $m \times m$ matrix field, we use $L(\phi^{(1)},\phi^{(2)})(t,x)$ to denote any function of the form
$$ L(\phi^{(1)}, \phi^{(2)})(t,x) = \sum_{i,j,k=1}^m \int K_{ijk}(y_1, y_2) \phi^{(1),i}(t,x-y_1) \phi^{(2),jk}(t, x-y_2)\ dy_1 dy_2$$
where for each $i,j,k$ the kernel $K_{ijk}$ (which may be itself scalar, vector, or matrix valued) has a bounded mass.

This $L$ notation will turn out to be extremely handy for suppressing matrix co-efficients, Littlewood-Paley multipliers, commutator expressions, null forms, etc., whenever these structures are not being exploited (in a manner similar to the big-$O()$ notation, which we shall also use).  The expression $L(\phi^{(1)}, \ldots, \phi^{(s)})$ should be thought of as a variant of the expression $O(\phi^{(1)} \ldots \phi^{(s)})$, and in general $L$ can be treated as if it were just the pointwise product operator.  For instance, the Fourier support of $L(\phi^{(1)}, \ldots, \phi^{(s)})$ is contained in the set-theoretic sum of the Fourier supports of the individual $\phi^{(i)}$.  Also, from Minkowski's inequality we immediately have

\begin{lemma}\label{minkowski}  Let $X_1, \ldots, X_s, X$ be spatially translation-invariant\footnote{By ``$X$ is spatially translation-invariant'' we mean that $f(t,x-x_0)$ has exactly the same $X$ norm as $f(t,x)$ for all $f \in X$ and $x_0 \in \R^n$.} Banach spaces such that we have the product estimate
$$ \|\phi^{(1)} \ldots \phi^{(s)}\|_X \leq A \| \phi^{(1)} \|_{X_1} \ldots \| \phi^{(s)} \|_{X_s}$$
for all scalar-valued $\phi^{(i)} \in X_i$ and some constant $A > 0$.  Then we have
$$ \|L(\phi^{(1)}, \ldots, \phi^{(s)})\|_X \lesssim (Cm)^{Cs} A \| \phi^{(1)} \|_{X_1} \ldots \| \phi^{(s)} \|_{X_s}$$
for all $\phi^{(i)} \in X_i$ which are scalars, $m$-dimensional vectors, or $m \times m$ matrices.  Similarly, if we have the null form estimate
$$ \|\phi^{(1)} \ldots \phi^{(s-2)} \phi^{(s-1)}_{,\alpha} \phi^{(s),\alpha} \|_X \leq A \| \phi^{(1)} \|_{X_1} \ldots \| \phi^{(s)} \|_{X_s}$$
for all scalar $\phi^{(i)} \in X_i$, then 
$$ \|L(\phi^{(1)}, \ldots, \phi^{(s-2)}, \phi^{(s-1),\alpha}, \phi^{(s),\alpha})\|_X \lesssim (Cm)^{Cs} A \| \phi^{(1)} \|_{X_1} \ldots \| \phi^{(s)} \|_{X_s}$$
for all $\phi^{(i)} \in X_i$ as previously.

Finally, if we have the product estimate 
$$ \| \phi \psi \|_X \leq A \| \phi \|_{X_1} \| \psi \|_{X_2}$$
for all scalar $\phi$, $\psi$, then we have
$$ \|L(\phi^{(1)}, \phi^{(2)}, \ldots, \phi^{(s)})\|_X \lesssim Cm A \| \phi^{(1)} \|_{X_1} \| L(\phi^{(2)}, \ldots, \phi^{(s)}) \|_{X_s}$$
for all $\phi^{(i)}$ as previously, where the right-hand side $L$ depends on the left-hand side $L$ and also on the functions $\phi^{(1)}, \ldots, \phi^{(s)}$.
\end{lemma}

In our applications (except in Step 2(c) of Section \ref{linear-sec}) $s$ will usually just be 2 or 3, so the $(Cm)^{Cs}$ term can be entirely neglected. 

The $L$ notation is invariant under permutations of the functions (with the kernel $K$ also being permuted, of course).  It also interacts well with Littlewood-Paley operators (which always have integrable kernels):
\begin{eqnarray}\label{l-lp} 
L(P_k \phi^{(1)}, \phi^{(2)}, \ldots, \phi^{(s)}) &= L(\phi^{(1)}, \phi^{(2)}, \ldots, \phi^{(s)});\\
P_k L(\phi^{(1)}, \phi^{(2)}, \ldots, \phi^{(s)}) &= L(\phi^{(1)}, \phi^{(2)}, \ldots, \phi^{(s)}).\nonumber
\end{eqnarray}
Similarly for $P_{<k}$, $P_{>k}$, etc.  A similar statement holds for derivatives of Littlewood-Paley operators:
\begin{eqnarray}\label{l-lp-d} 
L(\nabla_x P_{\leq k} \phi^{(1)}, \phi^{(2)}, \ldots, \phi^{(s)}) &= 2^k L(\phi^{(1)}, \phi^{(2)}, \ldots, \phi^{(s)}); \\
\nabla_x P_{\leq k} L(\phi^{(1)}, \phi^{(2)}, \ldots, \phi^{(s)}) &= 2^k L(\phi^{(1)}, \phi^{(2)}, \ldots, \phi^{(s)}).\nonumber
\end{eqnarray}
In particular, we have
\be{l-lp-d-2}
L(\nabla_x \phi^{(1)}, \phi^{(2)}, \ldots, \phi^{(s)}) = 2^{k_1} L(\phi^{(1)}, \phi^{(2)}, \ldots, \phi^{(s)})
\end{equation}
whenever $\phi^{(1)}$ has Fourier support on $D_0 \lesssim 2^k$.  Similarly for permutations.

A typical place where the $L$ notation is useful is in commutator expressions.  Specifically, we have:

\begin{lemma}[Leibnitz rule for $P_k$]\label{commutator}  We have the commutator identity
\be{commute}
P_k(fg) = f P_k g + L(\nabla_x f, 2^{-k} g).
\end{equation}
\end{lemma}

\begin{proof}
We may rescale $k=0$.  By the Fundamental theorem of Calculus we have
\begin{eqnarray*}
 (P_0(f g) - f P_0 g)(t,x) &=
\int_{\R^n} \check m(y) (f(t,x-y) - f(t,x)) g(t,x-y)\ dy\\
&=
-\int_0^1 \int_{\R^n} \check m(y) y \cdot \nabla_x f(t,x-sy) g(t,x-y)\ dy ds.
\end{eqnarray*}
The claim follows from the rapid decay of $\check m$. 
\end{proof}

This Lemma is most useful when $g$ has frequency $\sim 2^k$ and $f$ has frequency $\lesssim 2^k$.  In this case $P_k(fg) - f P_k g$ effectively shifts a derivative from the high-frequency function $g$ to the low-frequency function $f$.  This shift will generally ensure that all such commutator terms will be easily estimated.

We now recall a definition, essentially from \cite{tao:wavemap1}:

\begin{definition}\label{envelope-def}  A \emph{frequency envelope} is a sequence $c = \{c_k\}_{k \in \Z}$ of positive reals such that we have the $l^2$ bound
\be{l2-ass}
\|c\|_{l^2} \lesssim \eps
\end{equation}
and the local constancy condition
\be{local}
\chi^{(0)}_{k=k'} c_{k'} \lesssim c_k \lesssim \chi^{-(0)}_{k=k'} c_{k'}
\end{equation}
for all $k, k' \in \Z$.  In particular we have $c_k \sim c_k'$ whenever $k = k' + O(1)$.  If $c$ is a frequency envelope and $(f,g)$ is a pair of functions on $\R^n$, we say that $(f,g)$ \emph{lies underneath the envelope} $c$ if one has
\be{envelope}
\| P_k f \|_{\dot H^{n/2}} + \| P_k g \|_{\dot H^{n/2-1}} \leq c_k
\end{equation}
for all $k \in \Z$.
\end{definition}

One should think of \eqref{local} as asserting that $c_k$ is effectively constant.  In practice, the factors of $\chi^{\pm(0)}_{k=k'}$ shall always be dominated by factors such as $\chi^{(i)}_{k=k'}$ for $i > 0$.

We shall not use the full strength of \eqref{l2-ass} often, and shall usually rely instead on the weaker estimate 
\be{linfty-ass}
c_k \lesssim \eps
\end{equation}
for all integers $k$.  However, there will be occasions (especially in the construction of the co-ordinate frame $U$) in which we will have to fully exploit \eqref{l2-ass}.  Because of this, our argument will not apply to any Besov space beyond $\dot H^{n/2}$.

To prove Theorem \ref{main2}, we shall first prove

\begin{theorem}\label{main-envelope}
Let $T_0 > 0$ and $c$ be a frequency envelope, and suppose that $\phi$ is a classical wave map on $[-T_0,T_0] \times \R^n$ such that $\phi[0]$ lies underneath $\eps c$.  Then $\phi[t]$ lies underneath $C c$ for all $t \in [-T_0,T_0]$, where $C$ is an absolute constant depending only on $n$, $m$.
\end{theorem}

We now show how Theorem \ref{main2} follows from Theorem \ref{main-envelope}.  From the regularity theory in \cite{kman.selberg}, \cite{tataru:wave2} it suffices to show that \eqref{bound} holds in $[-T_0,T_0]$ for all classical wave maps $\phi$ on $[-T_0,T_0] \times \R^n$ whose initial data $\phi[0]$ has $\dot H^{n/2} \times \dot H^{n/2-1}$ norm  $\ll \eps^2$.  But if $\phi$ is such a wave map, we can define the envelope $c$ by
$$
c_k := \eps^{-1} \sum_{k' \in \Z} 2^{-\delta_0 |k-k'|} \| \phi_{k'}[0] \|_{\dot H^{n/2} \times \dot H^{n/2-1}}.$$
It is easy to verify \eqref{l2-ass}, \eqref{local}, \eqref{envelope}, so that $\phi[0]$ lies underneath $\eps c$.  Also, from Young's inequality on $l^2$ we have 
$$ (\sum_k  (2^{\sigma k} c_k)^2)^{1/2} \sim \| \phi[0] \|_{\dot H^{n/2+\sigma} \times \dot H^{n/2-1+\sigma}}$$
for all $|\sigma| \leq \delta_0$.  From Theorem \ref{main-envelope} we thus have that $\phi[t]$ lies underneath $C c$ for $t \in [-T_0,T_0]$, which implies again from \eqref{local} that
$$ \| \phi[t] \|_{\dot H^{n/2+\sigma} \times \dot H^{n/2-1+\sigma}} 
\lesssim (\sum_k  (2^{\sigma k} c_k)^2)^{1/2}.$$
This implies \eqref{bound}, and global regularity then follows.

Henceforth the envelope $c$ will be fixed.  The rest of the paper will be devoted to the proof of Theorem \ref{main-envelope}.
We remark that if the $l^2$ control in \eqref{l2-ass} was strengthened to $l^1$ then the proof of this Theorem is essentially contained in \cite{tataru:wave2}, and is based on a pure iteration argument in a rather intricate Banach space.  Our argument will also have an iterative flavor and uses similar spaces, but also requires the renormalizations in \cite{tao:wavemap1}.

\section{The ``iteration'' space, and key estimates}\label{iteration-sec}

In order to prove Theorem \ref{main-envelope} we shall need some Banach spaces $S(c)$, $S_k$, $N_k$ to ``iterate'' in.  The purpose of this section is to describe the important properties and estimates of these spaces; their exact construction and the verification of these properties will be deferred to the second half of the paper (Sections \ref{overview-sec}-\ref{o-lemma-sec}) as they are rather technical and lengthy. 

The space $S(c)$ shall contain the wave map $\phi$ and the gauge transform $U_{\leq k}$.  The spaces $S_k$ shall contain the Littlewood-Paley  pieces $\phi_k$ of $\phi$, as well as the renormalizations $w_k := U_{\leq k-10} \phi_k$ of these pieces.  Finally, $N_k$ contains the non-linearity $\Box \phi_k$ and similar expressions (basically any algebraic combination of $U$, $\phi$ with two derivatives and a frequency of $\sim 2^k$ should be placed inside $N_k$).

Before we describe the properties of these spaces in detail, we first give some motivation.  The proof of Theorem \ref{main-envelope} shall be based on the following (informal) scheme.

\begin{itemize}
\item (Bootstrap hypothesis) By a continuity argument, we shall assume a priori that the $\phi_k$ are already in $S_k$, and then bootstrap this to better control on $\phi_k$ in $S_k$.  
\item (Control of $\phi$) Since $\phi$ lies on the sphere, it will be in $L^\infty_t L^\infty_x$.  Combining this with the assumption $\phi_k \in S_k$, we shall obtain $\phi \in S(c)$.
\item (Construction of gauge) We then construct for each frequency $2^k$ a gauge transform $U_{\leq k}$, which turns out to be a polynomial in $\phi$, $\phi_k$.  Using the above control on $\phi_k$, $\phi$, we shall show\footnote{Actually, we shall first prove that $U_{\leq k} \in S(Cc)$ for some constant $C$, but the spaces $S(Cc)$ and $S(c)$ will end up being equivalent.} that $U_{\leq k} \in S(c)$ for all $k$. 
\item (Control of renormalized non-linearity) We define $w_k := U_{\leq k-10} \phi_k$, and use the above estimates to show that $\Box w_k \in N_k$.
\item (Energy estimates) We also have $w_k[0] \in \dot H^{n/2} \times \dot H^{n/2-1}$.  Using this and the previous step, we shall place $w_k \in S_k$.
\item (Inverting the gauge) We then show that $U_{\leq k-10}$ is invertible in $S(c)$, and use this and the previous step to show that $\phi_k = U_{\leq k-10}^{-1} w_k$ is in $S_k$, thus closing the bootstrap circle.
\item (Epilogue) Finally, we show that control of $\phi_k \in S_k$ implies that $\phi$ lies underneath the envelope $Cc$.
\end{itemize}

To make this scheme work we need several estimates connecting $S_k$, $N_k$, and $S(c)$, which we shall describe in detail in Theorem \ref{spaces}.
One of the major tasks in proving Theorem \ref{main-envelope} is thus to select spaces $S(c)$, $S_k$, $N_k$ which obey all of these estimates.

In the high-dimensional case \cite{tao:wavemap1}, $N_k$ was just the energy method space $L^1_t \dot H^{n/2-1}_x$ (localized to frequency $2^k$), and the spaces $S_k$ were Strichartz spaces.  The analogue of the $S(c)$ norm in \cite{tao:wavemap1} was essentially given by 
$$ \| \phi \|_{S(c)} \sim \|\phi\|_{L^\infty_t L^\infty_x} + \sup_k c_k^{-1} \| \phi_k \|_{S_k}.$$
The properties in Theorem \ref{spaces} turn out to be easily verifiable (after some minor modifications) from standard Strichartz estimates and H\"older's inequality as long as $n \geq 5$.

In the low-dimensional case one cannot place the non-linearity just in $L^1_t \dot H^{n/2-1}_x$, even after renormalization.  Thus we must enlarge the space $N_k$ somewhat, incorporating not only $L^1_t \dot H^{n/2-1}_x$ but also some $\dot X^{s,b}$ type norms, as well as some more complicated null frame spaces of Tataru \cite{tataru:wave2}.  Because we wish to prove energy estimates, this enlargement of $N_k$ also causes $S_k$ to become more complicated; it shall also involve $\dot X^{s,b}$ norms and null frame spaces.  This also forces the space $S(c)$ to become a more sophisticated Banach algebra - something of roughly the same strength as $\dot X^{n/2,1/2}$ but still controlling $L^\infty_t L^\infty_x$ and closed under multiplication.  Fortunately, examples of such algebras $S(c)$ exist (see \cite{kl-mac:algebra}, \cite{tataru:wave1}, \cite{tataru:wave2}, \cite{tao:algebra}).  The space we use shall be loosely based on that in \cite{tataru:wave2}.

We now summarize all the properties of the spaces $S(c)$, $S_k$, $N_k$ that we shall need, and discuss these properties in turn.

\begin{theorem}\label{spaces}  Let $T \geq 0$ and $c$ be a frequency envelope.  Then there exist Banach spaces $S(c) = S(c)([-T,T] \times \R^n)$, $S_k = S_k([-T,T] \times \R^n)$, and $N_k = N_k([-T,T] \times \R^n)$ for $k \in \Z$ of functions in $[-T,T] \times \R^n$, which satisfy the following properties for all integers $k$, $k_1$, $k_2$, $k_3$.

\begin{itemize}

\item (Quasi-continuity)
Suppose that $\phi$ is a classical wave map on $[-T_0,T_0] \times \R^n$.  Then the function 
\be{upstairs}
a(T) := \max(1,c_k^{-1} \sup_k\| \phi_k|_{[-T,T] \times \R^n} \|_{S_k([-T,T] \times \R^n)})
\end{equation}
defined on $[0,T_0]$ obeys the quasi-continuity property 
\be{continuous}
\limsup_{T' \to T} a(T') \lesssim \liminf_{T' \to T} a(T')
\end{equation}
for all $0 \leq T' \leq T_0$.

\item (Invariance properties)  
The spaces $S(c)$, $S_k$, $N_k$ are invariant under spatial translations.  For any $j \in \Z$, the scaling\footnote{In other words, the spaces $S(c)$, $S_k$ are dimensionless, while $N_k$ has the units $length^{-2}$.} $\phi(t,x) \mapsto \phi(2^j t,2^j x)$, $T \mapsto T/2^j$ maps $S(\{c_k\})$ to $S(\{c_{k-j}\})$ and $S_k$ to $S_{k+j}$.  Similarly, the scaling $F(t,x) \mapsto 2^{2j} F(2^j t, 2^j x)$, $T \mapsto T/2^j$ maps $N_k$ to $N_{k+j}$.

\item ($S(c)$ is an algebra) 
The space $S(c)$ is a Banach algebra with respect to pointwise multiplication, i.e. $S(c)$ contains the identity 1 and obeys the estimate
\be{algebra}
\| L(\phi, \psi) \|_{S(c)} \lesssim \| \phi \|_{S(c)} \| \psi \|_{S(c)}
\end{equation}
for all $\phi, \psi \in S(c)$.  Also, we have
\be{infty-control}
\| \phi \|_{L^\infty_t L^\infty_x} \lesssim \| \phi \|_{S(c)}
\end{equation}
for all $\phi \in S(c)$.

\item (Frequency-localized algebra property)  
If $\phi, \psi \in S_k$, then
\be{sk-sk}
\| L(\phi, \psi) \|_{S_k} \lesssim \| \phi \|_{S_k} \| \psi \|_{S_k}.
\end{equation}
Similarly, if $\phi \in S_k$, $\psi \in S(c)$ and $\psi$ has Fourier support in $D_0 \lesssim 2^k$, then
\be{sk-skp}
\| L(\phi, \psi) \|_{S_k} \lesssim \| \phi \|_{S_k} \| \psi \|_{S(c)}.
\end{equation}

\item ($S(c)$ insensitive to $c$) For all $\phi \in S(c)$ and $C > 0$ we have
\be{insensitive}
\| \phi \|_{S(Cc)} \sim \| \phi \|_{S(c)}
\end{equation}
with the implicit constants depending at most polynomially on $C$.

\item ($S(c)$ is built up from $S_k$)
Let $\phi$ be a smooth function on $[-T,T] \times \R^n$ which is constant outside of a compact set.  Suppose we have a decomposition $\phi = \sum_k \phi^{(k)}$, where each $\phi^{(k)}$ is in $S_k$.  Then we have
\be{s-sk}
\| \phi \|_{S(c)} \lesssim \| \phi \|_{L^\infty_t L^\infty_x} + \sup_k c_k^{-1} \| \phi^{(k)} \|_{S_k}.
\end{equation}

\item ($N_k$ contains $L^1_t \dot H^{n/2-1}_x$) 
Let $F$ be an $L^1_t L^2_x$ function on $[-T,T] \times \R^n$ which has Fourier support on the region $D_0 \sim 2^k$ for some integer $k$.  Then $F$ is in $N_k$ and
\be{l12}
\| F \|_{N_k} \lesssim \| F \|_{L^1_t \dot H^{n/2-1}_x} \sim 2^{(n-1)k} \| F \|_{L^1_t L^2_x}.
\end{equation}

\item (Adjacent $N_k$ are equivalent) 
We have the compatibility property
\be{compat}
\| F \|_{N_{k_1}} \sim \| F \|_{N_{k_2}}
\end{equation}
whenever $F \in N_{k_2}$ and $k_1 = k_2 + O(1)$.

\item (Energy estimate) 
For any Schwartz function $\phi$ on $[-T,T] \times \R^n$ with Fourier support in $D_0 \sim 2^k$, we have
\begin{eqnarray}\label{energy-est-2} 
\|\phi \|_{S_k} 
&\lesssim \| \Box \phi \|_{N_k} + 
\| \phi[0] \|_{\dot H^{n/2} \times \dot H^{n/2-1}}\\
&\sim \|\Box \phi \|_{N_k} + 2^{nk} \| \phi(0) \|_{L^2} + 2^{(n-1)k} \| \partial_t \phi(0) \|_{L^2}.\nonumber
\end{eqnarray}

\item (Product estimates) 
We have
\be{ur-algebra}
\| P_k L(\phi,  F) \|_{N_k} \lesssim \chi^{(1)}_{k \geq k_2} \| \phi \|_{S(c)} \| \psi \|_{N_{k_2}}
\end{equation}
for all $\phi \in S(c)$ and $F \in N_{k_2}$.  We also have the variant
\be{ur-algebra-sk}
\| P_k L(\phi, F) \|_{N_k} \lesssim \chi^{(1)}_{k \geq k_2} \| \phi \|_{S_{k_1}} \| F \|_{N_{k_2}}
\end{equation}
whenever $\phi \in S_{k_2}$ and $F \in N_{k_2}$.

\item (Null form estimates)
We have 
\be{null}
\| P_k L(\phi_{,\alpha}, \psi^{,\alpha}) \|_{N_k} \lesssim \chi^{(1)}_{k = \max(k_1,k_2)} \| \phi\|_{S_{k_1}} \| \psi \|_{S_{k_2}}
\end{equation}
for all $\phi \in S_{k_1}$, $\psi \in S_{k_2}$.

\item (Trilinear estimate)
We have
\be{o-lemma}
\| P_k L(\phi^{(1)}, \phi_{,\alpha}^{(2)}, \phi^{(3),\alpha}) \|_{N_k}
\lesssim  
\chi^{(1)}_{k = \max(k_1,k_2,k_3)}
\chi^{(1)}_{k_1 \leq \min(k_2,k_3)}
\prod_{i=1}^3
\| \phi^{(i)} \|_{S_{k_i}} 
\end{equation}
whenever $\phi^{(i)} \in S_{k_i}$ for $i=1,2,3$.

\item (Epilogue)  
For any $\phi \in S_k$ with Fourier support in $D_0 \lesssim 2^k$ we have
\be{energy-est}
\sup_t \| \phi[t] \|_{\dot H^{n/2}_x \times \dot H^{n/2-1}_x} \lesssim 2^{nk/2} \sup_t \| \phi[t] \|_{L^2_x \times L^2_x} \lesssim \| \phi \|_{S_k}.
\end{equation}
\end{itemize}
\end{theorem}

We now discuss each of the above properties in turn.  

\begin{itemize}

\item The estimate \eqref{continuous} is a technical fact needed in order to make the continuity argument work, and is proven in Section \ref{continuity-sec}, mainly using \eqref{energy-est-2} and \eqref{l12}.  Since we are assuming $\phi$ to be smooth and compact outside of a compact set, one would certainly expect the function $a$ to actually be continuous rather than just quasi-continuous, but we do not know how to prove this and in any event it is not needed for our argument.  In the high dimensional case this estimate was trivial as the spaces $S_k$ were just Lebesgue spaces, but more care is required here because $S_k$ will be defined by restriction from $\R^{1+n}$ and have a spacetime Fourier component in their norms.  We remark that the quantity $a(T)$ is necessarily finite for classical wave maps $\phi$, thanks to \eqref{energy-est-2} and \eqref{local}.

\item The invariance properties is unsurprising given the translation and scaling symmetries of the equation, and will be automatic from our construction of $S(c)$, $S_k$, $N_k$ in Section \ref{construction-sec}.  As a corollary of translation invariance we observe that the Littlewood-Paley  operators $P_k$, $P_{\leq k}$, etc. are bounded on the spaces $S(c)$, $S_k$, $N_k$.

\item The algebra property \eqref{algebra} is essential for us to invert the gauge transformation, and will be proven in Section \ref{algebra-sec}.  The spaces described in \cite{tataru:wave2} obey this algebra property if $c \in l^1$, but when $c \in l^2$ there is a logarithmic divergence in the estimates.  Fortunately, this divergence can be rectified (with some non-trivial effort) by adding $L^\infty_t L^\infty_x$ control to $S(c)$.  This is analogous to the well-known fact that $\dot H^{n/2}_x$ is not closed under multiplication, but $\dot H^{n/2}_x \cap L^\infty_x$ is.  The estimate \eqref{infty-control} thus will be an automatic consequence of our construction of $S(c)$ in Section \ref{algebra-sec}.  We shall be able to obtain \eqref{algebra} to some extent from \eqref{ur-algebra} via a convenient duality argument.

\item The estimates \eqref{sk-skp}, \eqref{sk-sk} are minor variants of \eqref{algebra}; indeed, all three estimates shall be treated in a unified manner in Section \ref{algebra-sec}.

\item The insensitivity property \eqref{insensitive} will be immediate from the construction of $S(c)$.  This property is required because it will turn out for induction purposes that it is more convenient to initially measure $U$ in $S(Cc)$ instead of $S(c)$.

\item The estimate \eqref{s-sk} shall turn out to be automatic, because we shall essentially use \eqref{s-sk} to define the space $S(c)$ in Section \ref{construction-sec}.

\item The estimate \eqref{l12} plays only a minor role in the main argument, ensuring that $N_k$ does indeed contain test functions and certain error terms.  Note that the space $L^1_t \dot H^{n/2-1}_x$ is the classical space which one would use to hold the non-linearity, if one attempted to apply the energy method (although this method of course fails at the critical regularity).  This estimate shall be an automatic consequence of our construction of $N_k$ in Section \ref{construction-sec}.

\item The compatibility property \eqref{compat} allows us to ignore certain technical ``frequency leakage'' problems arising from the fact that $\phi \psi$ does not quite have the same frequency as $\phi$, even when $\phi$ has much higher frequency than $\psi$.  It will be an automatic consequence of the construction of $N_k$ in Section \ref{construction-sec}.  A similar property for $S_k$ holds but will not be needed in our argument.  From \eqref{compat} and Littlewood-Paley  decomposition we observe the estimate 
\be{nk-split}
\| F \|_{N_k} \lesssim \sum_{k' = k + O(1)} \| P_{k'} F \|_{N_{k'}}
\end{equation}
whenever $F$ is supported on the region $D_0 \sim 2^k$.

\item The energy estimate \eqref{energy-est-2} is a bit lengthy, and is proven in Section \ref{nk-sec}.  One could try to make \eqref{energy-est-2} the definition of $S_k$, as is done in some other papers, but this makes the product estimate \eqref{ur-algebra} difficult to prove.

\item The estimates \eqref{ur-algebra}, \eqref{ur-algebra-sk} shall be proven  in Section \ref{ur-sec}. The factor $\chi^{(1)}_{k \geq k_2}$ is an indication that the high-high interactions in this problem are quite weak.  (A similar gain is implicit in \cite{tataru:wave2}).

\item We shall prove \eqref{null} in Section \ref{null-sec}.  The proof basically uses the estimates \eqref{ur-algebra-sk}, \eqref{algebra} described above, combined with the identity \eqref{null-form}. In practice we shall only apply \eqref{null} in the high-high interaction case (since we then obtain an exponential gain from $\chi^{(1)}_{k = \max(k_2,k_3)}$), or if a derivative has been transferred from the high-frequency term to the low-frequency term\footnote{Because of this, it is possible to 
lose a factor of up to (but not including) $2^{|k_1-k_2|}$ in \eqref{null} without affecting
 the argument.  This is for instance the case in the $n \geq 5$
 theory in \cite{tao:wavemap1}, where the high frequency term is estimated using the endpoint Strichartz space $L^2_t L^{2(n-1)/(n-3)}_x$ and the low frequency term in the companion space $L^2_t L^{n-1}_x$, thus losing a factor of $2^{|k_1-k_2|(n+1)/2(n-1)}$.}.
From \eqref{nk-split} we observe that the $P_k$ projection in the above lemma can be removed if the expression inside the $P_k$ already has frequency $\sim 2^k$.

\item The trilinear estimate \eqref{o-lemma} is the most difficult estimate in this Theorem to prove, and is handled in Section \ref{o-lemma-sec}.  The factor $\chi^{(1)}_{k = \max(k_1,k_2,k_3)}$ again reflects the fact that high-high interactions are weak.  The difficulty lies primarily in obtaining the small but crucial factor\footnote{For $n \geq 4$ this estimate can be obtained by estimating the two low frequencies in $L^2_t L^\infty_x$ and the high frequency in $L^1_t L^2_x$, and by moving these exponents around by an epsilon one can also cover the $n=3$ case by Strichartz estimates.  However in the $n=2$ case the Strichartz estimates are far too weak to prove this estimate, and we shall need to work much harder.} of $\chi^{(1)}_{k_1 \leq \min(k_2,k_3)}$.  Without this factor, \eqref{o-lemma} essentially follows from \eqref{ur-algebra-sk} and \eqref{null}.  The presence of this factor allows us to handle any non-linearity of cubic or higher degree in which at least one derivative falls on a low frequency term.  In order to obtain this key exponential gain we have to go beyond the arguments in \cite{tataru:wave2} and apply some other tools, notably some multiplier calculus to shift null forms from one function to another, and the use of Bernstein's inequality when the null forms are too degenerate for the multiplier calculus to be effective.
As with \eqref{null}, we remark that the $P_k$ can be removed if the expression inside the $P_k$ already has frequency $\sim 2^k$.

\item The estimate \eqref{energy-est} is basically a dual to \eqref{l12}, and shall be automatic from our construction of $S_k$ in Section \ref{algebra-sec}.  When $n=2$ it might be possible to use energy conservation to circumvent the need for this estimate, however this does not seem to achieve any substantial simplification in this paper.  From \eqref{energy-est} and Bernstein's inequality \eqref{bernstein} we observe the useful estimate
\be{infty-est}
\| \phi \|_{L^\infty_t L^\infty_x} \lesssim \| 2^{-k} \nabla_{x,t} \phi \|_{L^\infty_t L^\infty_x} \lesssim \| \phi \|_{S_k}
\end{equation}
whenever $\phi \in S_k$ has frequency $\sim 2^k$.
\end{itemize}

To close this section we informally discuss how the bilinear and trilinear estimates \eqref{ur-algebra}, \eqref{null}, \eqref{o-lemma} are to be used.  They cannot quite treat the original non-linearity $L(\phi, \phi_{,\alpha}, \phi^{,\alpha})$ in \eqref{wavemap-eq}, especially when the derivatives $\partial_\alpha$, $\partial^\alpha$ fall on high-frequency components of $\phi$.  However these estimates can treat these types of expressions when the derivatives are in more favorable locations.  Examples of such ``good'' non-linearities include

\begin{itemize}
\item (Derivative falls on a low frequency) An expression $L(\phi_{k_1}, \phi_{k_2,\alpha}, \phi^{,\alpha}_{k_3})$ with $k_1 \geq \min(k_2, k_3) + O(1)$.  For these expressions we use \eqref{o-lemma}.
\item (High-high interactions) An expression $L(\phi_{k_1}, \phi_{k_2,\alpha}, \phi^{,\alpha}_{k_3})$ with $k_2 = k_3 + O(1)$.  For these expressions we use \eqref{null} and \eqref{ur-algebra}.
\item (Derivative shifted from high-frequency to low, Type I) An expression of the form
$$ L(\nabla_x \phi_{k_1}, \phi_{k_2,\alpha}, \partial^\alpha \nabla_x^{-1} \phi_{k_3})$$
with $k_1 \leq k_3 + O(1)$.
This generalizes the high-high interaction non-linearity, and arises from commutator expressions via Lemma \ref{commutator}.  This non-linearity is treated by \eqref{o-lemma}.
\item (Derivative shifted from high-frequency to low, Type II) An expression of the form
$$
L(\phi_{k_1}, \nabla_x \phi_{k_2, \alpha}, \nabla_x^{-1} \phi_{k_3}^{,\alpha})$$
with $k_2 \leq k_3 + O(1)$. This type of non-linearity also arises from commutator expressions via Lemma \ref{commutator}, and is estimated by \eqref{null} and \eqref{ur-algebra}.
\item (Repeated derivatives avoiding the highest frequency) An expression of  the form $L(\Box \phi_{k_1}, \phi_{k_2}, \ldots, \phi_{k_s})$ with $k_1 \leq \max(k_2, \ldots, k_s) + O(1)$.  These types of expressions will arise once we apply the gauge transformation $U$.  In principle, one can use the equation \eqref{wavemap-eq} to break this expression up into combinations of the previous types of good non-linearity, although the computations become somewhat tedious in practice.
\end{itemize}

Note that in all cases we have retained the null structure of the non-linearity.  In the low dimensional cases $n=2,3,4$ this is vital to the above non-linearities being good.  In all of the above cases we obtain various exponential gains which will allow us to sum in the $k_i$ indices.

As a first approximation, one should treat these good non-linearities as being negligible errors.  The objective is then to gauge transform (Littlewood-Paley localized versions of) \eqref{wavemap-eq}, exploiting such geometric identities as $\phi^\dagger \phi_{,\alpha} = 0$ as well as Lemma \ref{commutator}, until all the non-linearities are negligible.  In this strategy the Littlewood-Paley  decomposition seems to play an indispensable role, as this decomposition allows us to easily separate the core component of the non-linearity (which for wave maps is a connection term where the connection $A_{\alpha; \leq k}$ has small curvature) from the remaining error terms which are good non-linearities and therefore negligible.

\section{The main proposition}\label{main-sec}

Let $S(c)$, $S_k$, $N_k$ be as in Theorem \ref{spaces}.  We now adapt the argument from \cite{tao:wavemap1}.

In the next section we shall prove the following ``bootstrap'' property of the $S_k$ norms:

\begin{proposition}[Main Proposition]\label{reduced}  Let $c$ be a frequency envelope, $0 < T < \infty$, and let $\phi$ be a classical wave map on $[-T,T] \times \R^n$, extended to $\R^{1+n}$ by the free wave equation, such that $\phi[0]$ lies underneath $\eps c$, and that 
\be{control-2}
\| \phi_k \|_{S_k} \lesssim c_k
\end{equation}
for all $k$.  Then we have
\be{control}
\| \phi_k \|_{S_k} \leq c_k
\end{equation}
holds for all $k$.
\end{proposition}

We now give the continuity argument which deduces Theorem \ref{main-envelope}  from this Proposition.

Let $T_0$, $c$, $\phi$ be as in Theorem \ref{main-envelope}, and let $a(T)$ be the quantity in \eqref{upstairs}.  From \eqref{energy-est-2} and the hypothesis that $\phi$ lies underneath $\eps c$ we see that $a(0) = 1$.  From Proposition \ref{reduced} we see that if $0 < T \leq T_0$ obeys $a(T) \lesssim 1$, then we can automatically bootstrap this bound to $a(T) = 1$.  From this and \eqref{continuous} we see that the set $\{ T \in [0,T_0]: a(T) = 1 \}$ is both open and closed in $[0,T]$.  Since this set contains the origin, we thus have $a(T_0) = 1$.  From this and \eqref{energy-est} we thus see that $\phi[t]$ lies underneath $Cc$ for all $0 \leq t \leq T_0$, as desired.

It only remains to prove Proposition \ref{reduced}.

\section{Renormalized iteration: The proof of Proposition \ref{reduced}}\label{linear-sec}

We shall divide this proof into several steps\footnote{We have decided to organize this paper as a tree as opposed to the more usual linear structure.  Hopefully this shall help keep the ``high-level'' ideas of the argument from being obscured by the dozens of cases and sub-cases which we shall eventually have to consider.}.

\divider{Step 0.  Scaling.}

Fix $c$, $T$, $\phi$, and suppose that the hypotheses of Proposition \ref{reduced} hold.  In this section all our functions and equations shall be on the slab $[-T,T] \times \R^n$. 

Since $\phi$ is on the sphere, it is bounded in $L^\infty_t L^\infty_x$.  From this and \eqref{s-sk} we have the bound
\be{phi-bomb}
\| \phi \|_{S(c)} \lesssim 1.
\end{equation}
Of course, the same bound then holds for all Littlewood-Paley projections of $\phi$, such as $\phi_k$, $P_{\leq k} \phi$, $P_{\geq k} \phi$, etc.

We need to show \eqref{control}.  By scale-invariance (scaling $T$, $c$, and $\phi$ appropriately) it suffices to show that
\be{psi-bound}
\| \phi_0 \|_{S_0} \leq c_0.
\end{equation}

By applying $P_0$ to \eqref{wavemap-eq} we obtain the evolution equation for $\phi_0$:
\be{psi-eq}
\Box \phi_0 = - P_0(\phi \phi_{,\alpha}^\dagger \phi^{,\alpha}).
\end{equation}

\divider{Step 1.  Linearize the $\phi_0$ evolution.}

\begin{definition}\label{ak-def}  For each integer $k$, define the connection $A_{\leq k; \alpha}$ by the formula\footnote{Morally speaking, $A_{\leq k;\alpha}$ is the connection on the tangent bundle of the sphere induced by $\phi_{\leq k}$; however this is a little inaccurate because $\phi_{\leq k}$ does not quite lie on the sphere.}
\be{a-def}
A_{\leq k; \alpha} := A_{<k+1; \alpha} := \phi_{\leq k,\alpha} \phi_{\leq k}^\dagger - \phi_{\leq k} \phi_{\leq k,\alpha}^\dagger.
\end{equation}
\end{definition}

\begin{definition}\label{error-def}  A function $F$ on $[-T,T] \times \R^n$ is said to be an \emph{acceptable error} if
$$ \|F\|_{N_0} \lesssim \eps c_0,$$
and we shall write $F = error$ to denote this.
\end{definition}

The purpose of this step is to convert the non-linear equation \eqref{psi-eq} into the linear transport equation (cf. \cite{tao:wavemap1})
\be{cancel}
\Box \phi_0 = 2 A_{\leq -10; \alpha} \phi_0^{,\alpha} + error.
\end{equation}

In order to convert \eqref{psi-eq} into \eqref{cancel} we need to show that
\be{error}
P_0(\phi \phi_{,\alpha}^\dagger \phi^{,\alpha}) + 2 A_{\leq -10; \alpha} \phi_0^{,\alpha}
\end{equation}
is an acceptable error.  We shall do this into two stages.

\divider{Step 1(a).  Decompose \eqref{error} into small pieces.} 

The purpose of this step is to prove the identity

\begin{eqnarray}
\eqref{error}
=
& \sum_{k_2,k_3: k_2 \geq O(1); k_3 = k_2 + O(1)}
P_0 L(\phi, \phi_{k_2,\alpha}, \phi^{,\alpha}_{k_3})\label{dec1}\\
+&
\sum_{k_2,k_3: k_3 \geq k_2; k_3 \geq O(1) }
P_{-10 < \cdot < 10} L(\phi_{> k_2 - 30}, \phi_{k_2,\alpha},
\phi_{k_3}^\alpha) \label{dec2}\\
+& 
L(\nabla_x \phi_{\leq -10}, \phi_{\leq -10,\alpha}, \phi_{-10 < \cdot < 10}^{,\alpha})\label{fire1}\\
+&
L(\phi_{\leq -10}, \nabla_x \phi_{\leq -10,\alpha}, \phi_{-10 < \cdot < 10}^{,\alpha}).\label{fire2} 
\end{eqnarray}
In the first two summations the kernel of $L$ may depend on $k_2$, $k_3$.

To obtain the above decomposition, we split
\begin{eqnarray*}
\phi \phi_{,\alpha}^\dagger \partial^\alpha \phi 
= &\sum_{k_2,k_3: \max(k_2, k_3) \geq 10; |k_2-k_3| < 5} \phi
\phi_{k_2,\alpha}^\dagger \phi^{,\alpha}_{k_3} \\
+&
\sum_{k_2,k_3: k_2 \geq 10, k_3 + 5}
\phi \phi_{k_2,\alpha}^\dagger \phi^{,\alpha}_{k_3}\\
+&
\sum_{k_2,k_3: k_3 \geq 10, k_2 + 5}
\phi \phi_{k_2,\alpha}^\dagger \phi^{,\alpha}_{k_3}\\
+&
\phi_{>-10} \phi_{<10,\alpha}^\dagger \phi^{,\alpha}_{<10}\\
+& \phi_{\leq -10} \phi_{-10 < \cdot < 10,\alpha}^\dagger \phi_{-10 < \cdot < 10}^{,\alpha}\\
+& \phi_{\leq -10} \phi_{\leq -10,\alpha}^\dagger \phi_{-10 < \cdot < 10}^{,\alpha}\\
+& \phi_{\leq -10} \phi_{-10 < \cdot < 10,\alpha}^\dagger \phi_{\leq -10}^{,\alpha}\\
+& \phi_{\leq -10} \phi_{\leq -10,\alpha}^\dagger \phi_{\leq -10}^{,\alpha}. 
\end{eqnarray*}
We apply $P_0$, and discard some terms which are now zero, to obtain

\begin{eqnarray*}
P_0(\phi \phi_{,\alpha}^\dagger \phi^{,\alpha})= 
& \sum_{k_2,k_3: \max(k_2, k_3) \geq 10; |k_2-k_3| < 5} P_0(\phi
\phi_{k_2,\alpha}^\dagger \phi^{,\alpha}_{k_3}) \\
+&
\sum_{k_2,k_3: k_2 \geq 10, k_3 + 5}
P_0(\phi_{k_2 - 5 < \cdot < k_2 + 5} \phi_{k_2,\alpha}^\dagger \phi^{,\alpha}_{k_3})\\
+&
\sum_{k_2,k_3: k_3 \geq 10, k_2 + 5}
P_0(\phi_{k_3 - 5 < \cdot < k_3 + 5} \phi_{k_2,\alpha}^\dagger \phi^{,\alpha}_{k_3})
\\ +&
P_0(\phi_{>-10} \phi_{<10,\alpha}^\dagger \phi^{,\alpha}_{<10})\\
+& 
P_0(\phi_{\leq -10} \phi_{-10 < \cdot < 10,\alpha}^\dagger \phi_{-10 < \cdot < 10}^{,\alpha})\\
+& P_0(\phi_{\leq -10} \phi_{\leq -10,\alpha}^\dagger \phi_{-10 < \cdot < 10}^{,\alpha})\\
+& P_0(\phi_{\leq -10} \phi_{-10 < \cdot < 10,\alpha}^\dagger \phi_{\leq -10}^{,\alpha}).
\end{eqnarray*}
The first term is of type \eqref{dec1}, while the next three terms are of type \eqref{dec2} by \eqref{l-lp} and dyadic decomposition (swapping $k_2$ and $k_3$ as necessary).  Also, the fifth term is of type \eqref{dec1} by \eqref{l-lp}.  The sixth and seventh terms are equal, so it remains to show that
$$
P_0(\phi_{\leq -10} \phi_{\leq -10,\alpha}^\dagger \phi_{-10 < \cdot < 10}^{,\alpha}) + A_{\leq -10; \alpha} \phi_0^{,\alpha}$$
is of the right form.
From Lemma \ref{commutator} and the Leibnitz rule we see that
$$
P_0(\phi_{\leq -10} \phi_{\leq -10,\alpha}^\dagger \phi_{-10 < \cdot < 10}^{,\alpha})
= \phi_{\leq -10} \phi_{\leq -10,\alpha}^\dagger \phi_0^{,\alpha} + \eqref{fire1} + \eqref{fire2}.
$$
From the definition \eqref{a-def} of $A_{\leq -10; \alpha}$ it thus remains to show that
\be{sph1}
 \phi_{\leq -10,\alpha} \phi_{\leq -10}^\dagger \phi_0^{,\alpha}
\end{equation}
is of the right form.

Since $\phi$ lies on the sphere, we have the geometric identity
$$ \phi^\dagger \phi^{,\alpha} = 0$$
for all $\alpha$.  We apply $P_0$ and multiply by $\phi_{\leq -10,\alpha}$ to obtain
$$ \phi_{\leq -10,\alpha} P_0(\phi^\dagger \phi_{,\alpha}) = 0.$$
On the other hand, from Lemma \ref{commutator} we have
$$ \phi_{\leq -10,\alpha} P_0(\phi_{\leq -10}^\dagger \phi^{,\alpha}) = \phi_{\leq -10,\alpha} \phi_{\leq -10}^\dagger \phi_0^{,\alpha} + \phi_{\leq -10,\alpha} L(\nabla_x \phi_{\leq -10}^\dagger, \phi^{,\alpha}).$$
Subtracting these two identities and re-arranging we obtain
$$
\eqref{sph1} = \phi_{\leq -10,\alpha} L(\nabla_x \phi_{\leq -10}, \phi_0^{,\alpha})
+ \phi_{\leq -10,\alpha} P_0 L(\phi_{>-10}, \phi^{,\alpha}).
$$
The first term is of the form \eqref{fire1} by \eqref{l-lp}.  We can exploit the $P_0$ projection to split the second term as
$$
\phi_{\leq -10,\alpha} P_0 L(\phi_{-10 < \cdot < 10}, \phi^{,\alpha}_{<20}) + 
\sum_{k_2 \geq 20} \sum_{k_1: |k_1-k_2| < 5}
\phi_{\leq -10,\alpha} P_0 L(\phi_{k_1}, \phi^{,\alpha}_{k_2}).$$
By \eqref{l-lp} the first term is of the form \eqref{dec2} and the second term is of the form \eqref{dec1}.

\divider{Step 1(b).  Show that the small pieces are all acceptable errors.}

To finish the proof of \eqref{cancel} we need to show that the four expressions \eqref{dec1}-\eqref{fire2} are acceptable errors. 

\divider{Step 1(b).1.  The contribution of \eqref{dec1}. (High-high interaction)}

By Lemma \ref{minkowski} and \eqref{phi-bomb} it suffices to show
$$ \| \sum_{k_2 \geq O(1)} \sum_{k_3 = k_2 + O(1)} P_0 (\phi^{(1)} \phi^{(2)}_{k_2,\alpha} \phi^{(3),\alpha}_{k_3}) \|_{N_0} \lesssim c_0
\| \phi^{(1)}\|_{S(c)} \| \phi^{(2)}\|_{S(c)} \| \phi^{(3)}\|_{S(c)}$$
for all scalar $\phi^{(i)}$, $i=1,2,3$.  By the triangle inequality and dyadic decomposition we can estimate the left-hand side by
$$ \lesssim \sum_{k_2 \geq O(1)} \sum_{k_3 = k_2 + O(1)} \sum_k \| P_0 (\phi^{(1)} P_k(\phi^{(2)}_{k_2,\alpha} \phi^{(3),\alpha}_{k_3})) \|_{N_0}.$$
From \eqref{ur-algebra} we may bound this by
$$
\lesssim \| \phi^{(1)} \|_{S(c)} \sum_{k_2 \geq O(1)} \sum_{k_3 = k_2 + O(1)} \sum_k \chi^{(1)}_{k \leq 0} \| P_k (\phi^{(2)}_{k_2,\alpha} \phi^{(3),\alpha}_{k_3}) \|_{N_k}.$$
We may restrict $k$ to $k \leq k_2 + O(1)$, since the summand vanishes otherwise. From \eqref{null}, \eqref{s-def}, \eqref{local} we may majorize the previous by
$$\lesssim \| \phi^{(1)} \|_{S(c)} \| \phi^{(2)} \|_{S(c)} \| \phi^{(3)} \|_{S(c)} \sum_{k, k_2, k_3: k_2 \geq O(1), k_3 = k_2 + O(1); k \leq k_2 + O(1)}
\chi^{(1)}_{k \leq 0} \chi^{(1)}_{k_2 = k} c_{k_2} c_{k_3}.$$
We first sum in $k_3$, then in $k_2$, and finally in $k$, using \eqref{local} and splitting into $k\leq 0$ and $k>0$ if desired.  The sum then converges to $O(c_0^2)$, which is acceptable by \eqref{linfty-ass}.  

\divider{Step 1(b).2.  The contribution of \eqref{dec2}. (Derivative falls on low frequency)}

We can bound the $N_0$ norm of \eqref{dec2} by 
$$ \sum_{k_3} \sum_{k_2 < k_3+O(1)} \sum_{k_1 > k_2 + O(1)} \| P_0 L(\phi_{k_1}, \phi_{k_2,\alpha}, \phi^{,\alpha}_{k_3}) \|_{N_0}.$$
By \eqref{o-lemma} and \eqref{control-2} this is bounded by
$$ \lesssim \sum_{k_3} \sum_{k_2 < k_3+O(1)} \sum_{k_1 > k_2 + O(1)}
\chi^{(1)}_{0 = \max(k_1,k_3)}
\chi^{(1)}_{k_1 \leq k_2}  c_{k_1}  c_{k_2}  c_{k_3}.$$
Split into $k_3 \geq k_1$ and $k_3 < k_1$.  If $k_3 \geq k_1$, we perform the $k_1$ summation using \eqref{local} to estimate this by
$$ \lesssim \sum_{k_3} \sum_{k_2} 
\chi^{(1)}_{k_3 = 0}
c_{k_2}^2 c_{k_3}$$
which is acceptable by \eqref{l2-ass}, \eqref{local}.  If $k_3 < k_1$, we perform the $k_2$ and $k_3$ summations using \eqref{local} to estimate this by
$$ \lesssim \sum_{k_1}
\chi^{(1)}_{k_1=0} c_{k_1}^3$$
which is acceptable by \eqref{local}, \eqref{linfty-ass}.

\divider{Step 1(b).3.  The contribution of \eqref{fire1}. (Derivative transferred from high frequency to low, Type I)}

We estimate the $N_0$ norm using the triangle inequality and \eqref{l-lp-d-2} by
$$
\sum_{k \leq -10} 2^k \| L(\phi_k, \phi_{\leq -10,\alpha}, \phi_{-10 < \cdot < 10}^{,\alpha}) \|_{N_0}.$$
We may freely add a $P_{-10 < \cdot < 10}$ projection to the expression inside the norm.  By \eqref{nk-split} and Lemma \ref{minkowski} we may thus bound the previous by
$$
\sum_{k \leq O(1)} \sum_{k' = O(1)} \sum_{k_2 \leq O(1)} \sum_{k_3 = O(1)}
2^k \| P_{k'} L(\phi_k, \phi_{k_2,\alpha}, \phi^{,\alpha}_{k_3}) \|_{N_{k'}}.$$
By \eqref{o-lemma} and \eqref{control-2} we may bound this by 
$$
\lesssim \sum_{k \leq O(1)} \sum_{k' = O(1)} \sum_{k_2 \leq O(1)} \sum_{k_3 = O(1)} 
\chi^{(1)}_{k \leq k_2}
2^k  c_k  c_{k_2}  c_{k_3}.$$
Using \eqref{local}, \eqref{linfty-ass} we can simplify this  to
$$
 \eps^2 c_0 \sum_{k \leq O(1)} \sum_{k_2 \leq O(1)} \chi^{-(0)}_{k=0} 2^k 
\chi^{(1)}_{k \leq k_2}$$
which is acceptable.  

\divider{Step 1(b).4.  The contribution of \eqref{fire2}. (Derivative transferred from high frequency to low, Type II)}

The expression $\phi_{\leq -10}$ is bounded in $S(c)$ by \eqref{phi-bomb}, while the expression $L(\nabla_x \phi_{\leq -10,\alpha}, \phi_{-10 < \cdot < 10}^{,\alpha})$ has Fourier support in $D_0 \sim 1$.  By \eqref{nk-split} and \eqref{ur-algebra}, we can thus bound the $N_0$ norm of \eqref{fire2} by
$$ \lesssim 
\sum_{k = O(1)}
\| P_k L(\nabla_x \phi_{\leq -10,\alpha}, \phi_{-10 < \cdot < 10}^{,\alpha}) \|_{N_k}.$$
By the triangle inequality and \eqref{l-lp-d-2} we can thus bound the previous by
$$ \lesssim \sum_{k = O(1)} \sum_{k_1 \leq O(1)} \sum_{k_2 = O(1)} 2^{k_1}
\| P_k L(\phi_{k_1,\alpha}, \phi^{,\alpha}_{k_2}) \|_{N_k}.$$
By \eqref{null}, \eqref{control-2} this is bounded by
$$ \lesssim \sum_{k = O(1)} \sum_{k_1 \leq O(1)} \sum_{k_2 = O(1)} 2^{k_1}
c_{k_1}  c_{k_2}.$$
But this is acceptable by \eqref{local}, \eqref{linfty-ass}.  This completes the proof of \eqref{cancel}.

\divider{Step 2.  Construct a gauge $U_{\leq -10}$.}

We continue the proof of \eqref{psi-bound}.  In Step 3 we shall apply a renormalization $w_0 = U_{\leq -10} \phi_0$ which will transform \eqref{cancel} into a much better form, namely $\Box w_0 = error$.

In order to transform \eqref{cancel} like this, we would like\footnote{This scheme is slightly modified from that in \cite{tao:wavemap1}; basically, we have replaced $U_{\leq -10}$ by $U^\dagger_{\leq -10}$ as it allows for some technical simplifications.} $U_{\leq -10}$ to approximately be an orthogonal co-ordinate frame given by adjoint parallel transport by $A_{\leq -10; \alpha}$ in all directions; in other words, we expect $U_{\leq -10} U_{\leq -10}^\dagger \approx 1$ and $U_{\leq -10,\alpha} \approx -U_{\leq -10} A_{\leq -10; \alpha}$.  

We now give the construction of $U_{\leq -10}$.

Let $M > 10$ be a large integer (depending on $T$, $n$, $m$, and the $\delta_i$) to be chosen later.  We define the real $m \times m$ matrix fields $U_k$ for integer $k$ by setting $U_k = 0$ for $k \leq -M$ and
\be{ack-accurate}
U_k := U_{<k} (\phi_{<k} \phi_k^\dagger - \phi_k \phi_{<k}^\dagger) 
\end{equation}
for $k \geq M$, where
$$ U_{<k} := U_{\leq k-1} := I + \sum_{k' < k} U_{k'}$$
and $I$ is the identity matrix.

An easy inductive argument shows that $U_{<k}$ has Fourier support on the region $D_0 \leq 2^{k+5}$.  From a heuristic point of view we have $U_k \approx P_k U_{\leq -10}$ and $U_{<k} \approx P_{<k} U_{\leq -10}$.  In particular, we expect $U_k$ will obey similar estimates to $\phi_k$, and $U_{<k}$ will obey similar estimates to $\phi_{<k}$.

We now prove the estimates on $U_{<k}$ that we will need.

\divider{Step 2(a).  Show that the matrix $U_{<k}$ is approximately orthogonal.}

The purpose of this step is to prove the estimates
\begin{eqnarray}
\| U_{\leq k} \|_{S(Cc)} &\leq C
 \label{u-s}\\
\| U_k \|_{S_k} &\leq C  c_k \label{uk-s}\\
\| U_{\leq k} U_{\leq k}^\dagger - I \|_{S(c)} &\lesssim  \eps \label{uu-s}
\end{eqnarray}
for all $k$ and some absolute constant $C$ depending only on $n$, $m$, and the $\delta_i$, assuming that $M$ is sufficiently large (depending on $T$, $\eps$, $n$, $m$, and the $\delta_i$).

We prove \eqref{u-s}, \eqref{uk-s} by induction on $k$.  As a by-product of this induction argument we shall also obtain \eqref{uu-s}.

The estimates \eqref{u-s}, \eqref{uk-s} are trivial when $k \leq -M$.  Now suppose that $k > -M$ and that \eqref{u-s}, \eqref{uk-s} held for all smaller $k$.

From \eqref{ack-accurate} and two applications of \eqref{sk-skp} (one with $c$ and one with $Cc$) we have 
$$
\| U_k \|_{S_k} \lesssim \| U_{<k}\|_{S(Cc)} \| \phi_{<k} \|_{S(c)}
 \| \phi_k \|_{S_k},$$
and \eqref{uk-s} will thus follow from \eqref{phi-bomb}, \eqref{control-2} and the inductive hypothesis for \eqref{u-s}.

From \eqref{ack-accurate} we have
$$ U_{<k} U^\dagger_k + U_k U_{<k}^\dagger = 0.$$
Telescoping this we obtain
$$ U_{\leq k} U_{\leq k}^\dagger - I  = \sum_{k' < k} U_{k'} U_{k'}^\dagger.$$
By \eqref{s-sk} we have 
$$ \| U_{\leq k} U_{\leq k}^\dagger - I\|_{S(c)}  \lesssim \| \sum_{k' < k} U_{k'} U_{k'}^{\dagger} \|_{L^\infty_t L^\infty_x} +
\sup_{k' < k} \| U_{k'}\|_{S_k}^2/c_k.$$
By the triangle inequality, \eqref{infty-est}, and the induction hypothesis of \eqref{uk-s} we thus see that
$$ \| U_{\leq k} U_{\leq k}^\dagger - I\|_{S(c)}  \lesssim \sum_{k' < k} c_{k'}^2 + \sup_{k' < k} c_k,$$
which implies \eqref{uu-s} by \eqref{l2-ass}, \eqref{linfty-ass}.
From \eqref{infty-control} we then have have
$$ \| U_{\leq k} U_{\leq k}^\dagger - I \|_{L^\infty_t L^\infty_x} \lesssim  \eps$$
which implies that
$$ \| U_{\leq k} \|_{L^\infty_t L^\infty_x} \leq C/2.$$
From this and \eqref{uk-s} (which has just been proven for this value of $k$ and all preceding values) we obtain \eqref{u-s} if $C$ is sufficiently large.  This concludes the induction, and also gives the proof of \eqref{uu-s}.

\divider{Step 2(b).  Show that $(U_{\leq -M,\alpha} + U_{\leq -M} A_{\leq -M; \alpha}) \phi_0^{,\alpha}$ is an acceptable error.}

Since $U_{\leq M} = I$, it suffices by \eqref{l12}, \eqref{a-def} to show
$$ \| L(\phi_{\leq -M}, \phi_{\leq -M,\alpha}, \phi_0^{,\alpha}) \|_{L^1_t L^2_x} \lesssim  \eps c_0.$$
By H\"older we may majorize the right-hand side by
$$ T 
\| \phi_{\leq -M} \|_{L^\infty_t L^\infty_x} 
\| \nabla_{x,t} \phi_{\leq -M} \|_{L^\infty_t L^\infty_x}
\| \nabla_{x,t} \phi_0 \|_{L^\infty_t L^2_x}.$$
The first norm can be discarded since $\phi$ lies on the sphere.  
The second norm is $O(2^{-M})$ by \eqref{energy-est}, \eqref{control-2}, and \eqref{bernstein}. The third norm is $O(c_0)$ by \eqref{energy-est}, \eqref{control-2}.  The claim then follows if $M$ is sufficiently large depending on $T$.

From \eqref{u-s}, \eqref{uk-s}, and \eqref{insensitive} we have
\begin{eqnarray}
\| U_{\leq k} \|_{S(c)} &\lesssim 1
 \label{u-s-2}\\
\| U_k \|_{S_k} &\lesssim c_k. \label{uk-s-2}
\end{eqnarray}

\divider{Step 2(c).  Show that $(U_{\leq -10,\alpha} + U_{\leq -10} A_{\leq -10; \alpha}) \phi_0^{,\alpha}$ is an acceptable error.}

From Step 2(b) and the triangle inequality it suffices by \eqref{l2-ass} to show the telescoping bound
$$
\|
(U_{k,\alpha} + U_{\leq k} A_{\leq k;\alpha} - U_{< k} A_{<k;\alpha}) \phi_0^{,\alpha} \|_{N_0} \lesssim  c_k^2 c_0
$$
for all $-M < k \leq -10$. 

We use \eqref{ack-accurate}, \eqref{a-def} to expand
\begin{eqnarray*}
U_{k,\alpha} =&
U_{<k} (\phi_{<k} \phi_{k,\alpha}^\dagger - \phi_{k,\alpha} \phi_{<k}^\dagger) \\
+& U_{<k,\alpha} (\phi_{<k} \phi_k^\dagger - \phi_k \phi_{<k}^\dagger) \\
+& U_{<k} (\phi_{<k,\alpha} \phi_k^\dagger - \phi_k \phi_{<k,\alpha}^\dagger) 
\end{eqnarray*}
and
\begin{eqnarray*}
U_{\leq k} A_{\leq k;\alpha} - U_{< k} A_{<k;\alpha} =&
U_{<k} (\phi_{k,\alpha} \phi_{<k}^\dagger -
\phi_{<k} \phi_{k,\alpha}^\dagger)\\
+&
U_{<k} (\phi_{\leq k,\alpha} \phi_k^\dagger -
\phi_k \phi_{\leq k,\alpha}^\dagger)\\
+&
U_k (\phi_{\leq k,\alpha} \phi_{\leq k}^\dagger -
\phi_{\leq k} \phi_{\leq k,\alpha}^\dagger).
\end{eqnarray*}
Comparing the two, we see that the first terms of both expressions cancel.  By several applications of \eqref{l-lp} we thus have
\begin{eqnarray*}
(U_{k,\alpha} + U_{\leq k} A_{\leq k;\alpha} - U_{< k} A_{<k;\alpha}) \phi_0^{,\alpha}
= &L(\phi_{\leq k}, \phi_k, U_{<k,\alpha}, \phi_0^{,\alpha})\\
+& L(U_{<k}, \phi_k, \phi_{\leq k,\alpha}, \phi_0^{,\alpha})\\
+& L(\phi_{\leq k}, U_k, \phi_{\leq k,\alpha}, \phi_0^{,\alpha}).
\end{eqnarray*}
Thus it remains to show that all the terms on the right-hand side
have a $N_0$ norm of $O( c_k^2 c_0)$. Note that all these non-linearities are of the good ``derivative falls on a low frequency'' type, as the derivative fails to fall on $\phi_k$ or $U_k$.  

The expressions $\phi_{\leq k}$, $U_{<k}$ are bounded in $S(c)$ by \eqref{phi-bomb}, \eqref{uk-s-2}.  By \eqref{ur-algebra}, \eqref{local} and dyadic decomposition it thus suffices to show that
$$\| L(\phi_k, U_{k',\alpha}, \phi_0^{,\alpha}) \|_{N_0}
+ \| L(\phi_k, \phi_{k',\alpha}, \phi_0^{,\alpha}) \|_{N_0}
+ \| L(U_k, \phi_{k',\alpha}, \phi_0^{,\alpha}) \|_{N_0}
\lesssim \chi^{(1)}_{k \leq k'}  c_k  c_{k'}  c_0$$
for all $k' \leq k$.  But this follows from \eqref{o-lemma}, \eqref{uk-s-2}, and \eqref{control-2}.

\divider{Step 2(d).  Show that $\Box U_{\leq -10} \phi_0$ is an acceptable error.}

Morally speaking, this is a non-linearity of ``repeated derivatives avoiding the high frequency'' type, although the treatment gets somewhat technical because the recursive construction of $U$.

Applying $\Box$ to the definition \eqref{ack-accurate} of $U_k$ we obtain the identity
\be{u-ident}
\Box U_k = \Box U_{<k} \Phi_k + F_k = \sum_{-M < k' < k} \Box U_{k'} \Phi_k + F_k
\end{equation}
for all $-M < k \leq -10$, where $\Phi_k$ is the matrix
$$ \Phi_k := \phi_{<k} \phi_k^\dagger - \phi_k \phi_{<k}^\dagger$$
and $F_k$ is the matrix
\begin{eqnarray*}
F_k := &U_{<k} (\Box \phi_{<k} \phi_k^\dagger - \phi_k \Box \phi_{<k}^\dagger) \\
+& U_{<k} (\phi_{<k} \Box \phi_k^\dagger - \Box \phi_k \phi_{<k}^\dagger) \\
+& 2 U_{<k,\alpha}
(\phi_{<k}^\alpha \phi_k^\dagger - \phi_k \phi_{<k}^{,\alpha\dagger}) \\
+& 2 U_{<k,\alpha}
(\phi_{<k} \phi_k^{,\alpha\dagger} - \phi_k^{,\alpha} \phi_{<k}^\dagger) \\
+& 2 U_{<k}
(\phi_{<k,\alpha} \phi_k^{,\alpha\dagger} - \phi_{k,\alpha} \phi_{<k}^{,\alpha\dagger}).
\end{eqnarray*}
We can simplify $F_k$ using \eqref{l-lp} as
\be{fk-form}
F_k = 
L(U_{<k}, \phi_{\leq k}, \Box \phi_{\leq k}) + 
L(\phi_{\leq k}, U_{<k,\alpha},  \phi_{\leq k}^{,\alpha})
+
L(U_{<k}, \phi_{\leq k,\alpha},  \phi_{\leq k}^{,\alpha}).
\end{equation}
Iterating \eqref{u-ident} completely we see that
$$
\Box U_{\leq -10} = \sum_{s=1}^{M+10} \sum_{-M < k_1 < \ldots < k_s \leq -10} F_{k_1} \Phi_{k_2} \ldots \Phi_{k_s}.
$$
At first glance this series seems to grow exponentially in $s$, but we shall eventually extract a $1/s!$ decay\footnote{The $1/s!$ decay is a benefit of the recursive construction of $U$, in that low frequencies affect high frequencies but not vice versa. It may be possible to obtain the desired estimate on $(\Box U_{\leq -10}) \phi_0$ without relying on the $1/s!$ decay; for instance, the $L^4_t L^\infty_x$ Strichartz estimate should be able to deal with the cases $s \geq 5$ by the techniques in \cite{tao:wavemap1}.  On the other hand, a similar decay is very useful in \cite{kl-mac:algebra}, and in a less related (but still critical regularity) context in \cite{christ:kiselev3}.} from the ordering of the $k_j$ to counteract this.

We multiply both sides by $\phi_0$ on the right and take the $N_0$ norm.  Using the $L$ notation (possibly conceding some powers of $C m$ because of matrix multiplication) we obtain
$$
\| (\Box U_{\leq -10}) \phi_0\|_{N_0}
\lesssim
\sum_{s=1}^{M+10} \sum_{-M < k_1 < \ldots < k_s \leq -10} 
(Cm)^{Cs} \| L(F_{k_1}, \Phi_{k_2}, \ldots, \Phi_{k_s}, \phi_0) \|_{N_0}.$$
We can re-arrange the right-hand side as
$$
\sum_{-M < k_1 \leq -10}
\sum_{s=1}^{M+10} \sum_{k_1 < k_2 < \ldots < k_s \leq -10} 
(Cm)^{Cs} \| L(F_{k_1}, \phi_0, \Phi_{k_2}, \ldots, \Phi_{k_s}) \|_{N_0}.$$
From \eqref{sk-skp}, \eqref{phi-bomb}, \eqref{control-2} we have
$$ \| \Phi_k \|_{S_k} \lesssim \| \phi_k \|_{S_k} \| \phi_{<k} \|_{S(c)}
\lesssim  c_k.$$
By repeated application of \eqref{ur-algebra-sk} we can thus estimate the previous by
$$
\lesssim \sum_{-M < k_1 \leq -10}
\sum_{s=1}^{M+10} \sum_{k_1 < k_2 < \ldots < k_s \leq -10} 
(C m)^{Cs} \| L(F_{k_1}, \phi_0) \|_{N_0}
c_{k_2} \ldots c_{k_s}.$$
By the binomial (or multinomial) theorem we can estimate this by
$$
\lesssim \sum_{-M < k_1 \leq -10}
\sum_{s=1}^{M+10} 
(C  m)^{Cs} \| L(F_{k_1}, \phi_0) \|_{N_0}
\frac{1}{s!} 
(\sum_{k_1 < k \leq -10} c_k)^s,$$
which in turn is estimated by
$$
\lesssim \sum_{-M < k_1 \leq -10}
\| L(F_{k_1}, \phi_0) \|_{N_0} \exp(C  m \sum_{k_1 < k \leq -10} c_k).$$
From \eqref{l2-ass} and Cauchy-Schwarz we have
$$ \sum_{k_1 < k \leq -10} c_k \lesssim \eps |k_1|^{1/2},$$
so that the previous is majorized by
$$
\lesssim \sum_{-M < k_1 \leq -10}
\| L(F_{k_1}, \phi_0) \|_{N_0} \exp(C  m \eps |k_1|^{1/2}).$$
It will thus suffice to show the bound
$$ \| L(F_k, \phi_0) \|_{N_0} \lesssim \chi^{(1)}_{k=0} \eps c_0$$
for all $-M < k \leq -10$.  

From \eqref{fk-form} we have
\begin{eqnarray*}
L(F_k,\phi_0) = &L(U_{<k}, \phi_{\leq k}, \Box \phi_{\leq k}, \phi_0) \\
+& L(\phi_{\leq k}, \phi_0, U_{<k,\alpha}, \phi_{\leq k}^{,\alpha}) \\
+& L(U_{<k}, \phi_0, \phi_{\leq k,\alpha}, \phi_{\leq k}^{,\alpha}) 
\end{eqnarray*}
so it suffices to show that these three expressions have an $N_0$ norm of $\lesssim \chi^{(1)}_{k=0}  \eps c_0$.

The quantities $U_{<k}$, $\phi_{\leq k}$ are all bounded in $S(c)$ by \eqref{u-s-2}, \eqref{phi-bomb}.  So by \eqref{ur-algebra} it suffices to show that
\begin{eqnarray}
\| &L(\Box \phi_{\leq k}, \phi_0) \|_{N_0} \label{box2}\\
\| &L(\phi_0, U_{<k,\alpha}, \phi_{\leq k}^{,\alpha}) \|_{N_0} \label{q1}\\
\| &L(\phi_0, \phi_{\leq k,\alpha}, \phi_{\leq k}^{,\alpha}) \|_{N_0} \label{q3}
\end{eqnarray}
are $\lesssim \chi^{(1)}_{k=0}  \eps c_0$.

\divider{Case 2(d).1. The treatment of \eqref{q1} (Derivative falls on low frequency).}

By dyadic decomposition and \eqref{o-lemma} we can bound \eqref{q1} by
$$ \sum_{k', k'' \leq k} \chi^{(1)}_{0 \leq \min(k',k'')}
\| \phi_0 \|_{S_0} \| U_{k'} \|_{S_{k'}} \| \phi_{k''} \|_{S_{k''}}.$$
But this is acceptable by \eqref{control-2}, \eqref{uk-s-2}, \eqref{control-2}.  

\divider{Case 2(d).2. The treatment of \eqref{q3} (Derivative falls on low frequency).}

This is identical to Case 2(d).1 except that \eqref{uk-s-2} is replaced with \eqref{control-2}.

\divider{Case 2(d).3. The treatment of \eqref{box2}.}

By dyadic decomposition it suffices to show
$$
\| L(\Box \phi_k, \phi_0) \|_{N_0} \lesssim \chi^{(1)}_{k=0} \eps c_0.$$
We may freely insert $P_{-10<\cdot<10}$ in front of the $L$.  Expanding $\Box \phi_k$ using the equation \eqref{wavemap-eq}, we can write the left-hand side as
$$ \| P_{-10 < \cdot < 10} L( P_k L(\phi, \phi_{,\alpha}, \phi^{,\alpha}), \phi_0)\|_{N_0}.$$
By dyadic decomposition and symmetry it suffices to show that
\be{dyad}
\sum_{k_2,k_3: k_2 \leq k_3} \| P_{-10 < \cdot < 10} L( P_k L(\phi, \phi_{k_2,\alpha}, \phi_{k_3}^{,\alpha}), \phi_0)\|_{N_0}
\lesssim \chi^{(1)}_{k=0} \eps c_0.
\end{equation}

We divide \eqref{dyad} into three contributions and treat each separately.

\divider{Case 2(d).3(a).  The contribution to \eqref{dyad} when $k_2 \leq O(1)$.  (Derivative falls on low frequency)}

We use \eqref{l-lp} to write this contribution as
$$ \sum_{k_2,k_3: k_2 \leq k_3, O(1)} \| P_{-10 < \cdot < 10} L( \phi, \phi_0, \phi_{k_2,\alpha}, \phi_{k_3}^{,\alpha})\|_{N_0}.$$
We split the first $\phi$ as $\phi_{<10} + \sum_{k_1 \geq 10} \phi_{k_1}$.  To deal with $\phi_{<10}$ we use \eqref{o-lemma}, \eqref{compat} to estimate this contribution  by
$$ 
 \sum_{k_2,k_3: k_2 \leq k_3, O(1)} \chi^{(1)}_{0 = \max(0, k_3)}
\chi^{(1)}_{0 \leq k_2}
\| L(\phi_{<10}, \phi_0) \|_{S_0} \| \phi_{k_2} \|_{S_{k_2}} \| \phi_{k_3} \|_{S_{k_3}}.$$
Using \eqref{control-2}, \eqref{phi-bomb}, and \eqref{sk-skp} we can bound this by
$$ 
\sum_{k_2,k_3: k_2 \leq k_3, O(1)} 
\chi^{(1)}_{0 = \max(0, k_3)}
\chi^{(1)}_{0 \leq k_2}
c_0 c_{k_2} c_{k_3}$$
which is acceptable by \eqref{linfty-ass} (performing the $k_3$ summation first).  To deal with
$\sum_{k_1 \geq 10} \phi_{k_1}$, we again use \eqref{o-lemma} to estimate this contribution by
$$ 
\sum_{k_1 \geq O(1)} \sum_{k_2,k_3: k_2 \leq k_3, O(1)} 
\chi^{(1)}_{0 = \max(k_1,k_3)}
\chi^{(1)}_{k_1 \leq k_2}
\| L(\phi_{k_1}, \phi_0) \|_{S_{k_1}} \| \phi_{k_2} \|_{S_{k_2}} \| \phi_{k_3} \|_{S_{k_3}}$$
which by the previous arguments is bounded by
$$ 
\sum_{k_1 \geq O(1)} \sum_{k_2,k_3: k_2 \leq k_3, O(1)} 
\chi^{(1)}_{0 = \max(k_1,k_3)}
\chi^{(1)}_{k_1 \leq k_2}
c_{k_1} c_{k_2} c_{k_3}.$$
But this is acceptable by \eqref{local} and \eqref{linfty-ass} (performing the $k_1$ summation, then the $k_3$, then the $k_2$).

\divider{Case 2(d).3(b).  The contribution to \eqref{dyad} when $k_2 \geq O(1)$ and $k_3 = k_2 + O(1)$. (High-high interaction)}

By \eqref{ur-algebra} and \eqref{control-2} we may control this contribution by
$$
\lesssim \sum_{k_2 \geq O(1)} \sum_{k_3 = k_2 + O(1)} c_0 \| P_k L(\phi, \phi_{k_2,\alpha}, \phi_{k_3}^{,\alpha}) \|_{N_k}.
$$
Since $\phi$ is bounded in $S(c)$ by \eqref{phi-bomb}, we apply another dyadic decomposition and \eqref{ur-algebra} again to bound this by
$$ \lesssim \sum_{k_2 \geq O(1)} \sum_{k_3 = k_2 + O(1)} \sum_{k_0} c_0 
\chi^{(1)}_{k \geq k_0} \| P_{k_0} L(\phi_{k_2,\alpha}, \phi_{k_3}^{,\alpha}) \|_{N_{k_0}}.$$
By \eqref{null} and \eqref{control-2} this is bounded by
$$ \lesssim  \sum_{k_2 \geq O(1)} \sum_{k_3 = k_2 + O(1)} \sum_{k_0} 
c_0 \chi^{(1)}_{k \geq k_0}
\chi^{(1)}_{k_0 = \max(k_2,k_3)} c_{k_2} c_{k_3},$$
which is acceptable by \eqref{linfty-ass}.

\divider{Case 2(d).3(c).  The contribution to \eqref{dyad} when $k_2 \geq O(1)$ and $k_3 < k_2 - 5$. (Derivative falls on low frequency)}

As in Case 2(d).3(b) we can use \eqref{ur-algebra}, \eqref{control-2} to control this contribution by 
$$
\lesssim \sum_{k_2 \geq O(1)} \sum_{k_3 < k_2-5} c_0 \| P_k L(\phi, \phi_{k_2,\alpha}, \phi_{k_3}^{,\alpha}) \|_{N_k}.
$$
By dyadic decomposition, \eqref{o-lemma} and \eqref{control-2} we can thus control the previous by
$$
\lesssim 
\sum_{k_2 \geq O(1)} \sum_{k_3 \leq k_2 + O(1)} \sum_{k_1}
c_0 
\chi^{(1)}_{k = \max(k_1,k_2)} \chi^{(1)}_{k_1 \leq k_3} c_{k_1} c_{k_2} c_{k_3}$$
which is acceptable by \eqref{linfty-ass}.  This concludes Step 2(c).

\divider{Step 3.  Renormalization}

We can now complete the proof of \eqref{psi-bound}. Define the renormalization $w_0$ of $\phi_0$ by $w_0 := U_{\leq -10} \phi_0$. 

\divider{Step 3(a).  Estimate $\Box w_0$.}

From the Leibnitz rule and \eqref{cancel} we have
\be{boxw}
\Box w_0 = 2 (U_{\leq -10, \alpha} + U_{\leq -10} A_{\leq -10; \alpha}) \phi_0^{,\alpha} 
+ \Box U_{\leq -10} \phi_0 + U_{\leq -10} error.
\end{equation}
Note that $U_{\leq -10} error$ must have Fourier support in the region $D_0 \sim 1$, since all the other terms of \eqref{boxw} do.  From \eqref{u-s-2}, \eqref{ur-algebra}, and \eqref{nk-split} we thus see that this term is an acceptable error.  From Steps 2(c), 2(d) the other two terms on the right-hand side of \eqref{boxw} are also acceptable errors.
We thus have
$$ \| \Box w_0 \|_{N_0} \lesssim \eps c_0.$$

\divider{Step 3(b).  Estimate $w_0[0]$.}

From \eqref{u-s-2}, \eqref{infty-control} we have
$$ \| U_{\leq -10}(0) \|_{L^\infty} \lesssim 1.$$
Also, from \eqref{u-s-2}, \eqref{energy-est}, and Bernstein's inequality \eqref{bernstein} we have
$$ \| P_k \partial_t U_{\leq -10}(0) \|_{L^\infty} \lesssim 2^k$$
for all $k$.  Summing over all $k \leq C$ and using the Fourier support of $U_{\leq -10}$ we obtain
$$ \| U_{\leq -10}[0] \|_{L^\infty \times L^\infty} \lesssim 1.$$
From the hypothesis that $\phi$ lies underneath the frequency envelope $\eps c$ we have
\be{psi-bang}
\| \phi_0[0] \|_{L^2 \times L^2} \sim \| \phi_0[0] \|_{\dot H^{n/2} \times \dot H^{n/2-1}} \lesssim \eps c_0.
\end{equation}
Combining the two we see that
$$ \| w_0[0] \|_{\dot H^{n/2} \times \dot H^{n/2-1}} \sim \| w_0[0] \|_{L^2 \times L^2} \lesssim \eps c_0.$$

\divider{Step 3(c).  Invert $U_{\leq 10}$.}

From the Fourier support of $U$ and $\phi$ we see that $w$ has Fourier support in the region $2^{-5} \leq D_0 \leq 2^5$.  From Steps 3(a), 3(b) and the energy estimate \eqref{energy-est-2} we thus have
\be{w-bound}
\| w_0 \|_{S_0} \lesssim \eps c_0.
\end{equation}

To obtain \eqref{psi-bound} from \eqref{w-bound} we need to invert $U_{\leq -10}$.  From \eqref{algebra} we know that $S(c)$ is a Banach algebra.  From \eqref{uu-s}  we thus see that $U_{\leq -10} U_{\leq -10}^\dagger$ is invertible in $S(c)$, so that 
$$ \| (U_{\leq -10}^{-1})^\dagger U_{\leq -10}^{-1} \|_{S(c)} \lesssim 1.$$
From \eqref{u-s-2} and Lemma \ref{algebra} we thus have
$$ \| U_{\leq -10}^{-1} \|_{S(c)} \lesssim 1.$$
To use this we observe from the Fourier support of $\phi_0$, $w_0$ that
$$ \phi_0 =P_{-5 < \cdot < 5} \phi_0 = P_{-5 < \cdot < 5} (U_{\leq -10}^{-1} w_0) = P_{-5 < \cdot 5} ( P_{< 10} U_{\leq -10}^{-1} w_0).$$
From \eqref{sk-skp} and the preceding bounds, we thus have
$$ \| \phi_0 \|_{S_0} 
\lesssim \| P_{< 10} U_{\leq -10}^{-1} w_0 \|_{S_0}
\lesssim \| P_{< 10} U_{\leq -10}^{-1} \|_{S(c)} \| w_0\|_{S_0}
\lesssim \| U_{\leq -10}^{-1} \|_{S(c)} \|w_0 \|_{S_0} \lesssim \eps c_0$$
and \eqref{psi-bound} follows.  This completes the proof of Proposition \ref{reduced} and hence Theorem \ref{main2}.
\endprf

We have now completed the iterative portion of the proof of Theorem \ref{main2}; it remains only to prove Theorem \ref{spaces}.  This shall occupy the remainder of the paper.

\section{Overview of spaces and estimates}\label{overview-sec}

We now begin the proof of Theorem \ref{spaces}.  Our spaces $S(c)$, $S_k$, $N_k$ are modeled on
 those of Tataru in \cite{tataru:wave2}.  Roughly speaking, these spaces can be divided into a ``Fourier'' or ``$\dot X^{s,b}$'' component, and a ``physical space'' or ``$L^q_t L^r_x$'' component, although the separation is not entirely clean, and the physical space component utilizes null frames in addition to the usual Euclidean space-time frame.

This section will only be an outline of the spaces and estimates to come, and will therefore be rather informal in nature.  Much of the notation used in this section will be defined later on.

In the sub-critical theory $s>n/2$ one usually chooses $N_k$ to be $X^{s-1,-1/2+}_k$ and $S(c)$ to be $\nabla_{x,t} X^{s-1,1/2+}$, in which case one essentially has Theorem \ref{spaces} for short times.  Thus at the critical level $s=n/2$ one expects the space $N_k$ to be something like $\dot X^{n/2-1,-1/2}_k$, and $S(c)$ to be something like $\nabla_{x,t}^{-1} \dot X^{n/2-1,1/2}$.  However, these spaces contain too many logarithmic divergences to be useful (for instance, $\nabla_{x,t}^{-1} \dot X^{n/2-1,1/2}$ will just barely fail to control $L^\infty_t \dot H^{n/2}_x$).  

In the high dimensional case $n \geq 4$ Tataru \cite{tataru:wave1} worked with the Besov counterparts $\dot X^{n/2-1,-1/2,1}_k$, $\dot X^{n/2,1/2,1}_k$ of the above spaces, and showed that if these spaces were combined with Strichartz spaces such as $L^\infty_t \dot B^{n/2,1}_2$ (for $S(c)$) and $L^1_t \dot B^{n/2-1,1}_2$ (for $N$), then one could recover all the estimates above provided that one localized each function to a single dyadic block.  This is enough to obtain a regularity and well-posedness theory for the Besov space $\dot B^{n/2,1}_2$.  When $n=2,3$
the Strichartz spaces are not strong enough to close all the above estimates, but in \cite{tataru:wave2} null frame Strichartz spaces (specifically, $L^2_{t_\omega} L^\infty_{x_\omega}$ and $L^\infty_{t_\omega} L^2_{x_\omega}$ type spaces for $S(c)$, and $L^1_{t_\omega} L^2_{x_\omega}$ type spaces for $N$) were introduced in order to resolve this problem.

Our spaces shall be similar to those in \cite{tataru:wave2}, but we shall need to prove stronger estimates than what was shown in that paper, and so we have modified the spaces somewhat in order to accomplish this.  Specifically, we need to obtain an exponential gain (already implicit in \cite{tataru:wave2}) in all high-high interactions, and we also need an exponential gain in the trilinear estimate $S S_{,\alpha} S^{,\alpha} \subset N$ when one or more derivatives fall on low frequency terms (which is a substantially more difficult gain to obtain, especially when $n=2$).  Both of these gains are essential in this paper because these are the two types of interactions which are not removed or mollified by the renormalization.  Finally, we need to select $S(c)$ so that it is a genuine algebra, as opposed to merely being obeying product-type estimates when all functions are restricted to dyadic blocks; this is essential in order for us to invert the renormalization, among other things.  Heuristically the algebra property should be obtainable just by adding $L^\infty_t L^\infty_x$ control to the space $S(c)$, but the actual proof requires some care, mainly because one needs to commute multiplication by $L^\infty_t L^\infty_x$ functions with frequency projections such as $P_k$ or $Q_j$ at various junctures.

There is a certain amount of discretion in how to select the $S$ and $N$ spaces.  However, any attempt to make one of the above estimates easier to prove (for instance, by making a space $S$ or $N$ stronger or weaker) often causes two or three other estimates to become more difficult to prove.  The spaces that the author eventually settled upon were obtained via a lengthy trial-and-error procedure, in which extreme difficulties in one of the above estimates were traded (via one or more modifications of the $S$ and $N$ spaces) for slightly less extreme difficulties in other estimates, with this process being iterated until all estimates could eventually be proved by the author.  As in other work at or near the critical problem, a partial duality relationship \eqref{sf-duality} between $S$ and $N$ can be exploited to slightly reduce the amount of computation required.

The rest of the paper is organized as follows.  In Section \ref{xsb-def-sec} we 
define the $Q_j$ family of multipliers which is used to define spaces of $\dot X^{s,b}$ type, and prove some technical lemmata allowing us to manipulate the $Q_j$, as well as a Strichartz estimate for these spaces. In 
Section \ref{frames-sec} we set out the notation for spherical caps and null frames, define two key spaces $NFA[\kappa]$, $S[k,\kappa]$ adapted to these null frames and caps, and set out their properties.   With these preliminaries, and some sector projections defined in Section \ref{sector-sec}, we shall be able to define the spaces $S_k$, $N_k$, $S(c)$ in Section \ref{construction-sec}, and set out some basic estimates for these spaces, including a partial duality between the $S$ and $N$ spaces.  We then prove the energy estimate \eqref{energy-est-2} in Section \ref{nk-sec}, which then allows us to demonstrate the quasi-continuity estimate \eqref{continuous} in Section \ref{continuity-sec}.

The remaining sections are devoted to bilinear and trilinear estimates.  We begin by a discussion on the geometry of the cone in Section \ref{geometry-sec}.  We then prove a basic product estimate, Lemma \ref{core}, in Section \ref{prod-sec}.  Using this Lemma and some additional effort, we shall be able to prove the product estimates \eqref{ur-algebra}, \eqref{ur-algebra-sk} in Section \ref{ur-sec}, the algebra estimates \eqref{algebra}, \eqref{sk-skp}, \eqref{sk-sk} in Section \ref{algebra-sec}, and the null form estimate \eqref{null} in Section \ref{null-sec}.  Finally in Section \ref{o-lemma-sec} we use the estimates just proven, plus some additional arguments, to prove the trilinear estimate \eqref{o-lemma}.

\section{$\dot X^{s,b}$ type spaces}\label{xsb-def-sec}

We first use a convenient trick (going back at least to \cite{ginibre:survey}) to work in the global Minkowski space $\R^{1+n}$ instead of just the time-localized slab $[-T,T] \times \R^n$.  Namely, if $X(\R^{1+n})$ is a Banach space of functions on $\R^{1+n}$ and $T \geq 0$, we define $X([-T,T] \times \R^n)$ to be the restriction of functions in $X$ to $[-T,T] \times \R^n$, with the norm
$$ \| \phi \|_{X([-T,T] \times \R^n)} := \inf \{ \| \Phi\|_{X(\R^{1+n})}: \Phi|_{[-T,T] \times \R^n} = \phi \}.$$
We shall construct the above Banach spaces first on $\R^{1+n}$, prove the estimates in Theorem \ref{spaces} in this global setting, and then restrict these spaces to $[-T,T] \times \R^n$ by the above procedure.  This restriction will not affect the above estimates (except for some minor technicalities involving time cutoffs in proving \eqref{energy-est}, which we shall address in Section \ref{nk-sec}).  For a further discussion see e.g. \cite{selberg:thesis}. 

Now that we are working globally in spacetime, we now have available the space-time Fourier transform
$$ {\cal F} \phi(\tau,\xi) := \int\int e^{-2\pi i (x \cdot \xi + t \tau)} \phi(t,x)\ dt dx.$$

It will be convenient to define the non-negative quantities $D_0$, $D_+$, $D_-$ on frequency space $\{(\tau,\xi): \tau \in \R, \xi \in \R^n \}$ by
$$ D_0 := |\xi|; \quad D_+ := |\xi| + |\tau|; \quad D_- := ||\tau| - |\xi||.$$
Note that $D_+ \sim D_0 + D_-$.  We shall also define the direction $\Theta \in S^{n-1}$ by
$$ \Theta := \tau \xi / |\tau\xi|$$
if $\tau \xi \neq 0$, and $\Theta := 0$ otherwise.  We sometimes refer to $D_0$ as the \emph{frequency}, $D_-$ as the \emph{modulation}, and $\Theta$ as the \emph{direction}.  We shall often use $D_0$, $D_-$, $D_+$, and $\Theta$ to define various regions in frequency space, thus for instance $\{ D_0 \sim 2^k, D_- \sim 2^j \}$ denotes the frequency region $\{ (\tau,\xi): |\xi| \sim 2^k, ||\xi|-|\tau|| \sim 2^j \}$.

Most of the action shall take place in the region $D_0 \sim D_+$, and one can think of these quantities as heuristically equivalent.  In practice, the region $D_+ \gg D_0$ will generate numerous minor sub-cases in the estimates, which can usually be disposed of quickly.

Using the spacetime Fourier transform we can define Littlewood-Paley projections adapted to the light cone.  For any integer $j$, define the projection operator $Q_{\leq j} = Q_{<j+1}$ by
$$ {\cal F} (Q_{\leq j} \phi)(\tau,\xi) := m_0(2^{-j}D_-) {\cal F} \phi(\tau,\xi).$$
Similarly define $Q_j$, $Q_{\geq j}$, $Q_{j_1 \leq \cdot \leq j_2}$, etc. in exact analogy to the corresponding $P$ multipliers.  Thus, for instance, $Q_j$ has symbol $m(2^{-j}D_-)$.

We observe the decomposition\footnote{Of course, this identity is not true for global free solutions, since $Q_j \phi$ always vanishes then, but all our functions shall be time localized, so this issue does not arise.}
$\phi = \sum_j Q_j \phi$
for Schwartz functions $\phi$.  

For any sign $\pm$, define the operator $Q^\pm_{\leq j}$ to be the restriction of $Q_{\leq j}$ to the frequency region $\pm \tau \geq 0$, thus $Q_{\leq j} = Q^+_{\leq j} + Q^-_{\leq j}$ and
$$ {\cal F} (Q^\pm_{\leq j} \phi)(\tau,\xi) := \chi_{[0,\infty)}(\pm \tau) m_0(2^{-j}(\pm \tau - |\xi|)) {\cal F}  \phi(\tau,\xi).$$
Similarly define $Q^\pm_j$, $Q^\pm_{\geq j}$, $Q^\pm_{j_1 \leq \cdot \leq j_2}$, etc.  The reader should be cautioned that $Q^\pm_{\leq j} + Q^\pm_{> j}$ does not add up to the identity, but rather to the Riesz projection to the half-space $\{ \pm \tau \geq 0 \}$.

Note that the operators $P_k$, $Q_j$, $Q^\pm_j$, $Q_{\leq j}$, etc. are all spacetime Fourier multipliers and thus commute with each other and with constant co-efficient differential operators.  The operator $P_k Q_j$ is essentially a projection onto the functions with spatial frequency $2^k$ and modulation $2^j$.  

\begin{definition}\label{disposable-def} A spacetime Fourier multiplier is said to be \emph{disposable} if its (distributional) convolution kernel is given by a measure with total mass $O(1)$.
\end{definition}

By Minkowski's inequality, a disposable multiplier is bounded on all spacetime-translation-invariant Banach spaces; in practice, this means that these operators can be discarded whenever one wishes.  Note also that the composition of two disposable multipliers is also disposable. 

The Littlewood-Paley operators $P_k$, $P_{\leq k}$, $P_{\geq k}$, etc. can easily be seen to be disposable.  Unfortunately, the operators $Q_j$, $Q^\pm_j$, $Q_{\leq j}$, etc. are not disposable, however the following two lemmas will serve as adequate (and very useful!) substitutes.

\begin{lemma}\label{disposable}  If $j, k$ are integers such that $j \geq k+O(1)$, then $P_{\leq k} Q_{\leq j}$, $P_k Q_{\leq j}$, $P_{\leq k} Q_j$, and $P_k Q_{\leq j}$ are disposable multipliers.

If in addition we assume $|k-j|>10$, then $P_k Q^\pm_j$, $P_k Q^\pm_{\leq j}$, $P_k Q^\pm_{\geq j}$ are disposable for either choice of sign $\pm$.
\end{lemma}

\begin{proof}
Observe that when $|k-j|>10$, then $P_k Q^\pm_{\leq j}$ can be factorized as the product of $P_k Q_{\leq j}$ and a multiplier whose symbol is a bump function adapted to $D_+ \sim 2^k$, and is therefore disposable.  Thus it suffices to prove the claims in the first paragraph.  In fact, it suffices to verify the claim for $P_{\leq k} Q_{\leq j}$, as the other projections are linear combinations of this type of multiplier.

We may rescale $k=0$.  The symbol $m$ of $P_{\leq 0} Q_{\leq j}$ supported on the region $|\xi| \lesssim 1, |\tau| \lesssim 2^j$, is smooth except when $\xi = 0$ or $\tau = 0$, and obeys the derivative bounds
$$ |\partial^s_\tau \partial^\alpha_\xi m(\tau,\xi)| \lesssim 2^{-js} |\xi|^{1-|\alpha|}$$
for arbitrarily many indices $s$, $\alpha$ away from the singular regions $\xi = 0$, $\tau = 0$.  Standard calculations then show that the kernel $K(t,x)$ obeys the bounds
$$ |K(t,x)| \lesssim (1 + |x|)^{-n-1} (1 + 2^j |t|)^{-2}$$
which is integrable as desired.
\end{proof}

\begin{lemma}\label{q-trunc}
The operators $Q_j$, $Q_{\leq j}$, $Q_{\geq j}$, $Q_{j_1 \leq \cdot \leq j_2}$, etc. are bounded on the spaces $L^p_t L^2_x$ for all $1 \leq p \leq \infty$.
\end{lemma}

\begin{proof}
It suffices to check $Q_{\leq j}$.  By scaling we may take $j=0$.  We then decompose
$$ Q_{\leq 0} = P_{\leq 10} Q_{\leq 0} + P_{>10} Q^+_{\leq 0} + P_{>10} Q^-_{\leq 0}.$$
The first term is disposable by Lemma \ref{disposable}, so by conjugation symmetry it suffices to demonstrate the boundedness of $P_{>10} Q^+_{\leq 0}$.  It suffices to show that the multiplier with symbol $m_0(\tau - |\xi|)$ is bounded on $L^p_t L^2_x$.  However, the transformation $U$ given by
$$ {\cal F} (U\phi)(\tau,\xi) := {\cal F} \phi(\tau-|\xi|,\xi)$$
is easily seen to be an isometry on $L^p_t L^2_x$ (indeed, for each time $t$, $U$ is just the unitary operator $e^{2\pi i t \sqrt{-\Delta}}$).  The claim then follows by conjugating by $U$ and observing that the multiplier with symbol $m_0(\tau)$ is bounded on $L^p_t L^2_x$.
\end{proof}

Of course, on $L^2_t L^2_x$ one does not need the above Lemmata, and can just use Plancherel to discard all multipliers with bounded symbols.

In the sequel we shall frequently be faced with estimating bilinear expressions such as $Q(\phi \psi)$, where $Q$ is some disposable multiplier and $\phi$ has a larger frequency support than $\psi$.  Often we shall use the $Q$ projection to localize $\phi$ to some frequency region (e.g. replacing $\phi$ by $\tilde Q \phi$ for some variant $\tilde Q$ of $Q$), and then use one of the above lemmata to discard the $Q$.  This effectively moves the multiplier $Q$ from $\phi \psi$ onto $\phi$, albeit at the cost of enlarging $Q$ slightly.  We shall exploit this trick in much the same way that Lemma \ref{commutator} was exploited in previous sections.

Because we wish to treat the two dimensional case $n=2$, we will not have access to good $L^2_t L^\infty_x$ or $L^2_t L^4_x$ Strichartz estimates, in contrast to \cite{tao:wavemap1}.  However, we have the following partial substitute (which basically comes from the $L^4_t L^\infty_x$ estimate):

\begin{lemma}[Improvement to Bernstein's inequality]\label{wimpy-strichartz} If $\phi$ has Fourier support in the region $D_0 \sim 2^k$, $D_- \sim 2^j$, then
$$ \| \phi \|_{L^2_t L^\infty_x} \lesssim \chi^{(4)}_{j \geq k} 2^{nk/2} \| \phi \|_{L^2_t L^2_x}.$$
\end{lemma}

\begin{proof}
We may rescale $j=0$.  We may assume that $k \gg 1$ since the claim follows from Bernstein's inequality \eqref{bernstein} otherwise. 

From the Poisson summation formula we can construct a Schwartz function $a(t)$ whose Fourier transform is supported on interval $\tau \ll 1$ and which satisfies 
$$ 1 = \sum_{s \in \Z} a^3(t-s)$$
for all $t \in \R$.  From H\"older we have
$$ \| \phi \|_{L^2_t L^\infty_x} = \| \sum_{S(c)} a^3(\cdot-s) \phi \|_{L^2_t L^\infty_x} \lesssim (\sum_{S(c)} \| a^2(\cdot-s) \phi \|_{L^2_t L^\infty_x}^2)^{1/2}
\lesssim (\sum_{S(c)} \| a(\cdot-s) \phi \|_{L^4_t L^\infty_x}^2)^{1/2}.$$
On the other hand, from the $L^4_t L^\infty_x$ Strichartz estimate (see e.g. \cite{tao:keel}) and the Fourier support of $a(\cdot-s) \phi$ we have
$$ \| a(\cdot - s) \phi \|_{L^4_t L^\infty_x} \lesssim 
\| a(\cdot - s) \phi \|_{X^{n/2-1/4,1/2+}_k} \sim
2^{-k/4} 2^{nk/2} \| a(\cdot - s) \phi \|_{L^2_t L^2_x},$$
and the claim follows by square-summing in $s$.
\end{proof}

One can improve this bound substantially when $n > 2$, but we shall not need to do so here.  The factor $\chi^{(4)}_{j \geq k}$ is a gain over Bernstein's inequality when $j < k$, and reflects the $\Lambda_p$ properties of the cone (as already indicated by Strichartz estimates). 

For any $s,b \in \R$, $k \in \Z$ and $1 \leq q \leq \infty$, we define $\dot X^{s,b,q}_k$ to be the completion of space of all Schwartz functions with Fourier support in $2^{k-5} \leq D_0 \leq 2^{k+5}$ whose norm
$$ \| \phi \|_{\dot X^{s,b,q}_k} := 2^{sk} [(\sum_j 2^{bj} \| Q_j \phi \|_{L^2_{t,x}})^q]^{1/q}$$
is finite, with the usual supremum convention when $q=\infty$.  

As a first approximation, one should think of $\phi_k$ as belonging to spaces of strength
comparable to $\nabla_{x,t}^{-1} \dot X^{n/2-1,1/2,q}_k$ for some $q$, while $\Box \phi_k$ belongs to spaces of strength comparable to $\dot X^{n/2-1,1/2-1,q}_k$.

The space $\dot X^{s,b,1}_k$ is an atomic space, whose atoms are functions with spacetime Fourier support in the region $\{ 2^{k-5} \leq D_0 \leq 2^{k+5}, D_- \sim 2^j \}$ and have an $L^2_t L^2_x$ norm of $O(2^{-sk} 2^{-j/2})$ for some integer $j$.  The verification of multi-linear estimates on these spaces then reduces to verifying the estimates on atoms.  Unfortunately only a portion of our waves $\phi$ can be placed into such a nice space.

From Plancherel we observe the duality relationship
\be{xsb-dual}
\sup\{ | \langle \phi, \psi \rangle | : \| \phi \|_{\dot X^{s,b,1}_k} \leq 1 \} \sim 2^{-k(s+s')} \| \phi \|_{\dot X^{s',-b,\infty}_k}
\end{equation}
for all $k, s, s', b$ and Schwartz $\psi \in \dot X^{s',-b,\infty}_k$, where 
$$ \langle \phi, \psi \rangle := \int\int \phi \overline \psi\ dx dt$$
is the usual inner product.

\section{Null frames}\label{frames-sec}

We now set out the machinery of null frames which we shall need to define our spaces.  In particular we shall describe the complementary null frame spaces $NFA[\kappa]$ and $S[k,\kappa]$, which shall play a key role in the construction
of $N_k$ and $S_k$ respectively.

We define a \emph{spherical cap} to be any subset $\kappa$ of $S^{n-1}$ of the form
$$ \kappa = \{ \omega \in S^{n-1}: |\omega - \omega_\kappa| < r_\kappa \}$$
for some $\omega_\kappa \in S^{n-1}$ and $0 < r_\kappa < 2$.  We call $\omega_\kappa$ and $r_\kappa$ the center and radius of $\kappa$ respectively.  If $\kappa$ is a cap and $C > 0$, we use $C\kappa$ to denote the cap with the same center but $C$ times the radius.  If $\pm$ is a sign, we use $\pm \kappa$ to denote the cap with the same radius but center $\pm\omega_\kappa$.

For any direction $\omega \in S^{n-1}$, we define the null direction $\theta_\omega$ by
$$ \theta_\omega := \frac{1}{\sqrt{2}}(1, \omega)$$
and the null plane $NP(\omega)$ by
$$ NP(\omega) := \{ (t,x) \in \R^{1+n}: (t,x) \cdot \theta_\omega = 0 \}$$
where $(t,x) \cdot (t',x') := tt' + x \cdot x'$ is the usual Euclidean inner product.

We can parameterize physical space $\R^{1+n}$ by the null co-ordinates $(t_\omega, x_\omega) \in \R \times NP(\omega)$ defined by
$$ t_\omega := (t,x) \cdot \theta_\omega; \quad x_\omega := (t,x) - t_\omega \theta_\omega.$$
One can similarly parameterize frequency space by 
$(\tau_\omega, \xi_\omega) \in \R \times NP(\omega)$ defined by
$$ \tau_\omega := (\tau,\xi) \cdot \theta_\omega; \quad \xi_\omega := (\tau,\xi) - \tau_\omega \theta_\omega.$$
We can then define Lebesgue spaces $L^q_{t_\omega} L^r_{x_\omega}$, $L^q_{\tau_\omega} L^r_{\xi_\omega}$ in the usual manner.  The trivial identity
\be{222}
\| \phi \|_{L^2_t L^2_x} \sim \| \phi \|_{L^2_{t_\omega} L^2_{x_\omega}}
\end{equation}
shall be crucial in linking the null frame estimates to the Euclidean frame estimates.

The null plane $NP(\omega)$ contains the null direction $(1, -\omega)$.  As such, it is not conducive to good energy estimates; for instance, control of $\Box \phi$ in $L^1_{t_\omega} L^2_{x_\omega}$ does not give satisfactory control of $\phi$ in general.  However, if one also knows that $\phi$ has Fourier support in a sector $\{\Theta \in \kappa\}$ and that $\omega$ is outside of $2\kappa$, then one can recover good energy estimates in this co-ordinate system (but losing a factor of $1/\dist(\omega,\kappa)$).  This motivates the definition of a Banach space $NFA[\kappa]$ of null frame atoms oriented away from $\kappa$.  More precisely:

\begin{definition}\label{nfa-def}
For any cap $\kappa$, we define $NFA[\kappa]$ to be the atomic Banach space whose atoms are functions $F$ with 
$$ \| F \|_{L^1_{t_\omega} L^2_{x_\omega}} \leq \dist(\omega,\kappa) $$
for some $\omega \not \in 2\kappa$.
\end{definition}

The spaces $NFA[\kappa]$ will be one of the more exotic building blocks for our nonlinearity space $N_k$, and appears implicitly in \cite{tataru:wave1}. 

Since the norm $NFA[\kappa]$ is defined entirely using Lebesgue spaces, we have the estimate
\be{lif}
\| \phi \psi \|_{NFA[\kappa]} \lesssim \| \phi \|_{L^\infty_t L^\infty_x}
\| \psi \|_{NFA[\kappa]}
\end{equation}
for all Schwartz $\phi$, $\psi$.
We also observe the nesting inequality
\be{inscribed}
\| F \|_{NFA[\kappa']} \leq \| F \|_{NFA[\kappa]}
\end{equation}
whenever $\kappa' \subset \kappa$.

We now construct a Banach space $S[k,\kappa]$ which is something of a dual space to $NFA[\kappa]$, and will be a key part of the construction of the spaces $S_k$ and $S(c)$.

\begin{proposition}\label{sklw-prop}  For each integer $k$ and cap $\kappa$ we can associate a Banach space $S[k,\kappa]$, such that the following properties hold for all integers $k$, $k'$, caps $\kappa$, $\kappa'$, and Schwartz functions $\phi$, $\psi$.

\begin{itemize}

\item ($S[k,\kappa]$ is stable under $L^\infty$ multiplications) We have
\be{liff}
\| \phi \psi \|_{S[k,\kappa]} \lesssim \| \phi \|_{L^\infty_t L^\infty_x} \| \psi \|_{S[k,\kappa]}.
\end{equation}

\item (Nesting property)
We have
\be{inscribed-sklw}
\| \phi \|_{S[k,\kappa]} \lesssim \| \phi \|_{S[k,\kappa']}
\end{equation}
whenever $\kappa' \subset \kappa$.

\item ($S[k,\kappa]$ is dimensionless)  For any integer $j$, we have
\be{dimensionless}
\| \phi(2^j t, 2^j x) \|_{S[k+j,\kappa]} = \| \phi(t,x) \|_{S[k,\kappa]}.
\end{equation}

\item (Energy estimate) We have
\be{sklw-energy}
\| \phi \|_{L^\infty_t L^2_x} \lesssim 2^{-nk/2} \| \phi \|_{S[k,\kappa]}.
\end{equation}

\item (Duality)  We have
\be{duality}
|\langle \phi, \psi\rangle| \lesssim 2^{-nk/2} \| \phi \|_{S[k,\kappa]} \| \psi \|_{NFA[\kappa]}
\end{equation}

\item (Consistency)  We have
\be{consistency}
\| \phi \|_{S[k,\kappa]} \sim \| \phi \|_{S[k',\kappa]}
\end{equation}
whenever $k' = k + O(1)$.

\item (Product estimates)  We have the estimates\footnote{The angular separation condition $2\kappa \cap 2\kappa' = \emptyset$ is only really needed when $n=2$.  When $n>2$ one the Strichartz estimates allow one to make $S(c)$ contain $L^4_t L^4_x$ type spaces, and $N_k$ contain $L^{4/3}_t L^{4/3}_x$ type spaces, in which case one can just use H\"older as a substitute for \eqref{NFAPW-dual}, \eqref{NFAPW}.  The factor $|\kappa'|^{1/2} \sim r_{\kappa'}^{(n-1)/2}$ allows us to obtain a gain in the small angle interaction case when $n > 1$.}

\be{NFAPW-dual}
\| \phi \psi \|_{NFA[\kappa]} \lesssim \frac{|\kappa'|^{1/2}}{2^{k'/2} \dist(\kappa,\kappa')}
\| \phi\|_{L^2_t L^2_x} \| \psi \|_{S[k',\kappa']}
\end{equation}
and
\be{NFAPW}
\| \phi \psi \|_{L^2_t L^2_x} \lesssim \frac{|\kappa'|^{1/2}}{
2^{nk/2} 2^{k'/2} \dist(\kappa, \kappa')}
\| \phi\|_{S[k,\kappa]} \| \psi \|_{S[k',\kappa']}
\end{equation}
whenever $2\kappa \cap 2\kappa' = \emptyset$.  

\item (Square-summability)  Let $0 < r \leq 2^{-5} r_\kappa$ be a radius, and suppose that $K$ is a finitely overlapping collection of caps of radius $r$ which are contained in $\kappa$.  For each $\kappa' \in K$ we associate a Schwartz function $\phi_{\kappa'}$ with Fourier support in the region $D_0 \sim 2^k$, $\Theta \in \kappa'$, $D_- \lesssim r^2 2^k$. Then we have
\be{sklw-ortho}
\| \sum_{\kappa \in K} \phi_{\kappa'} \|_{S[0,\kappa]}
\lesssim
(\sum_{\kappa \in K} \| \phi_{\kappa'} \|_{S[0,\kappa']}^2)^{1/2}
\end{equation}

\item ($S[k,\kappa]$ contains $\dot X^{n/2,1/2,1}_k$ waves with direction in $\kappa$) We have the estimate
\be{into-sklw}
\|\phi \|_{S[k,\kappa]} \lesssim \| \phi \|_{\dot X^{n/2,1/2,1}_k}
\end{equation}
whenever $\phi \in \dot X^{n/2,1/2,1}_k$ has Fourier support on $\{ \Theta \in \kappa \}$.

\end{itemize}
\end{proposition}

\begin{proof}

\divider{Step 1.  Construction of $S[k,\kappa]$.}

Define $NFA^*[\kappa]$ to be the space of functions $\phi$ whose norm
$$ \| \phi \|_{NFA^*[\kappa]} := \sup_{\omega \not \in 2\kappa} \dist(\omega,\kappa) 
\| \phi \|_{L^\infty_{t_\omega} L^2_{x_\omega}}$$
is finite.  Observe that $NFA^*[\kappa]$ is the dual of $NFA[\kappa]$.

We also define the plane wave space $PW[\kappa]$ to be the atomic Banach space whose atoms are functions $\phi$ with 
$$ \| \phi \|_{L^2_{t_{\omega}} L^\infty_{x_{\omega}}} \leq 1$$
for some $\omega \in \kappa$.

We now define $S[k,\kappa]$ by the norm
\be{sklw-def}
\| \phi \|_{S[k,\kappa]} := 
2^{nk/2} \| \phi \|_{NFA^*[\kappa]}
+ |\kappa|^{-1/2} 2^{k/2} \| \phi \|_{PW[\kappa]}
+ 2^{nk/2} \| \phi \|_{L^\infty_t L^2_x}.
\end{equation}

The estimates \eqref{liff}, \eqref{inscribed-sklw}, \eqref{duality}, \eqref{sklw-energy}, \eqref{consistency} and \eqref{dimensionless} are now straightforward and are left to the reader.
It remains to prove \eqref{NFAPW-dual}, \eqref{NFAPW}, \eqref{sklw-ortho}, and \eqref{into-sklw}.

\divider{Step 2.  Proof of the product estimates \eqref{NFAPW-dual}, \eqref{NFAPW}.}

The estimate \eqref{NFAPW} will follow from \eqref{NFAPW-dual}, \eqref{duality}, and duality, so we need only prove \eqref{NFAPW-dual}.  By \eqref{sklw-def} it suffices to show
$$
\| \phi \psi \|_{NFA[\kappa]} \lesssim \frac{1}{\dist(\kappa,\kappa')}
\| \phi\|_{L^2_t L^2_x} \| \psi \|_{PW[\kappa']}$$
But this follows by reducing $\psi$ to a $PW[\kappa']$ atom, then using H\"older and \eqref{222}.

\divider{Step 3.  Proof of the square-summability estimate \eqref{sklw-ortho}.}

Fix $\kappa$.  By \eqref{dimensionless} we may rescale $k=0$.  By \eqref{inscribed-sklw} we may replace the $S[0,\kappa']$ norms on the right-hand side by $S[0,\kappa]$.  We expand the left-hand side of \eqref{sklw-ortho} using \eqref{sklw-def} and treat the three components separately.

\divider{Step 3(a).  The contribution of $PW[\kappa]$.}

By triangle inequality and Cauchy-Schwarz we can bound this contribution by 
$$
\lesssim |\kappa|^{-1/2} |K|^{1/2}
(\sum_{\kappa' \in K} \| \phi_{\kappa'} \|_{PW[\kappa]}^2)^{1/2}.
$$
Since the caps in $K$ are finitely overlapping, we have $|K| \lesssim |\kappa|/r^{n-1}$, and the claim then follows from \eqref{sklw-def}.

\divider{Step 3(b).  The contribution of $L^\infty_t L^2_x$.}

This follows from Plancherel and the observation that the $\phi_{\kappa'}$ have finitely overlapping $\xi$-support as $\kappa$ varies.

\divider{Step 3(c).  The contribution of $NFA^*[\kappa]$.}

Suppose that $\omega$ is a direction such that $\omega \not \in 2\kappa$.  Then from elementary geometry we see that the functions $\phi_{\kappa'}$ have finitely overlapping $\xi_\omega$-Fourier support as $\kappa'$ varies.  We thus have from Plancherel that
$$
\| \sum_{\kappa' \in K}
\phi_{\kappa'}(t_\omega) \|_{L^2_{x_\omega}}
\lesssim
(\sum_{\kappa' \in K}
\| \phi_{\kappa'}(t_\omega) \|_{L^2_{x_\omega}}^2)^{1/2}
$$
for all $t_\omega \in \R$.
Taking suprema in $t_\omega$, then taking suprema in $\omega$ and using \eqref{sklw-def} we obtain the claim.

\divider{Step 4.  Proof of the embedding \eqref{into-sklw}.}

By \eqref{dimensionless} we may rescale $k=0$.  By the atomic nature of $\dot X^{n/2,1/2,1}_0$ and conjugation symmetry we may assume that $\phi$ has Fourier support in the region $\tau > 0, D_0 \sim 1, D_- \sim 2^j, \Theta \in \kappa$ for some integer $j$, in which case we need to show
$$ \| \phi \|_{S[k,\kappa]} \lesssim 2^{j/2} \| \phi \|_{L^2_t L^2_x}.$$

Fix $j$.  We now prove the claim for the three components of \eqref{sklw-def} separately.

\divider{Step 4(a).  The estimation of the $NFA^*[\kappa]$ norm.}

By Plancherel and Minkowski we have
\be{icky}
\| \phi \|_{L^\infty_t L^2_x} \lesssim \| {\cal F}  \phi \|_{L^2_\xi L^1_\tau}.
\end{equation}
For fixed $\xi$, the function ${\cal F}  \phi$ has $\tau$-support in an interval of length $O(2^j)$.  The claim then follows by using H\"older and Plancherel.

\divider{Step 4(b).  The estimation of the $NFA^*[\kappa]$ norm.}

Let $\omega \in S^{n-1}$ be such that $\omega \not \in 2\kappa$.  
By Plancherel and Minkowski we have
$$ \| \phi \|_{L^\infty_{t_\omega} L^2_{x_\omega}} \lesssim \| {\cal F}  \phi \|_{L^2_{\xi_\omega} L^1_{\tau_\omega}}.$$
From elementary geometry we see that for fixed $\xi_\omega$, the function ${\cal F}  \phi$ has $\tau_\omega$-support in an interval of length $O(2^j/\dist(\omega,\kappa)^2)$.  The claim then follows by using H\"older and Plancherel, and then taking suprema in $\omega$.

\divider{Step 4(c).  The estimation of the $PW[\kappa]$ norm.}

From the Fourier inversion formula in polar co-ordinates and the Fourier support of $\phi$ we have
$$ \phi(t,x) = C \int_{\omega \in \kappa} \int_{|a| \lesssim 2^j}
e^{2\pi i a t} \int_{r \sim 1} {\cal F} \phi(r+a,r\omega) 2^{\pi i r (t,x) \cdot (1,\omega)} r^{n-1}\ dr da d\omega.$$
From Minkowski's inequality and the definition of $PW[\kappa]$ we thus have
$$ \| \phi \|_{PW[\kappa]} \lesssim \int_{\omega \in \kappa} \int_{|a| \lesssim j}
\| e^{2\pi i a t} \int_{r \sim 1} {\cal F}  \phi(r+a,r\omega) 2^{\pi i r (t,x) \cdot (1,\omega)} r^{n-1}\ dr\|_{L^2_{t_\omega} L^\infty_{x_\omega}} da d\omega.$$
The $e^{2\pi i at}$ factor is bounded and can be discarded.
Since $(t,x) \cdot (1,\omega) = \sqrt{2} t_\omega$, we can estimate the previous by
$$
\lesssim \int_{\omega \in \kappa} \int_{|a| \lesssim 2^j}
\| \int_{r \sim 1} {\cal F}  \phi(r+a,r\omega) e^{2\pi i \sqrt{2} r t_\omega} r^{n-1}\ dr\|_{L^2_{t_\omega}} da d\omega.$$
By Plancherel this is bounded by
$$
\lesssim \int_{\omega \in \kappa} \int_{|a| \lesssim 2^j}
(\int_{r \sim 1} ({\cal F} \phi(r+a,r\omega) r^{n-1})^2\ dr)^{1/2} da d\omega.$$
By Cauchy-Schwarz this is bounded by
$$
\lesssim |\kappa|^{1/2} 2^{j/2} 
(\int_{\omega \in \kappa} \int_{|a| \lesssim 2^j}
\int_{r \sim 1} ({\cal F} \phi(r+a,r\omega) r^{n-1})^2\ dr da d\omega)^{1/2}.$$
Undoing the polar co-ordinates we can estimate this by
$$
\lesssim |\kappa|^{1/2} 2^{j/2} \| {\cal F}  \phi \|_{L^2_\tau L^2_\xi}.$$
By Plancherel we thus see that the contribution of the $PW[\kappa]$ norm is acceptable.  This finishes the proof of \eqref{into-sklw}.
\end{proof}

\section{Sector projections}\label{sector-sec}

We now define some sector projection operators which will be needed in the definition of $N_k$ and $S(c)$.

For every real $l > 10$ we let $\Omega_l$ be a maximal $2^{-l}$-separated subset of the unit sphere $S^{n-1}$.  We let $K_l$ denote the space of caps
$$K_l := \{ \kappa: \omega_\kappa \in \Omega_l; r_\kappa = 2^{-l+5} \}.$$

For any real $l > 10$ and any cap $\kappa \in K_l$, we defnine a bump function $m''_\kappa$ on $S^{n-1}$ such that $1 = \sum_{\kappa \in K_l} m''_\kappa$.  For any integer $k$, we then define $P_{k,\kappa}$ to be the spatial Fourier multiplier with symbol
$$m_{k,\kappa}(\xi) := m'_k(|\xi|) m''_\kappa(\xi/|\xi|),$$ 
where $m'_k$ is a bump function adapted to $2^{k-4} \leq D_0 \leq 2^{k+4}$ which equals 1 on $2^{k-3} \leq D_0 \leq D^{k+3}$.  Thus $m_{k,\kappa}$ is a bump function adapted to the tube 
\be{tube}
\{ \xi \in \R^n: 2^{k-4} \leq D_0 \leq 2^{k+4}; \frac{\xi}{|\xi|} \in \frac{1}{2}\kappa \}
\end{equation}
andt $\phi = \sum_{\kappa \in K_l} P_{k,\kappa} \phi$ for all functions $\phi$ with Fourier support in $2^{k-3} \leq D_0 \leq 2^{k+3}$. We shall also need the variant $\tilde P_{k,\kappa}$, given by a spatial Fourier multiplier $\tilde m_{k,\kappa}$ which is a bump function equals 1 on \eqref{tube} and is adapted to the enlargement
$$
\{ \xi \in \R^n: 2^{k-5} \leq D_0 \leq 2^{k+5}; \frac{\xi}{|\xi|} \in \kappa \}
$$
of \eqref{tube}.  Observe that the multipliers $P_{k,\kappa}$ and $\tilde P_{k,\kappa}$ are disposable.

\begin{lemma}\label{proj-disposable}
Let $j, k \in \Z, l \in \R$ be such that $l > 10$ and $j \geq k-2l + O(1)$, and let $\kappa = K_l$.  Then $\tilde P_{k,\kappa} Q_j$, $\tilde P_{k,\kappa} Q_{\leq j}$, and $\tilde P_{k,\kappa} Q_{\geq j}$ are disposable multipliers.  If we also assume that $|j-k| > 10$, then $\tilde P_{k,\kappa} Q^\pm_j$, $\tilde P_{k,\kappa} Q^\pm_{\leq j}$, and $\tilde P_{k,\kappa} Q^\pm_{\geq j}$ are also disposable.  The above results also hold when $\tilde P_{k,\kappa}$ is replaced by $P_{k,\kappa}$.
\end{lemma}

\begin{proof}
We shall just prove these claims for $\tilde P_{k,\kappa}$, as the $P_{k,\kappa}$ claims then follow from the factorization $P_{k,\kappa} = P_{k,\kappa} \tilde P_{k,\kappa}$.

If $|j-k| \leq 5$ then we factorize $\tilde P_{k,\kappa} Q_{\leq j} = \tilde P_{k,\kappa} (P_{k-10 < \cdot < k+10} Q_{\leq j})$, which is disposable by Lemma \ref{disposable}.  Similarly for $\tilde P_{k,\kappa} Q_j$ and $\tilde P_{k,\kappa} Q_{\geq j}$.

Now suppose that $|j-k| > 5$.  It suffices to show that $\tilde P_{k,\kappa} Q^\pm_{\leq j}$ is disposable, since the remaining claims then follow by taking suitable linear combinations of these multipliers.  By conjugation symmetry we may take $\pm = +$.

It is then an easy matter to verify that the symbol of $\tilde P_{k,\kappa} Q^+_{\leq j}$ is a bump function adapted to the parallelopiped
$$ \{ (\tau,\xi): \xi \cdot \omega_\kappa \sim 1; |\xi \wedge \omega_\kappa| \lesssim 2^{-l};
\tau = \xi \cdot \omega_\kappa + O(2^{-2l}) \}.$$
The kernel is then rapidly decreasing away from the dual of this parallelopiped, and the claim follows.
\end{proof}

\section{Construction of $N_k$, $S_k$, $S(c)$.}\label{construction-sec}

In this section we construct the spaces $N_k$, $S_k$, $S(c)$ required for Theorem \ref{spaces}, and prove some of the easier properties about them.  The energy estimate \eqref{energy-est-2}, the bilinear estimates \eqref{algebra}, \eqref{sk-skp}, \eqref{sk-sk}, \eqref{null}, \eqref{ur-algebra}, \eqref{ur-algebra-sk}, and the trilinear estimate \eqref{o-lemma} are quite lengthy to prove and will be deferred to later sections.

The spaces $N_k$, $S_k$, $S(c)$ will be defined in terms of frequency localized variants $N[k]$, $S[k]$.  We begin with the definition of the $S(c)$-type spaces.

\begin{definition}\label{sk-space-def}
For any integer $k$, we define $S[k]$ to be the completion of the space of Schwartz functions with Fourier support in the region $\{ \xi \in \R^n: 2^{k-3} \leq D_0 \leq 2^{k+3}\}$ with respect to the norm
\begin{eqnarray}\label{sk-def} 
\| \phi \|_{S[k]} := &
\| \nabla_{x,t} \phi \|_{L^\infty_t \dot H^{n/2-1}_x} +
\| \nabla_{x,t} \phi \|_{\dot X^{n/2-1,1/2,\infty}_k} \\
+ &
\sup_\pm \sup_{l>10} (\sum_{\kappa \in K_l} \| P_{k,\pm \kappa} Q^\pm_{<k-2l} \phi \|_{S[k,\kappa]}^2)^{1/2}.\nonumber
\end{eqnarray}
\end{definition}

The first two terms in \eqref{sk-def} are standard, but are not sufficient by themselves to obtain good product estimates\footnote{More precisely, a logarithmic divergence occurs whenever one considers the imbalanced modulation case (as defined in Section \ref{geometry-sec}).}.  This will be rectified by the addition of the third term in \eqref{sk-def}, which also essentially appears in \cite{tataru:wave2}.  Note that the square sum in the third term of \eqref{sk-def} is essentially increasing in $l$ thanks to \eqref{sklw-ortho}.

It is not a priori clear that the third term in \eqref{sk-def} is finite even for Schwartz functions.  However, this follows from

\begin{lemma}\label{f-lemma}
Suppose $\phi$ is a Schwartz function with Fourier support in $\{ 2^{k-5} \leq D_0 \leq 2^{k+5} \}$.  Then we have
\be{f-bound}
(\sum_{\kappa \in K_l} \| P_{k,\pm \kappa} Q^\pm_{<-2l+O(1)} \phi \|_{S[k,\kappa]}^2)^{1/2}
\lesssim \| \phi \|_{\dot X^{n/2,1/2,1}_k}
\end{equation}
for all $l > 0$ and any sign $\pm$.
\end{lemma}

\begin{proof}
By \eqref{into-sklw} it suffices to show
$$
 (\sum_{\kappa \in K_l} \| P_{k,\pm \kappa} Q^\pm_{<-2l+O(1)} \phi \|_{\dot X^{n/2,1/2,1}_k}^2)^{1/2}
\lesssim \| \phi \|_{\dot X^{n/2,1/2,1}_k}
$$
By the atomic nature of $\dot X^{n/2,1/2,1}_k$ we may assume that $\phi$ has Fourier support in the region $D_- \sim 2^j$ for some integer $j$.  But then the claim then follows from Plancherel's theorem and the finite overlap of the $P_{k,\kappa}$ Fourier supports.
\end{proof}

A function $\phi$ in $S[k]$ enjoys several estimates.  From \eqref{sk-def} and Lemma \ref{q-trunc} we clearly have
\be{sk-energy}
\| \phi \|_{L^\infty_t L^2_x}, \| Q_{\leq j} \phi \|_{L^\infty_t L^2_x}, \| Q_j \phi \|_{L^\infty_t L^2_x}, etc. \lesssim 2^{-nk/2} \| \phi \|_{S[k]}
\end{equation}
and hence from Bernstein's inequality \eqref{bernstein} that
\be{sk-infty}
\| \phi \|_{L^\infty_t L^\infty_x}, \| Q_{\leq j} \phi \|_{L^\infty_t L^\infty_x}, \| Q_j \phi \|_{L^\infty_t L^\infty_x}, etc.
 \lesssim \| \phi \|_{S[k]}.
\end{equation}
Also, from \eqref{sk-def} we have
\be{qbound}
\| Q_j \phi \|_{L^2_t L^2_x} \lesssim 2^{-nk/2} 2^{-j/2} \frac{2^k}{2^k + 2^j} \| \phi \|_{S[k]}
\end{equation}
since $\nabla_{x,t}$ behaves on $Q_j \phi$ like the quantity $D_+ \sim 2^k + 2^j$, at least as far as the $L^2_t L^2_x$ norm is concerned.  From this and Lemma \ref{wimpy-strichartz} we have
\be{qbound-strichartz}
\| Q_j \phi \|_{L^2_t L^\infty_x} \lesssim \chi^{(4)}_{j \geq k} \min(2^{-(j-k)},1) 2^{-j/2}
\| \phi \|_{S[k]}.
\end{equation}
These four estimates will be used very frequently in the sequel, and between them are capable of handling all cases except when two of $\phi$, $\psi$, $\phi \psi$ are close to the light cone\footnote{More precisely, in the language of Section \ref{geometry-sec}, we shall use \eqref{sk-energy}, \eqref{sk-infty}, \eqref{qbound}, \eqref{qbound-strichartz} when the modulations of $\phi$, $\psi$, $\phi \psi$ are balanced, and rely on Lemma \ref{sklw-prop} when the modulations are imbalanced.}, in which case one must resort to the null frame estimates in Lemma \ref{sklw-prop} instead.  It is worth mentioning that of the above four estimates, the $L^2_x$ estimates are more effective for high frequencies, while the $L^2_t$ estimates are more effective for large modulations.  For instance, if $\phi$ has low frequency and small modulation, then \eqref{sk-infty} is probably the best of the above four estimates to use.

The space $S[k]$ also contains $\dot X^{n/2,1/2,1}_k$ functions:

\begin{lemma}\label{xn11 contain}
If $\phi$ is in $\dot X^{n/2,1/2,1}_k$, then
$$ \| \phi \|_{S[k]} \lesssim  \| 2^{-k} \nabla_{x,t} \phi \|_{\dot X^{n/2,1/2,1}_k}.$$
\end{lemma}

\begin{proof}
We may rescale $k=0$.  By the atomic nature of $\dot X^{n/2,1/2,1}_0$ and conjugation symmetry we may assume that $\phi$ has Fourier support in the region $\tau > 0, D_0 \sim 1, D_- \sim 2^j, \Theta \in \kappa$ for some integer $j$, in which case we need to show
$$ \| \phi \|_{S[0,\kappa]} \lesssim (1 + 2^j) 2^{j/2} \| \phi \|_{L^2_t L^2_x}.$$
To control the first term of \eqref{sk-def}, we use \eqref{icky} to obtain 
$$
\| \nabla_{x,t} \phi \|_{L^\infty_t L^2_x}
\lesssim (1 + 2^j) \| \phi \|_{L^\infty_t L^2_x}
\lesssim (1 + 2^j) \| {\cal F}  \phi \|_{L^2_\xi L^1_\tau}.
$$
Since ${\cal F}  \phi$ is supported in an interval of length $O(2^j)$, the claim then follows from H\"older and Plancherel.

The second term of \eqref{sk-def} is trivial, so it remains only to show
$$(\sum_{\kappa \in K_l} \| P_{0,\pm \kappa} Q^\pm_{<-2l} \phi \|_{S[0,\kappa]}^2)^{1/2}
\lesssim (1 + 2^j) 2^{j/2} \| \phi \|_{L^2_t L^2_x}$$
for all $l > 10$ and signs $\pm$.

Fix $l$, $\pm$.  We may assume that $j \leq -2l+O(1)$ since the left-hand vanishes otherwise.  In particular we have $1 + 2^j \sim 1$.  By \eqref{into-sklw} we may estimate the left-hand side by
$$\lesssim 2^{j/2} (\sum_{\kappa \in K_l} \| P_{0,\pm \kappa} \phi \|_{L^2_t L^2_x}^2)^{1/2},$$
and the claim follows from  Plancherel.
\end{proof}

We also have the technical estimate 

\begin{lemma}\label{tech-lemma}
For all $\phi \in S[k']$ and $k' = k+O(1)$ we have
\be{tech}
\| P_k \phi \|_{S[k]} \lesssim \| \phi \|_{S[k']}.
\end{equation}
\end{lemma}

\begin{proof}
This estimate is easily verified for the first two factors of \eqref{sk-def}.  Thus it only remains to show 
$$(\sum_{\kappa \in K_l} \| P_{k,\pm \kappa} Q^\pm_{<k-2l} P_k \phi \|_{S[k,\kappa]}^2)^{1/2}
\lesssim \| \phi \|_{S[k']}$$
for all $l > 10$ and signs $\pm$.

Fix $\pm$, $l$.  We split $Q^\pm_{<k-2l} = Q^\pm_{<k'-2l} + (Q^\pm_{<k-2l}-Q^\pm_{<k'-2l})$.  To treat the second term we use \eqref{into-sklw} to estimate 
this contribution by
$$  \lesssim 2^{(k-2l)/2} (\sum_{\kappa \in K_l} \| P_{k,\pm \kappa} (Q^\pm_{<k-2l} - Q^\pm_{<k'-2l}) \phi \|_{L^2_t L^2_x}^2)^{1/2},$$
which by Plancherel is bounded by
$$  \lesssim 2^{(k-2l)/2}  \| (Q^\pm_{<k-2l} - Q^\pm_{<k'-2l}) \phi \|_{L^2_t L^2_x}.$$
But this is acceptable by the second term in \eqref{sk-def}.

It remains to estimate
$$(\sum_{\kappa \in K_l} \| P_{k,\pm \kappa} Q^\pm_{<k'-2l} P_k \phi \|_{S[k,\kappa]}^2)^{1/2}.$$
From the construction of $P_{k,\pm \kappa}$ and the Fourier support of $\phi$ we observe the identity
$$ P_{k,\pm \kappa} P_k \phi = P_{k',\pm \kappa} P_k \phi.$$
Discarding the $P_k$ we then see that this contribution is acceptable by the third term in \eqref{sk-def}.
\end{proof}

We are now in a position to define the full space $S(c)(\R^{1+n})$.

We define the space $S(c)(\R^{1+n})$ to be the closure of the space of functions $\phi$ such that $\phi-e$ is a Schwartz function for some constant $e$ (so that $\phi$ is asymptotically constant) and whose norm
\be{s-def}
\| \phi \|_{S(c)(\R^{1+n})} := \| \phi \|_{L^\infty_t L^\infty_x} + \sup_k c_k^{-1} \| \phi_k \|_{S[k]}
\end{equation}
is finite.  We can then form the restricted space $S(c) := S(c)([-T,T] \times \R^n)$ in the usual manner.

We then define the space $S_k(\R^{1+n})$ to be the subspace of $S(c)(\R^{1+n})$ given by the norm
\be{S_k-def}
\| \phi \|_{S_k(\R^{1+n})} := \sup_{k'} 2^{\delta_1|k-k'|} \| \phi_k \|_{S[k]},
\end{equation}
and define $S_k := S_k([-T,T]\times \R^n)$ in the usual manner.  Thus $S_k$ is a larger variant of $S[k]$ which allows for some leakage outside of the frequency region $D_0 \sim 2^k$.  Most of the estimates in this half of the paper will contain a decay of $\chi^{(2)}$ or stronger, which shall be more than enough to overcome this $2^{\delta_1|k-k'|} = \chi^{-(1)}_{k=k'}$ leakage. 

We shall also define the slightly weaker norm $S(1)(\R^{1+n})$ by
\be{stil-def}
\| \phi \|_{S(1)(\R^{1+n})} := \| \phi \|_{L^\infty_t L^\infty_x} + \sup_k  \| \phi_k \|_{S[k]}.
\end{equation}
From \eqref{s-def}, \eqref{sk-infty}, \eqref{S_k-def}, and Lemma \ref{tech} we observe the embeddings
\be{stil-s}
S(c)(\R^{1+n}), S_k(\R^{1+n}), S[k] \subseteq S(1)(\R^{1+n}).
\end{equation}
The space $S(1)$ will therefore be convenient for unifying the treatment of $S(c)(\R^{1+n})$ and $S[k]$ in the estimates in Theorem \ref{spaces}. 

Certain portions of Theorem \ref{spaces} can now be easily verified.  It is easy to see that $S(c)$ and $S_k$ satisfy the required invariance properties, as well as \eqref{infty-control}, \eqref{insensitive}, and \eqref{s-sk}.  It is also clear that $S(c)$ contains the identity function 1.

Having defined the Banach algebra $S(c)$, which will be used to hold the wave map $\phi$, we now define the spaces $N[k]$ and $N_k$, which are used to hold the renormalized frequency pieces of $\Box \phi$.  The invariance properties of $N_k$ and \eqref{l12} shall be immediate, however the proof of \eqref{energy-est-2} is far more involved, and shall be deferred to the next section.

\begin{definition}\label{atom-def} Let $k$ be an integer, and let $F$ be a Schwartz function with Fourier support in the region $2^{k-4} \leq D_0 \leq 2^{k+4}$.  We say that $F$ is an \emph{$L^1_t \dot H^{n/2-1}_x$-atom at frequency $2^k$} if
$$ \| F \|_{L^1_t L^2_x} \leq 2^{-(\frac{n}{2}-1)k}.$$
If $j \in \Z$, we say that $F$ is a \emph{$\dot X^{n/2-1,-1/2,1}$-atom with frequency $2^k$ and modulation $2^j$} if $F$ has Fourier support in the region $2^{k-4} \leq D_0 \leq 2^{k+4}$, $2^{j-5} \leq D_- \leq 2^{j+5}$ and
$$ \| F \|_{L^2_t L^2_x} \leq 2^{j/2} 2^{-(\frac{n}{2}-1)k}.$$
Finally, if $l > 10$ is a real number and $\pm$ is a sign, we say that $F$ is a \emph{$\pm$-null frame atom with frequency $2^k$ and angle $2^{-l}$} if there exists a decomposition $F = \sum_{\kappa \in K_l} F_\kappa$ such that each $F_\kappa$ has Fourier support in the region
$\{ (\tau, \xi): \pm \tau > 0; D_- \leq 2^{k-2l-50}; 2^{k-4} \leq D_0 \leq 2^{k+4}; \Theta \in \frac{1}{2} \kappa \}$
and
$$ (\sum_{\kappa \in K_l} \| F_\kappa \|_{NFA[\kappa]}^2)^{1/2} \leq 2^{-(\frac{n}{2}-1)k}.$$
We let $N[k]$ be the atomic Banach space generated by the $N[k]$ atoms.
\end{definition}

For future reference we observe that if $F$ is a $\pm$-null frame atom with frequency $2^k$ and angle $2^{-l}$, and $F_\kappa$ is as above, then
\be{ftpf}
F_\kappa = \tilde P_{k,\kappa} F_\kappa.
\end{equation}

From our definition we clearly have the estimates
\be{f-l12}
\| F \|_{N[k]} \lesssim \| F \|_{L^1_t \dot H^{n/2-1}_x} 
\sim 2^{(n/2-1)k} \| F\|_{L^1_t L^2_x}
\end{equation}
and
\be{f-xsb}
\| F \|_{N[k]} \lesssim \| F\|_{\dot X^{n/2-1,-1/2,1}_k}
\end{equation}
whenever $F$ is a Schwartz function with Fourier support in $2^{k-4} \leq D_0 \leq 2^{k+4}$.  Furthermore, we have
\be{f-null}
\| \sum_{\kappa \in K_l} F_\kappa \|_{N[k]}
\lesssim 2^{(n/2-1)k}
(\sum_{\kappa \in K_l} \| F_\kappa \|_{NFA[\kappa]}^2)^{1/2}
\end{equation}
for all signs $\pm$, $l > 10$ and $F_\kappa$ with Fourier support in the region $\{ (\tau, \xi): \pm \tau > 0; D_- \leq 2^{k-2l-50}; 2^{k-4} \leq D_0 \leq 2^{k+4}; \Theta \in \frac{1}{2}\kappa \}$.

Often we will need to estimate the $N_k$ norm of a function $F$ which is a bilinear expression of two other functions.  In such cases, the estimate \eqref{f-xsb} is the most favorable when $F$ has modulation at least as large as its component functions, while \eqref{f-l12} is favorable when $F$ has modulation much smaller than its components.  The remaining bound \eqref{f-null} is needed when $F$ and one of its components have small modulation, while the other component has large modulation. In this case one must use the geometry of the cone to create some angular separation between the the two small modulation functions and then use \eqref{NFAPW-dual}.

\begin{definition}\label{nk-def}  Let $F$ be a Schwartz function and $k$ be an integer.  We say that $F$ is an \emph{$N_k$ atom} if there exists a $k' \in \Z$ such that $2^{100n|k-k'|} F$ is a $N[k]$ atom.  We define $N_k(\R^{1+n})$ to be the atomic Banach space generated by the $N_k$ atoms.  Finally, we define $N_k := N_k([-T,T] \times \R^n)$ to be the restriction of $N_k(\R^{1+n})$ to the slab $[-T,T] \times \R^n$.
\end{definition}

Thus $N_k$ is a very slight enlargement of $N[k]$ which (reluctantly!) allows for some frequency leakage outside the region $D_0 \sim 2^k$.  Observe that the invariance properties of $N_k$ required in Theorem \ref{spaces} are immediate, as are \eqref{l12} and \eqref{compat}.

From \eqref{duality}, \eqref{xsb-dual}, and H\"older, we observe the useful duality property
\be{sf-duality}
|\langle \phi, F \rangle| \lesssim 2^{-(n-1)k} \| \phi \|_{S[k]} \| F \|_{N[k']} \end{equation}
whenever $k' = k+O(1)$ and $\phi \in S[k]$, $F \in N[k']$.  Thus up to a scaling factor, $S[k]$ is contained in the dual of $N[k]$ and vice versa.

From Lemma \ref{xn11 contain}, \eqref{sf-duality}, and duality we obtain the useful embedding

\begin{lemma}\label{n-xsb}
For all $F \in N[k]$ we have
$$ \| F \|_{\dot X^{n/2-1,-1/2,\infty}_k} \lesssim \| F \|_{N[k]}.$$
\end{lemma}

In particular, we may reverse \eqref{f-xsb} if $D_-$ is restricted to a single dyadic block.

\section{Energy estimates: The proof of \eqref{energy-est-2}}\label{nk-sec}

In this section we prove \eqref{energy-est-2}.  We remark that this section can be read independently since the techniques used here are not needed elsewhere in the paper.

\divider{Step 0.  Scaling.}

By \eqref{S_k-def} and Definition \ref{nk-def} it suffices to show that
$$
\|\phi_{k'} \|_{S[k']([-T,T] \times \R^n)} \lesssim \|\Box \phi \|_{N[k] ([-T,T] \times \R^n)} + \| \phi[0] \|_{\dot H^{n/2} \times \dot H^{n/2-1}}.
$$
for all $k$, $k'$ with $k' = k + O(1)$.  By rescaling $k$, $k'$, and $T$ we may make $k' = 0$, so that $k = O(1)$.  We may assume that $T \geq 1$, since the case $T < 1$ then follows by restricting the $T \geq 1$ estimate to a smaller slab.  For any $M > 0$, we let $\eta_M(t)$ be a bump function adapted $[-10M, 10M]$ and equals 1 on $[-5M, 5M]$.

Fix $k$.  By linearity we can split into $\phi$ into a free solution and a solution to the inhomogeneous problem, and we address the two cases separately.

\divider{Step 1.  (Free solutions) Prove \eqref{energy-est-2} when $\Box \phi = 0$ on $[-T,T] \times \R^n$.}

We may assume that $\phi[0]$ is Schwartz.  Of course we may then extend $\phi$ to be a free solution on all of $\R^{1+n}$.  By Plancherel's theorem we then observe the Fourier representation
\be{fourier}
{\cal F} \phi_0(\tau,\xi) = f_+(\xi) \delta(\tau - |\xi|) + f_-(\xi) \delta(\tau + \xi)
\end{equation}
where $f_+$, $f_-$ are supported on $2^{-3} \leq D_0 \leq 2^3$ and obey
$$ \| f_+ \|_2 + \|f_-\|_2 \lesssim \| \phi[0] \|_{\dot H^{n/2} \times \dot H^{n/2-1}}.$$
It thus suffices to show that
$$ \| \eta_T(t) \phi_0 \|_{S[0](\R^{1+n})} \lesssim \|f_+\|_2 + \|f_-\|_2.$$
By Lemma \ref{xn11 contain} it suffices to show
$$ \| \eta_T(t) \phi_0 \|_{\dot X^{n/2,1/2,1}_0} \lesssim \|f_+\|_2 + \|f_-\|_2.$$
But this follows by applying the identity
$$
{\cal F} (\eta_T \phi_0)(\tau,\xi) = f_+(\xi) \hat \eta_T(\tau - |\xi|) + f_-(\xi) \hat \eta_T(\tau + |\xi|)
$$
coming from \eqref{fourier}, then using Plancherel's theorem and a routine computation.  

\divider{Step 2.  (Inhomogeneous solutions) Prove \eqref{energy-est-2} when $\phi[0] = 0$.}

As in \cite{tataru:wave2}, the idea shall be to use Duhamel's formula to write the inhomogeneous solution as an average of free solutions truncated along half-spaces.  Because $N[k]$ contains null frame atoms, we shall need to perform some truncation along null planes $NP(\omega)$, which becomes a little messy.

We turn to the details.  By Duhamel's formula it suffices to show that 
$$ \| P_0 \int_{-\infty}^t \frac{\sin((t-s)\sqrt{-\Delta})}{\sqrt{-\Delta}} \chi_{[0,\infty)}(s) F(s)\ ds \|_{S[0]([-T,T] \times \R^n)} \lesssim \| F \|_{N[k]([-T,T] \times \R^n)}$$
for all $F \in N[k]([-T,T] \times \R^n)$. 

We may of course extend $F$ to be an element of $N[k](\R^{1+n})$ with comparable norm.  

To prove this estimate we begin by disposing of the $\chi_{[0,\infty)}$ cutoff:

\divider{Step 2(a).  Show that $\|\chi_{[0,\infty)}(t) F\|_{N[k]} \lesssim \| F \|_{N[k]}$.}

By scaling and linearity it suffices to show that
$$ \| \sgn(t) F \|_{N[0]} \lesssim \| F\|_{N[0]}.$$
We may of course assume that $F$ is a $N[0]$ atom.  We then split into three cases depending on what type of atom $F$ is.

\divider{Case 2(a).1.  $F$ is an $L^1_t \dot H^{n/2-1}_x$ atom with frequency 1.}

This embedding is trivial since multiplication by $\sgn(t)$ does not affect the $L^1_t L^2_x$ norm.

\divider{Case 2(a).2.  $F$ is a $\dot X^{n/2-1,-1/2,1}$ atom with frequency 1 and modulation $2^j$.}

In this case $F$ has Fourier support on the region $\{ D_0 \sim 1, D_- \sim 2^j \}$ and
$$ \| F \|_{L^2_t L^2_x} \lesssim 2^{-j/2}.$$
We partition $\sgn(t) = \sgn_{<j-C}(t) + \sgn_{\geq j-C}(t)$, where $\sgn_{<j-C}(t)$ is $\sgn(t)$ smoothly cut off in Fourier space to the region $|\tau| \lesssim 2^{j-C}$, and similarly for $\sgn_{\geq j-C}$.  Observe that $\sgn_{<j-C}(t) F$ has Fourier support in the region $D_0 \sim 1, D_- \sim 2^j$, and hence by \eqref{f-xsb} and the boundedness of $\sgn_{<j-C}(t)$ we have
$$ \| \sgn_{<j-C}(t) F \|_{N_0} \lesssim 2^{j/2} \| \sgn_{<j-C}(t) F \|_{L^2_t L^2_x} \lesssim 2^{j/2} \| F\|_{L^2_t L^2_x} \lesssim 1$$
as desired.

To control $\sgn_{\geq j-C}(t) F$ we use \eqref{f-l12} and H\"older to estimate 
$$ \| \sgn_{\geq j-C}(t) F \|_{N[0]} \lesssim \| \sgn_{\geq j-C}(t) \|_{L^2_t} \| F \|_{L^2_t L^2_x} \lesssim 2^{j/2} \| \sgn_{\geq j-C} \|_2.$$
But a routine computation shows that $\sgn_{\geq j-C}$ is bounded and is rapidly decreasing away from the interval $[-2^{-j+C}, 2^{-j+C}]$, and the claim follows.

\divider{Case 2(a).3.  $F$ is a $\pm$-null frame atom with frequency 1 and angle $2^{-l}$.}

By conjugation symmetry we may take $\pm = +$.  We then have a decomposition $F = \sum_{\kappa \in K_l} F_\kappa$ such that each $F_\kappa$ has Fourier support in the region

$$
\{ (\tau, \xi): \tau > 0; D_0 \sim 1, D_- \lesssim 2^{-2l}; \Theta \in \frac{1}{2} \kappa \}
$$
and
\be{chunky}
(\sum_{\kappa \in K_l} \| F_\kappa \|_{NFA[\kappa]}^2)^{1/2} \lesssim 1.
\end{equation}

Let $C$ be a large constant.  The function $Q_{>-2l-5C} F$ has Fourier support in the region $D_0 \sim 1, D_- \sim 2^{-2l}$.  By Lemma \ref{n-xsb} we can therefore write $Q_{>-2l-5C} F$ as a bounded linear combination of $X^{n/2-1,-1/2,1}$ atoms with frequency $1$ and modulation $\sim 2^{-2l}$.  By Case 2 this contribution is acceptable, so we need only consider the contribution of $Q_{<-2l-5C} F$.

We split
\begin{eqnarray}
 \sgn(t) Q_{<-2l-5C} F = \sum_{\kappa'\in K_{l+C}} &Q^+_{<-2l-3C} P_{0,\kappa'} (\sgn(t) Q_{<-2l-5C} F) \label{2a}\\
+& (1 - Q^+_{< -2l-3C})(\sgn(t) Q_{<-2l-5C} F)\label{2b}
\end{eqnarray}
and consider the two terms separately.

\divider{Case 2(a).3(a).  The $N[0]$ norm of \eqref{2a}.}

By \eqref{f-null} this is bounded by
$$
\lesssim 
(\sum_{\kappa'\in K_{l+C}} \| Q^+_{<-2l-3C} P_{0,\kappa'} (\sgn(t) Q_{<-2l-5C} F) \|_{NFA[\kappa']}^2)^{1/2}.$$
We expand $F$ in terms of $F_\kappa$ and observe that for fixed $\kappa'$ the contribution of $F_\kappa$ vanishes unless $\kappa' \subset \kappa$.  Thus for each $\kappa'$ there are only $O(1)$ $\kappa$ which contribute, so we may estimate the previous by
$$
\lesssim 
(\sum_{\kappa'\in K_{l+C}} \sum_{\kappa \in K_l: \kappa' \subset \kappa} \| P_{0,\kappa'} Q^+_{<-2l-3C} (\sgn(t) Q_{<-2l-5C} F_\kappa) \|_{NFA[\kappa']}^2)^{1/2}.$$
From \eqref{ftpf} we may insert $\tilde P_{0,\kappa}$ in front of $F_\kappa$.  By \eqref{inscribed} we may replace the $NFA[\kappa']$ norm
by the $NFA[\kappa]$ norm.  We then use Lemma \ref{proj-disposable} to discard $P_{0,\kappa'} Q^+_{<-2l-3C}$.  By \eqref{lif} we can then discard $\sgn(t)$.  By another application of Lemma \ref{proj-disposable} we then discard $\tilde P_{0,\kappa} Q_{<-2l-5C}$.  We have thus estimated the previous by
$$
\lesssim 
(\sum_{\kappa'\in K_{l+C}} \sum_{\kappa \in K_l: \kappa' \subset \kappa} \| F_\kappa \|_{NFA[\kappa]}^2)^{1/2}.$$
But this is acceptable by \eqref{chunky}.

\divider{Case 2(a).3(b).  The $N[0]$ norm of \eqref{2b}.}

We may freely replace $\sgn$ by $\sgn_{>-2l-4C}$ since the contribution of $\sgn_{\leq -2l-4C}$ vanishes.  We then estimate this term by
$$ \leq
\| \sgn_{>-2l-4C}(t) Q_{<-2l-5C} F \|_{N[0]}
+
\| 
Q^+_{< -2l-3C}(\sgn_{>-2l-4C}(t) Q_{<-2l-5C} F)
\|_{N[0]}.$$
For the first term we use \eqref{f-l12} and for the second term we use \eqref{f-xsb}.  Since the latter term has Fourier support in the region $D_0 \sim 1, D_- \sim 2^{-2l}$, we can thus bound the previous by
$$ \lesssim
\| \sgn_{>-2l-4C}(t) Q_{<-2l-5C} F \|_{L^1_t L^2_x}
+ 2^{l}
\| 
Q^+_{< -2l-3C}(\sgn_{>-2l-4C}(t) Q_{<-2l-5C} F)
\|_{L^2_t L^2_x}.$$
In the latter term we can use Plancherel to discard $Q^+_{<-2l-3C}$, and then discard the bounded function $\sgn_{>-2l-4C}$.  In the former term we use H\"older.  We can thus bound the previous by
$$ \lesssim
(\| \sgn_{>-2l-4C}\|_2 + 2^l) \| Q_{<-2l-5C} F \|_{L^2_t L^2_x}.$$
However, a routine calculation (as in Case 2) shows that $\| \sgn_{>-2l-4C} \|_2 \lesssim 2^l$, while from Lemma \ref{n-xsb} we have $\| Q_{<-2l-5C} F \|_{L^2_t L^2_x} \lesssim 2^{-l}$.  The claim follows.

\divider{Step 2(b).  Proof of the retarded estimate.}

By Step 2(a) it suffices to show that
$$ \| P_0 \int_{-\infty}^t \frac{\sin((t-s)\sqrt{-\Delta})}{\sqrt{-\Delta}} \eta(t-s) F(s)\ ds \|_{S[0]} \lesssim \| F \|_{N[k]}$$
for all $F \in N[k]$.

We are thus reduced to showing that
\be{splash}
\| \phi \|_{S[0]} \lesssim 1
\end{equation}
whenever
\be{phi-duhamel}
\phi := P_0 \int_\R \eta_T^+(t-s) \frac{\sin((t-s)\sqrt{-\Delta})}{\sqrt{-\Delta}} F(s)\ ds,
\end{equation}
$F$ is an $N[k]$ atom and $\eta_T^+$ is the restriction of $\eta_T$ to $[0,\infty)$.

Observe the identity
\be{fourier-ident}
{\cal F} \phi(\tau,\xi) = C \frac{m(\xi)}{|\xi|} (\hat \eta_T^+(\tau - |\xi|) - \hat \eta_T^+(\tau + |\xi|)) {\cal F}  F(\tau,\xi)
\end{equation}
and the estimate 
\be{est}
\hat \eta_T^+(\tau - |\xi|) - \hat \eta_T^+(\tau + |\xi|) = O(\frac{D_-}{1 + D_-}),
\end{equation}
which follows from a routine computation of $\hat \eta_T^+$ and the hypothesis $T \geq 1$.

From \eqref{fourier-ident} and Lemma \ref{n-xsb} we see that 
$$ \| \nabla_{x,t} \phi \|_{\dot X^{n/2,1/2,\infty}_0} \lesssim \| P_0 F \|_{\dot X^{n/2-1, -1/2,\infty}_0} \lesssim
\| F \|_{\dot X^{n/2-1, -1/2,\infty}_k} \lesssim
 \| F \|_{N[k](\R^n)} \lesssim 1.$$
Thus the second component of the $S[0]$ norm in \eqref{sk-def} always makes an acceptable contribution to \eqref{splash}.

We now divide into three cases depending on which type of $N[k]$ atom $F$ is.

\divider{Case 2(b).2(a). $F$ is an $\dot X^{n/2-1,-1/2,1}$ atom with frequency $2^{k}$ and modulation $2^j$.}

In this case $F$ has Fourier support in $2^{j-5} \leq D_- \leq 2^{j+5}$ and
$$ \| F \|_{L^2_t L^2_x} \lesssim 2^{j/2}.$$
From this, \eqref{est} and \eqref{fourier-ident} we see that $\phi$ is bounded in $\dot X^{n/2,1/2,1}_0$, and \eqref{splash} follows from Lemma \ref{xn11 contain}.

\divider{Case 2(b).2(b). $F$ is an $L^1_t \dot H^{n/2-1}_x$ atom with frequency $2^k$.}

By Minkowski's inequality and time translation invariance it suffices to show that
$$
\| P_0  \eta_T^+(t) \frac{\sin(t\sqrt{-\Delta})}{\sqrt{-\Delta}} f \|_{S[0]} \lesssim \|f\|_2$$
for all $f \in L^2(\R^n)$.

The $L^\infty_t L^2_x$ component of \eqref{sk-def} is acceptable by standard energy estimates.  The $\dot X^{n/2,1/2,\infty}_0$ component is acceptable by the remarks made previously, so we turn to the last component.
The idea is to treat the expression inside the norm as a free solution, smoothly truncated to the upper half-space $\{t>0\}$.  

By conjugation symmetry it suffices to show
$$
(\sum_{\kappa \in K_l}
\| P_{0,\kappa} Q^+_{<-2l} P_0  \eta_T^+(t) \frac{\sin(t\sqrt{-\Delta})}{\sqrt{-\Delta}} f \|_{S[0,\kappa]}^2)^{1/2}
\lesssim \|f\|_2.
$$
for all $l > 10$.  

Fix $l$.  By the finite overlap of the $P_{0,\kappa}$ it suffices to show that
\be{chook}
\| P_{0,\kappa} Q^+_{<-2l} P_0  \eta_T^+(t) \frac{\sin(t\sqrt{-\Delta})}{\sqrt{-\Delta}} f \|_{S[0,\kappa]}
\lesssim \|f\|_2
\end{equation}
for each $\kappa \in K_l$ and all $f \in L^2$.

Fix $\kappa$, $f$.  The expression inside the norm has a spacetime Fourier transform of
$$ C m_{0,\kappa}(\xi) m_0(2^{2l} (\tau - |\xi|)) \frac{m(\xi)}{|\xi|} (\hat \eta_T^+(\tau - |\xi|) + O(1)) \hat f(\xi)$$
since $\hat \eta_T^+(\tau + |\xi|) = O(1)$ on the support of $m_0(2^{2l} (\tau - |\xi|)) \frac{m(\xi)}{|\xi|}$.

We treat the main term $\eta_T^+(\tau - |\xi|)$ and the error $O(1)$ separately.

\divider{Case 2(b).2(b).1. The contribution of $O(1)$.}

By Lemma \ref{f-lemma} we may estimate this contribution by
$$
\lesssim 
\| {\cal F}^{-1} O(m_{0,\kappa}(\xi) m(2^{2l} (\tau - |\xi|)) \frac{m(\xi)}{|\xi|} \hat f(\xi)) \|_{\dot X^{n/2,1/2,1}_1},$$
which by Plancherel and the support of $m(2^{2l} (\tau - |\xi|))$ is bounded by
$$
\lesssim 2^{-l}
\| m_{0,\kappa}(\xi) m(2^{2l} (\tau - |\xi|)) \frac{m(\xi)}{|\xi|} \hat f(\xi) \|_{L^2_\xi L^2_\tau},$$
which is acceptable (with a factor of $2^{-2l}$ to spare).

\divider{Case 2(b).2(b).2. The contribution of $\hat \eta_T^+(\tau - |\xi|)$.}

This contribution can be rewritten as
$$
C \| [\eta_T^+ * \check m(2^{-2l}\cdot)](t) P_{0,\kappa} P_0  \frac{e^{it\sqrt{-\Delta}}}{\sqrt{-\Delta}} f \|_{S[0,\kappa]}.
$$
The expression inside the norm is Schwartz, so the above expression can be estimated by
$$
\lesssim 
\| [\eta_T^+ * \check m(2^{-2l}\cdot)](t) \eta_M(t) P_{0,\kappa} P_0  \frac{e^{it\sqrt{-\Delta}}}{\sqrt{-\Delta}} f \|_{S[0,\kappa]}^2
$$
for some sufficiently large $M > 0$ (depending on $f$, $l$, $T$).
By \eqref{liff} we can discard the bounded function $\eta_T^+ * \check m_0(2^{-2l}\cdot)$. We then revert to the frequency space formulation, and apply Lemma \ref{f-lemma} to bound this contribution by
$$
\lesssim \| {\cal F}^{-1} [m_{0,\kappa}(\xi) \frac{m(\xi)}{|\xi|} \hat \eta_M(\tau - |\xi|) \hat f(\xi)] \|_{\dot X^{n/2,1/2,1}_1},$$
which is acceptable by Plancherel and a direct computation.

\divider{Case 2(b).2(c). $F$ is a $\pm$-null frame atom with frequency $2^{k''}$ and angle $2^{-l}$.}

By conjugation symmetry we may take $\pm = +$.  Then there exists a decomposition
$F = \sum_{\kappa \in K_l} F_\kappa$ such that each $F_\kappa$ has Fourier support in 
$
\{ \tau > 0; D_- \lesssim 2^{-2l}; D_0 \sim 1; \Theta \in \frac{1}{2} \kappa \}
$
and
\be{chopper}
(\sum_{\kappa \in K_l} \| F_\kappa \|_{NFA[\kappa]}^2)^{1/2} \lesssim 1.
\end{equation}

From \eqref{fourier-ident} we see that $\phi$ has the same Fourier support as $F$.  Since $F$ is supported on the region $D_0 \sim 1, D_- \lesssim 2^{-2l}$, we have $D_+ \sim 1$ and therefore that
$$ \| \nabla_{x,t} \phi \|_{L^\infty_t L^2_x} \lesssim \| \phi \|_{L^\infty_t L^2_x}.$$
By orthogonality and \eqref{sklw-def} we thus see that the first component of $\|\phi \|_{S[0]}$ in \eqref{sk-def} is bounded by the third.  

The second component is controlled by the remarks made previously, so it remains to control the third component of $\| \phi \|_{S[0]}$, i.e. that
\be{third}
 (\sum_{\kappa' \in K_{l'}} \| P_{0,\pm\kappa'} Q^\pm_{<-2l'} \phi  \|_{S[0,\kappa']}^2)^{1/2}
\lesssim 1
\end{equation}
for all $l' > 10$ and choices of sign $\pm$.  Actually, we can take $\pm = +$ since the $\pm = -$ term vanishes completely.

We can decompose $\phi = \sum_{\kappa \in K_l} \phi_{\kappa}$, where $\phi_{\kappa}$ is given by the Fourier representation
\be{phiw}
{\cal F} \phi_{\kappa} := C \frac{m(\xi)}{|\xi|} (\hat \eta_T^+(\tau - |\xi|) - \hat \eta_T^+(\tau + |\xi|)) {\cal F} F_{\kappa}(\tau,\xi).
\end{equation}

We shall first show \eqref{third} for $l' \geq l+C$, and then later bootstrap this to the case $l+C > l' > 10$.  

\divider{Case 2(b).2(c).1.  Proof of \eqref{third} when $l' \geq l+C$.}

In this case each $\phi_{\kappa}$ only contributes to those $\kappa' \in K_{l'}$ for which $\kappa' \subset \kappa$.  Thus each $\kappa' \in K_{l'}$ has contributions from at most $O(1)$ values of $\kappa \in K_l$.  It therefore suffices by \eqref{chopper} to show that
\be{cheetah}
 (\sum_{\kappa' \in K_{l'}: \kappa' \subset \kappa} \| P_{0,\kappa'} Q^+_{<-2l'} \phi_{\kappa}  \|_{S[0,\kappa']}^2)^{1/2}
\lesssim \| F_{\kappa} \|_{NFA[\kappa]}
\end{equation}
for all $\kappa \in K_l$.

Fix $\kappa$.  Because of the multiplier $Q^+_{<-2l'}$ we may estimate $\hat \eta^+_T(\tau + |\xi|)$ by $O(1)$ in \eqref{phiw}, as in Case 2(b).2(b).  From Definition \ref{nfa-def} and \eqref{phiw} it suffices to show
\begin{eqnarray*}
 (\sum_{\kappa' \in K_{l'}: \kappa' \subset \kappa } \| &P_{0,\kappa'} Q^+_{<-2l'} {\cal F}^{-1} [ \frac{m(\xi)}{|\xi|}(\hat \eta_T^+(\tau - |\xi|) + O(1)) {\cal F}  F(\tau,\xi) ]
  \|_{S[0,\kappa']}^2)^{1/2}\\
&\lesssim \frac{1}{\dist(\omega, \kappa)} \| F \|_{L^1_{t_{\omega}} L^2_{x_{\omega}} }
\end{eqnarray*}
for all $\omega \not \in 2\kappa$ and $F \in L^1_{t_{\omega}} L^2_{x_{\omega}}$.

Fix $\omega$.  By rotation invariance we may take $\omega = e_1$.
By Minkowski's inequality it suffices to prove this estimate with $F$ replaced by $\delta(t_{e_1} - t_0) f(x_{e_1})$ for some $t_0 \in \R$ and $f \in L^2(NP(e_1))$, and with 
$\| F \|_{L^1_{t_{e_1}} L^2_{x_{e_1}} }$ replaced by $\|f\|_2$.  By translation invariance we may let $t_0 = 0$, so we reduce to
$$
 (\sum_{\kappa' \in K_{l'}: \kappa' \subset \kappa } \| P_{0,\kappa'} Q^+_{<-2l'} {\cal F}^{-1} [ \frac{m(\xi)}{|\xi|} (\hat \eta_T^+(\tau - |\xi|) + O(1)) \hat f(\xi_{e_1}) ]
  \|_{S[0,\kappa']}^2)^{1/2}
\lesssim \frac{1}{\dist(e_1,\kappa)} \| f \|_2.
$$
For each $\kappa'$, the only portion of $\hat f$ which contributes is that which is supported on the region 
\be{sigma}
\{ \xi_{e_1}: \tau > 0, D_+ \sim 1, D_- \lesssim 2^{-2l'}, \Theta \in \frac{1}{2} \kappa' \}.
\end{equation}
From elementary geometry we see that these regions are finitely overlapping as $\kappa'$ varies.  Thus by Plancherel we need only show that
$$
\| P_{0,\kappa'} Q^+_{<-2l'} {\cal F}^{-1} [ \frac{m(\xi)}{|\xi|} (\hat \eta_T^+(\tau - |\xi|) + O(1)) \hat f(\xi_{e_1}) ]
  \|_{S[0,\kappa']}
\lesssim \frac{1}{\dist(e_1,\kappa')} \| f \|_2.
$$
for all $\kappa' \in K_{l'}$ such that $\kappa' \subset \kappa$, and all $f \in L^2$ (using the trivial estimate $\dist(e_1,\kappa') \geq \dist(e_1,\kappa)$).

Fix $\kappa'$. We rewrite the left-hand side of the previous as
\be{fillip-2}
\|
{\cal F}^{-1} [ 
m_{0,\kappa'}(\xi)
m_0(2^{2l'}(\tau-|\xi|))
\frac{m(\xi)}{|\xi|}
[\hat \eta_T^+(\tau - |\xi|) + O(1)] \hat f(\xi_{e_1}) ]\|_{S[0,\kappa']}.
\end{equation}
We shall now replace these Euclidean multipliers with null frame counterparts.

Observe that we may freely insert the multiplier $\tilde m_{0,\kappa'}(\xi_{e_1})$, where $\tilde m_{0,\kappa'}$ equals 1 on the region
$$ \{ \xi_{e_1}: \tau > 0; D_0 \sim 1, D_- \lesssim 2^{-2l'}, \Theta \in \kappa \}$$
and is a bump function adapted to a dilate of this region.

Next, we observe the identity
\be{tau-ident}
\tau - |\xi| = \frac{1}{\tau+|\xi|} (\tau^2 - |\xi|^2)
= \frac{2 \xi^1_{e_1}}{\tau+|\xi|} (\tau_{e_1} - h(\xi_{e_1}))
\end{equation}
where
$$ \xi^1_{e_1} := \xi_{e_1} \cdot \frac{1}{\sqrt{2}}(1,-e_1); \quad
\xi'_{e_1} := \xi_{e_1} - \xi^1_{e_1} \frac{1}{\sqrt{2}} (1, -e_1)$$
and $h$ is the function
$$ h(\xi_{e_1}) := |\xi'_{e_1}|^2 / 2 \xi^1_{e_1}.$$
On the Fourier support of \eqref{fillip-2} we see from elementary geometry that $|\xi'_{e_1}| \sim \dist(e_1,\kappa')$ and $|\xi^1_{e_1}| \sim \dist(e_1,\kappa')^2$.  We thus have
\be{taugraph}
\tau_{e_1} = h(\xi_{e_1}) + O(2^{-2l'}/\dist(e_1,\kappa')^2).
\end{equation}
We can therefore freely insert the multiplier
$$ m_0(2^{2l'-C} \dist(e_1,\kappa')^2 (\tau_{e_1} - h(\xi_{e_1}))$$
in \eqref{fillip-2}.

By Lemma \ref{proj-disposable} the multiplier with symbol $m_{0,\kappa'}(\xi) 
m_0(2^{2l'}(\tau-|\xi|)) m(\xi)/|\xi|$ is disposable.  We may therefore estimate \eqref{fillip-2} by
$$ \lesssim 
\|
{\cal F}^{-1} [ 
m_0(2^{2l'-C} \dist(e_1,\kappa')^2 (\tau_{e_1}-h(\xi_{e_1})))
[\hat \eta_T^+(\tau - |\xi|) + O(1)]
\tilde m_{0,\kappa'}(\xi_{e_1}) \hat f(\xi_{e_1}) ]\|_{S[0,\kappa']}.$$
It still remains to convert $\hat \eta_T^+(\tau - |\xi|)$ to a null frame equivalent.

A direct computation gives the estimate 
$$ \hat \eta_T^+(\tau) = \frac{C}{\tau} + O( T(T\tau)^{-100} )$$
for $|\tau| \gtrsim T^{-1}$, and
$$ \hat \eta_T^+(\tau) = O(T)$$
for $|\tau| \lesssim T^{-1}$.  Combining these facts with \eqref{tau-ident} we have
$$ \hat \eta_T^+(\tau - |\xi|) + O(1) = \frac{\tau+|\xi|}{2\xi^1_{e_1}} 
\hat \eta_{T/\dist(e_1,\kappa')^2}^+(\tau_\theta - h(\xi_{e_1}))
+ O(T \min(1, (T D_-)^{-100}) + O(1).$$

We deal with the main term and the two $O()$ terms separately.

\divider{Case 2(b).2(c).1(a). The contribution of the $O()$ terms.}

We use \eqref{into-sklw} to estimate this contribution by
\begin{eqnarray*}
\lesssim \|&
{\cal F}^{-1} [ 
m_0(2^{2l'-C} \dist(e_1,\kappa')^2 (\tau_{e_1}-h(\xi_{e_1})))
[O(T \min(1, (T D_-)^{-100}) + O(1)] \\
&\tilde m_{0,\kappa'}(\xi_{e_1}) \hat f(\xi_{e_1}) ]\|_{\dot X^{n/2,1/2,1}_0}.
\end{eqnarray*}
The expression inside the norm has Fourier support on $\tau = |\xi| + O(2^{-2l'})$.  From \eqref{tau-ident}, Plancherel and the triangle inequality we can thus estimate the previous by
$$ \lesssim \sum_{j \leq -2l'+C} 2^{j/2}
[T \min(1, (T 2^j)^{-100}) + 1]
\|
{\cal F}^{-1} [ 
\chi_{\tau = |\xi| + O(2^{j})}
\tilde m_{0,\kappa'}(\xi_{e_1}) \hat f(\xi_{e_1}) ]\|_{L^2_t L^2_x}.$$
By Plancherel we may estimate the previous by
$$
\lesssim \sum_{j \leq -2l'+C} 2^{j/2}
[T \min(1, (T 2^j)^{-100}) + 1]
\|
\chi_{\tau = |\xi| + O(2^j)} \tilde m_{0,\kappa'}(\xi_{e_1}) \hat f(\xi_{e_1})
\|_{L^2_{\xi_{e_1}} L^2_{\tau_{e_1}} }.$$
For each $\xi_{e_1}$ in the support of $\tilde m_{0,\kappa'}$, the expression $\chi_{\tau = |\xi| + O(2^j)}$ has an $L^2_{\tau_{e_1}}$ norm of $2^{j/2}/\dist(e_1,\kappa')$ by \eqref{tau-ident}.  Thus we can estimate the previous by
$$
\lesssim \sum_{j \leq -2l'+C} \frac{2^j}{\dist(e_1,\kappa')}
[T \min(1, (T 2^j)^{-100}) + 1]
\|f\|_2$$
which is acceptable.

\divider{Case 2(b).2(c).1(b). The contribution of the main term.}

We need to show 
\begin{eqnarray*}
\|
{\cal F}^{-1} [ &
m_0(2^{2l'-C} \dist(e_1,\kappa')^2 (\tau_{e_1}-h(\xi_{e_1})))
\frac{\tau+|\xi|}{2\xi^1_{e_1}}
\hat \eta_{T/\dist(e_1,\kappa')^2}^+(\tau_{e_1} - h(\xi_{e_1})) \\
&\tilde m_{0,\kappa'}(\xi_{e_1}) \hat f(\xi_{e_1}) ]\|_{S[0,\kappa']}
\lesssim \frac{1}{\dist(e_1,\kappa')} \| f \|_2.
\end{eqnarray*}
Since the above expression has Fourier support on the region $|\tau| + |\xi| \sim 1$, the multiplier $(\tau + |\xi|)/2$ is effectively disposable and can be discarded.  We can thus estimate the previous by
$$\lesssim 
\|
{\cal F}^{-1} [ 
(\varphi \hat \eta_{T/\dist(e_1,\kappa')^2}^+)(\tau_{e_1} - h(\xi_{e_1})) \hat f_{\kappa'}(\xi_{e_1})
]\|_{S[0,\kappa']}$$
where
$$ \hat f_{\kappa'}(\xi_{e_1}) := \frac{1}{\xi^1_{e_1}} \tilde m_{0,\kappa'}(\xi_{e_1}) \hat f(\xi_{e_1})$$
and
$$ \varphi(\tau) := m_0(2^{2l'-C} \dist(e_1,\kappa')^2 \tau).$$
We can rewrite this further as
$$\lesssim 
\| (\check \varphi * \eta_{T/\dist(e_1,\kappa')^2}^+)(t_{e_1}) \eta_M
{\cal F}^{-1} [ 
\delta(\tau_{e_1} - h(\xi_{e_1})) \hat f_{\kappa'}(\xi_{e_1})
]\|_{S[0,\kappa']}$$
for some large $M$ (depending on $f$, $l$, $T$).  

By \eqref{liff} we may discard the $(\check \varphi * \eta_{T/\dist(e_1,\kappa')^2}^+)(t_{e_1})$ factor.  We then move $\eta_M$ inside the Fourier transform and use \eqref{into-sklw} to estimate the previous by
$$\lesssim 
\| {\cal F}^{-1} [ 
\hat \eta_M(\tau_{e_1} - h(\xi_{e_1})) \hat f_{\kappa'}(\xi_{e_1})
]\|_{\dot X^{n/2,1/2,1}_0}.$$
By Plancherel we can estimate this by
$$\lesssim \sum_{j} 2^{j/2}
\|  
\chi_{\tau = |\xi| + O(2^j)}
\eta_M(\tau_{e_1} - h(\xi_{e_1})) |\xi^1_{e_1}|^{-1}
m_{0,\kappa'}(\xi_{e_1}) \hat f(\xi_{e_1})
\|_{L^2_{\xi_{e_1}} L^2_{\tau_{e_1}} }.$$
We use the estimates $|\tau_{e_1} - h(\xi_{e_1})| \sim 2^j / \dist(e_1,\kappa')^2$ (from \eqref{tau-ident}) and $\hat \eta_M(\tau) = O(M \min(1, (M|\tau|)^{-100}))$ to estimate this by
\begin{eqnarray*}
\lesssim \sum_{j}& 2^{j/2} M \min(1, 
(M 2^j / \dist(e_1,\kappa')^2)^{-100})\\
&\|  
\chi_{\tau_{e_1} = h(\xi_{e_1}) = 2^j / \dist(e_1,\kappa')^2 }
|\xi^1_{e_1}|^{-1}
m_{0,\kappa'}(\xi_{e_1}) \hat f(\xi_{e_1})
\|_{L^2_{\xi_{e_1}} L^2_{\tau_{e_1}} }.
\end{eqnarray*}
Since $|\xi^1_{e_1}| \sim \dist(e_1,\kappa')^2$, we can estimate this by
$$\lesssim \dist(e_1,\kappa')^{-1} \sum_j \frac{2^j M}{\dist(e_1,\kappa')^2} \min(1, 
(M 2^j / \dist(e_1,\kappa')^2)^{-100})
\| \hat f\|_2$$
which is acceptable.  This completes the proof of \eqref{third} when $l' > l+C$.

\divider{Case 2(b).2(c).2.  Proof of \eqref{third} when $10 < l' < l+C$}

We divide $Q_{<-2l'} \phi = Q_{<-2l-4C} \phi + Q_{-2l-4C \leq \cdot < -2l} \phi$.

\divider{Case 2(b).2(c).2(a).  The contribution of $Q_{-2l-4C \leq \cdot < -2l} \phi$.}

In this case we use \eqref{into-sklw} to bound this contribution by
$$
\lesssim (\sum_{\kappa' \in K_{l'}} \| P_{0,\kappa'} Q^+_{-2l-4C \leq \cdot <-2l'} \phi  \|_{\dot X^{n/2,1/2,1}_0}^2)^{1/2}.$$
From \eqref{fourier-ident} we see that the expression inside the norm has Fourier support in $D_- \sim 2^{-2l}$, hence we can estimate the previous by
$$
\lesssim 2^{-l} (\sum_{\kappa' \in K_{l'}} \| P_{0,\kappa'} Q^+_{-2l-4C \leq \cdot <-2l'} \phi  \|_{L^2_t L^2_x}^2)^{1/2},$$
which by Plancherel is bounded by
$$
\lesssim 2^{-l} \| Q_{-2l-4C \leq \cdot <-2l'} \phi  \|_{L^2_t L^2_x}^2)^{1/2},$$
which is acceptable since we have already bounded the second component of \eqref{sk-def}.

\divider{Case 2(b).2(c).2(b).  The contribution of $Q_{<-2l-4C} \phi$.}

In this case we write the contribution to \eqref{third} as
$$
(\sum_{\kappa' \in K_{l'}} \| \sum_{\kappa'' \in K_{l+2C}}
P_{0,\kappa'} P_{0,\kappa''} Q^+_{<-2l-4C} \phi  \|_{S[0,\kappa']}^2)^{1/2}.
$$
We may restrict the inner summation to the case when $\kappa'' \subset \kappa'$ since the contribution vanishes otherwise.  We then discard $P_{0,\kappa'}$, estimating the previous by
$$
\lesssim
(\sum_{\kappa' \in K_{l'}} \| \sum_{\kappa'' \in K_{l+2C}: \kappa'' \subset \kappa'}
P_{0,\kappa''} Q^+_{<-2l-4C} \phi  \|_{S[0,\kappa']}^2)^{1/2}.
$$
By \eqref{sklw-ortho} we may estimate this by
$$
\lesssim 
(\sum_{\kappa' \in K_{l'}} \sum_{\kappa \in K_{l+2C}: \kappa \subset \kappa'}
\| P_{0,\kappa''} Q^+_{<-2l-4C} \phi  \|_{S[0,\kappa'']}^2)^{1/2}
$$
which simplifies to
$$
\lesssim 
(\sum_{\kappa \in K_{l+2C}}
\| P_{0,\kappa''} Q^+_{<-2l-4C} \phi  \|_{S[0,\kappa'']}^2)^{1/2}
$$

Split $\phi = \sum_{\kappa \in K_l} \phi_{\kappa}$ as before.  For each $\kappa''$ we can restrict the $\kappa$ summation to those caps for which $\kappa'' \subset \kappa$, so that only $O(1)$ values of $\kappa$ contribute for each $\kappa''$.  Thus we may estimate the previous by
$$
\lesssim 2^{-(n-1)l'/2} 2^{(n-1)(l'-l)/2}
(\sum_{\kappa'' \in K_{l+2C}} \sum_{\kappa \in K_l: \kappa'' \subset \kappa}
\| P_{0,\kappa''} Q^+_{<-2l-4C} \phi_{\kappa}  \|_{S[0,\kappa'']}^2)^{1/2}.
$$
By \eqref{cheetah} (applied with $l'$ replaced by $l+2C$) and \eqref{sklw-def} we can bound this by
$$
\lesssim 2^{-(n-1)l'/2} 2^{(n-1)(l'-l)/2} 2^{(n-1)l/2}
(\sum_{\kappa'' \in K_{l+2C}} \sum_{\kappa \in K_l: \kappa'' \subset \kappa }
\| F_{\kappa} \|_{NFA[\kappa'']}.
$$
By \eqref{inscribed} we may replace $NFA[\kappa'']$ with $NFA[\kappa]$.  But this is then acceptable by \eqref{chopper}.
This concludes the proof of \eqref{energy-est-2}.
\endprf

\section{The continuity of \eqref{upstairs}}\label{continuity-sec}

In this section we prove the continuity of \eqref{upstairs}.  Our main tools will be the estimates \eqref{energy-est-2}, \eqref{l12}, which have already been proven. We remark that this section can be read independently since the arguments used here are not needed elsewhere in the paper.

Our arguments here are somewhat inelegant, but we have not been able to find a more natural approach for this problem.

Fix $\phi$, $T$, $T_0$; we may extend $\phi$ to all of $\R^{1+n}$ by the free wave equation outside of $[-T_0,T_0] \times \R^n$.  

Since the spaces $S_k([-T,T] \times \R^n)$ are defined by restriction we see that the function $a(T)$ is monotone non-decreasing in $T$.  

Let $A$ denote the quantity $A := \liminf_{T' \to T} a(T')$.  Clearly $A \geq 1$.  By monotonicity we need to show that
\be{chips}
\| \phi_k|_{[-T-\epsilon,T+\epsilon] \times \R^n}  \|_{S_k([-T-\epsilon,T+\epsilon] \times \R^n)} \lesssim A c_k
\end{equation}
for all $k$, where $0 < \epsilon \ll 1$ is some quantity depending on $\phi$, $T$, $A$ but is independent of $k$.

By \eqref{energy-est-2} we may majorize the left-hand side of \eqref{chips} by
\be{chips-triv}
2^{nk/2} \| \phi_k(0) \|_2 + 2^{(\frac{n}{2}-1)k} \| \partial_t \phi_k(0) \|_2 + 2^{(\frac{n}{2}-1)k} \| \Box \phi_k \|_{L^1_t L^2_x([-T-\epsilon,T+\epsilon] \times \R^n)}.
\end{equation}
Since $\phi$ is a classical wave map, we see that this quantity decays like $O(2^{-|k|})$ or faster as $k \to \pm \infty$, uniformly in $T' \in [0,T_0]$.  From \eqref{local} we thus see that \eqref{chips} holds whenever $|k|$ is sufficiently large.  Thus we only have a finite number of $k$ left to deal with, which implies that we only need to show \eqref{chips} for each $k$ separately (with $\epsilon$ now allowed to depend on $k$).

Fix $k$; we may rescale $k=0$.  We first dispose of the case $T=0$.  In this case we must have $A = a(0)$.  By \eqref{energy-est} we then have
$$ \| \phi_0(0) \|_2 + \| \partial_t \phi_0(0) \|_2 \lesssim A c_0.$$
Since $\phi$ is a classical wave map, we may therefore bound \eqref{chips-triv} by $\lesssim A c_0$ for sufficiently small $\epsilon$, as desired.

Now suppose $T > 0$.  If $\epsilon$ is sufficiently small, we have 
$$ a(T - \epsilon) \sim A$$
by monotonicity.  We may therefore find $\tilde \phi_0 \in S_0(\R^{1+n})$ which agrees with $\phi_0$ on $[-T+\epsilon, T-\epsilon]$ and satisfies the estimate
\be{agog}
\| \tilde \phi_0 \|_{S_0} \lesssim A c_0.
\end{equation}
By replacing $\tilde \phi_0$ with $P_{-5 < \cdot < 5} \tilde \phi_0$ if necessary we may assume that $\tilde \phi_0$ has Fourier support on the region $D_0 \sim 1$.

From \eqref{energy-est} and the previous we see that
\be{milk}
\| \tilde \phi_0(t) \|_2 \lesssim \| \nabla_{x,t} \tilde \phi[t] \|_2 \lesssim A c_0
\end{equation}
for all $t \in \R$.  Unfortunately, this is not enough regularity for us to apply \eqref{energy-est-2}, since we need to control $\Box \tilde \phi$.  To resolve this we shall regularize $\tilde \phi$.  Specifically, define the smoothing operator $S$ by
$$ S \tilde \phi(t) := \int  \tilde \phi(t+\epsilon s) \varphi(s)\ ds$$
where $\varphi$ is a bump function on $[-1,1]$ of mass 1.  From Minkowski's inequality and \eqref{agog} we have
$$ \| S \tilde \phi_0 \|_{S_0} \lesssim A c_0,$$
which of course implies
$$ \| S \tilde \phi_0|_{[-T-\epsilon,T+\epsilon] \times \R^n}  \|_{S_0([-T-\epsilon,T+\epsilon] \times \R^n)} \lesssim A c_0.$$
From the identity
$$ \phi_0 = S \tilde \phi_0 + S (\phi_0 - \tilde \phi_0) + (1-S) \phi_0$$
and \eqref{energy-est} we see that \eqref{chips} will follow from
\begin{eqnarray}\label{chips-2} 
&\| \Box S(\phi_0 - \tilde \phi_0) \|_{L^1_t L^2_x([-T-\epsilon,T+\epsilon])}
+ \| (1-S) \phi_0[0] \|_{L^2 \times L^2}\\
& + \| (1-S) \Box \phi_0 \|_{L^1_t L^2_x([-T-\epsilon,T+\epsilon] \times \R^n)}\nonumber
\lesssim A c_0
\end{eqnarray}
(note that $S(\phi_0 - \tilde \phi_0)$ vanishes at time 0).

Since $\phi$ is a classical wave map on $[-T_0,T_0]$ and is extended by the free wave equation, we have the estimates
\be{cheese}
\| \nabla_{x,t}^j \phi_0(t) \|_{L^2} \leq C_{\phi}
\end{equation}
for all $t \in \R$ and $j=0,1,2$, and some quantity $C_{\phi} < \infty$ which is independent of $\epsilon$.  From this and Lebesgue differentiation theorem we see that the last two terms of \eqref{chips-2} go to zero as $\epsilon \to 0$, and are therefore acceptable if $\epsilon$ is sufficiently small.

Since $S(\phi_0 - \tilde \phi_0)$ vanishes on $[-T+2\epsilon, T-2\epsilon]$ we may bound the first term of \eqref{chips-2} by
$$
\lesssim \epsilon
\| \Box S(\phi_0 - \tilde \phi_0) \|_{L^\infty_t L^2_x(\{ |t| = T + O(\epsilon) \})}.$$
Expanding out $\Box$ and $S(c)$ and using the fact that $\phi_0 - \tilde \phi_0$ has Fourier support in $D_0 \sim 1$, we may bound this by
$$
\lesssim 
\| \nabla_{x,t}(\phi_0 - \tilde \phi_0) \|_{L^\infty_t L^2_x(\{ |t| = T + O(\epsilon) \})}.$$
The contribution of $\tilde \phi_0$ is acceptable by \eqref{milk}.  To control the contribution of $\phi_0$, we observe from \eqref{cheese} and the Fundamental Theorem of Calculus that
$$ \| \nabla_{x,t} \phi_0(t) \|_{L^2} \lesssim \| \nabla_{x,t} \phi_0(T-\epsilon) \| + C_\phi \epsilon$$
whenever $t = T + O(\epsilon)$.  Since $\phi_0(T-\epsilon) = \tilde \phi_0(T-\epsilon)$, we thus see from \eqref{milk} that this term is also acceptable if $\epsilon$ is sufficiently small.  This concludes the proof of \eqref{chips}.

\section{The geometry of the cone}\label{geometry-sec}

In the next few sections we shall be proving a number of bilinear estimates.  In all of these estimates it will be important to understand the relationship between the modulation of $\phi$, $\psi$, and $\phi \psi$, and with the angular separation of $\phi$ and $\psi$. 

\begin{definition}\label{unilateral-def}  
If $N_0, N_1, N_2$ are numbers, we write $N_{max} \geq N_{med} \geq N_{min}$ for the maximum, median, and minimum of $N_0$, $N_1$, $N_2$ respectively.  We say that $N_0, N_1, N_2$ are \emph{balanced} if $N_{med} \sim N_{max}$ and
\emph{imbalanced} otherwise.
\end{definition}

The Littlewood-Paley  trichotomy thus says that $P_{k_0}(\phi_{k_1} \psi_{k_2})$ vanishes unless the frequencies $2^{k_0}, 2^{k_1}, 2^{k_2}$ are balanced.  There is a similar relationship concerning modulations, but is more complicated.  For simplicity we restrict ourselves to the case when all frequencies lie near the cone:

\begin{lemma}\label{geometry-lemma}
Let $C$ be a large constant, let $2^{k_0}, 2^{k_1}, 2^{k_2}$ be balanced, let $j_0, j_1, j_2$ be integers such that $j_i \leq k_i - C$ for $i=0,1,2$, and let $\kappa_1$, $\kappa_2$ be caps of radius $0 < r \leq 2^{-5}$.  Let $\phi$ have Fourier support on the region $\{ D_0 \sim 2^{k_1}, D_- \sim 2^{j_1}, \Theta \in \kappa_1\}$ and $\psi$ have Fourier support on the region $\{ D_0 \sim 2^{k_2}, D_- \sim 2^{j_2}, \Theta \in \kappa_2\}$.  

\begin{itemize}
\item  If $2^{j_0}, 2^{j_1}, 2^{j_2}$ are imbalanced, then $P_{k_0} Q_{j_0} L( \phi, \psi )$ vanishes unless 
\be{pm-c}
2^{j_{max}} \lesssim 2^{k_{min}} \hbox{ and } \dist(\kappa_1,\kappa_2) + r \sim 2^{k_0 - k_{max}}  2^{(j_{max} - k_{min})/2} + r
\end{equation}

\item If $2^{j_0}, 2^{j_1}, 2^{j_2}$ are balanced, then $P_{k_0} Q_{j_0} L( \phi, \psi )$ vanishes unless 
\be{imbal}
\dist(\kappa_1,\kappa_2) \lesssim 2^{k_0 - k_{max}}  2^{(j_{max} - k_{min})/2} + r.
\end{equation}
\end{itemize}
\end{lemma}

\begin{proof}
If $P_k Q_j L(\phi,\psi)$ does not vanish, then by the Fourier transform there must exist $(\tau_0, \xi_0) = (\tau_1, \xi_1) + (\tau_2, \xi_2)$ such that
$||\tau_i| - |\xi_i|| \sim 2^{j_i}$ and $|\xi_i| \sim 2^{k_i}$ for $i=0,1,2$, and that $\tau_i \xi_i / |\tau_i \xi_i| \in \kappa_i$ for $i=1,2$.  Similarly for $P_k Q_j L(\phi_{,\alpha},\psi^{,\alpha})$. From our hypothesis $j_i \leq k_i - C$ we see that $|\tau_i| \sim |\xi_i|$ for $i=0,1,2$.

By conjugation symmetry we may assume $\tau_0 \geq 0$.  By symmetry it suffices to consider three cases.

\divider{Case 1. ($(++)$ case)  $k_0 \leq k_{max}-C$ and $\tau_1, \tau_2 > 0$.}

Observe the identity
$$ |\xi_1| + |\xi_2| - |\xi_1 + \xi_2| = (|\tau_0| - |\xi_0|) - (|\tau_1| - |\xi_1|) - (|\tau_2| - |\xi_2|).$$
The left hand side is $\sim 2^{k_{max}}$, but the left hand side can at most be $O(2^{j_{max}})$.  Since we are assuming $j_i \leq k_i - C$, this case is therefore impossible.

\divider{Case 2. ($(+-)$ case)  $k_0 \leq k_{max}-C$ and $\tau_1> 0$, $\tau_2 < 0$.}

This forces $k_0 = k_{min}$ and $k_1, k_2 = k_{max} + O(1)$.

Now observe the identity
$$ -|\xi_1| + |\xi_2| + |\xi_1 + \xi_2| = -(\tau_0 - |\xi_0|) + (\tau_1 - |\xi_1|) - (|\tau_2| - |\xi_2|).$$
The right-hand side has magnitude $O(2^{j_{max}})$, and in the imbalanced case has magnitude $\sim 2^{j_{max}}$.  The left-hand side can be rewritten as
$$ -|\xi_1| + |\xi_2| + |\xi_1 + \xi_2| = 2 \frac{|\xi_1 + \xi_2| |\xi_2| + (\xi_1 + \xi_2) \cdot \xi_2}{|\xi_1| + |\xi_2| + |\xi_1 + \xi_2|}$$
and therefore
$$ -|\xi_1| + |\xi_2| + |\xi_1 + \xi_2| \sim \angle(\xi_1 + \xi_2, -\xi_2)^2 2^{k_0}.$$
By the sine rule we thus have
$$ -|\xi_1| + |\xi_2| + |\xi_1 + \xi_2| \sim \angle(\xi_1, -\xi_2)^2 2^{-k_0} 2^{-k_{max}}.$$
Combining this with the previous we see that we are in either \eqref{pm-c} or \eqref{imbal}.

\divider{Case 3. (Low-high case)  $k_0 = k_{max}+O(1)$ and $k_2 = k_{max} + O(1)$.}

If $\tau_2 < 0$, then $\tau_1 > 0$.  Since $|\tau_i| \sim |\xi_i|$, we must then have $k_1 = k_{max} + O(1)$, at which point we could swap $k_1$ and $k_2$.  Thus we may assume that $\tau_2 > 0$.

If $\tau_1 > 0$, we have the identity
$$ |\xi_1| + |\xi_2| - |\xi_1 + \xi_2| = (|\tau_0| - |\xi_0|) - (|\tau_1| - |\xi_1|) - (|\tau_2| - |\xi_2|).$$
The right-hand side has magnitude $O(2^{j_{max}})$, and in the imbalanced case has magnitude $\sim 2^{j_{max}}$.  The left-hand side can be rewritten as
$$ |\xi_1| + |\xi_2| - |\xi_1 + \xi_2| = 2 \frac{|\xi_1| |\xi_2| - \xi_1 \cdot \xi_2}{|\xi_1| + |\xi_2| + |\xi_1 + \xi_2|}$$
so that
$$ |\xi_1| + |\xi_2| - |\xi_1 + \xi_2| \sim \angle(\xi_1, \xi_2)^2 2^{k_1}.$$
Combining this with the previous we see that we are in either \eqref{pm-c} or \eqref{imbal}.

If $\tau_1 < 0$, we have the identity
$$ |\xi_1| - |\xi_2| + |\xi_1 + \xi_2| = -(|\tau_0| - |\xi_0|) - (|\tau_1| - |\xi_1|) + (|\tau_2| - |\xi_2|).$$
The right-hand side has magnitude $O(2^{j_{max}})$, and in the imbalanced case has magnitude $\sim 2^{j_{max}}$.  The left-hand side can be rewritten as
$$ |\xi_1| - |\xi_2| + |\xi_1 + \xi_2| = 2 \frac{|\xi_1| |\xi_1 + \xi_2| - \xi_1 \cdot (\xi_1 + \xi_2)}{|\xi_1| + |\xi_2| + |\xi_1 + \xi_2|}$$
so that
$$ |\xi_1| + |\xi_2| - |\xi_1 + \xi_2| \sim \angle(-\xi_1, \xi_1 + \xi_2)^2 2^{k_1}.$$
By the sine rule we thus have
$$ |\xi_1| + |\xi_2| - |\xi_1 + \xi_2| \sim \angle(-\xi_1, \xi_2)^2 2^{k_1}.$$
Combining this with the previous we see that we are in either \eqref{pm-c} or \eqref{imbal}.
\end{proof}

One could of course formulate variants of this Lemma when we drop the $j_i \leq k_i - C$ hypothesis, but the number of additional cases becomes excessive, and we shall just treat these cases by hand whenever they arise.  Most of these cases are rather simple, as the modulations are so large that the geometry of the cone becomes irrelevant.  However, there is one case in this category, namely the $(++)$ case mentioned above, when $\phi$, $\psi$ have opposing high frequencies and low modulation, and $\phi \psi$ has low frequency and high modulation, which warrants some special attention in Case 2(a).3 of Section \ref{null-sec}.

In the situations considered in Lemma \ref{geometry-lemma}, we refer to $\phi$, $\psi$ as \emph{inputs} and $P_{k_0} Q_{j_0} L(\phi, \psi)$ or $P_{k_0} Q_{j_0} L(\phi_{,\alpha}, \psi^{,\alpha})$ as the \emph{output}.  We refer to the inputs and output collectively\footnote{Since one can use duality to swap an input with the output, it is natural to consider both inputs and outputs on the same footing.} as \emph{functions}.  We say that one function \emph{dominates} another if its modulation of the former is much greater than that of the latter. The imbalanced modulation case thus occurs when one of the functions dominate the other two.  In such a case, our estimates will usually be proven by dividing up into caps $\kappa$ of the radius suggested by the above proposition, and then applying Lemma \ref{sklw-prop} or the third term in \eqref{sk-def} (for the inputs) and \eqref{f-null} or \eqref{sk-def} (for the output).  The balanced modulation case occurs when none of the functions dominate each other, and in this case one usually uses H\"older's inequality and the estimates \eqref{sk-energy}, \eqref{sk-infty}, \eqref{qbound}, \eqref{qbound-strichartz} or Lemma \ref{n-xsb} (for the inputs), and \eqref{f-l12}, \eqref{f-xsb}, or \eqref{sk-def} (for the output).

\section{The core product estimate}\label{prod-sec}

The purpose of this section is to state and prove the core product estimate Lemma \ref{core}, which lies at the heart of the proof of all the bilinear and trilinear estimates in Theorem \ref{spaces} (\eqref{algebra}, \eqref{sk-skp}, \eqref{sk-sk}, \eqref{ur-algebra}, \eqref{ur-algebra-sk}, \eqref{null}, \eqref{o-lemma}):

\begin{lemma}\label{core}
Let $j, k, k_1, k_2$ be integers such that $j \leq \min(k_1,k_2)+O(1)$.  Then we have
\be{core-est}
\| P_k (F \psi) \|_{N[k]} \lesssim \chi^{(4)}_{k=\max(k_1,k_2)} \chi^{(4)}_{j=\min(k_1,k_2)}
\| F \|_{\dot X^{n/2-1,-1/2,\infty}_{k_1}} \| \psi \|_{S[k_2]}
\end{equation}
for all Schwartz functions $F$ on $\R^{1+n}$ with Fourier support in $2^{k_1-5} \leq D_0 \leq 2^{k_1+5}$, $D_- \sim 2^j$ and Schwartz functions $\psi \in S[k_2]$. 
\end{lemma}

The factor $\chi^{(4)}_{k=\max(k_1,k_2)}$ reflects the fact that high-high interactions (which are the only case when $k$ is significantly less than $\max(k_1,k_2)$) are weak.  The factor $\chi^{(4)}_{j = \min(k_1,k_2)}$ indicates some decay when $F$ is close to the light cone; note that if $F$ is normalized in $\dot X^{n/2-1,-1/2,\infty}_{k_1}$, then the $L^2_t L^2_x$ norm of $F$ decays like the square root of the distance to the light cone.  This gain near the light cone shall be crucial for eliminating various logarithmic divergences arising from certain types of low-high interactions, and is a reflection of the fact that small-angle interactions are weaker than large-angle interactions when $n > 1$.  (Small-angle interactions are the only interactions in which $F$, $\psi$, and $F \psi$ can all stay near the cone).

The estimates \eqref{ur-algebra}, \eqref{ur-algebra-sk} shall largely be obtained from Lemma \ref{core} via Lemma \ref{n-xsb}.  By \eqref{sf-duality} we can convert Lemma \ref{core} to an estimate of the form $S \cdot S \subset \dot X^{n/2,1/2,1}$, at which point \eqref{algebra}, \eqref{sk-skp}, \eqref{sk-sk} can largely be obtained from Lemma \ref{xn11 contain}.  Then, \eqref{null} will be largely obtained from the previous estimates and \eqref{null-form}.  Finally, \eqref{o-lemma} will be obtained from all the previous estimates, plus some additional arguments to deal with some difficult sub-cases.

The rest of this section is devoted to the proof of Lemma \ref{core}.  By the Littlewood-Paley  product trichotomy we may divide into the high-high interaction case $k \leq k_2 + O(1)$, $k_1 = k_2 + O(1)$, the high-low interactions $k = k_1 + O(1)$, $k_2 \leq k_1 + O(1)$, and the low-high interactions $k = k_2 + O(1)$, $k_1 \leq k_2 + O(1)$.

\divider{Case 1.  (High-high interactions) $k \leq k_2 + O(1)$, $k_1 = k_2 + O(1)$, $j \leq k_2 + O(1)$.}

We rescale $k_1 = 0$, so $k_2 = O(1)$ and $k \leq O(1)$, and we reduce to showing
$$
\| P_k (F \psi) \|_{N[k]} \lesssim \chi^{(4)}_{k=0}
\chi^{(4)}_{k=j}
2^{-j/2} \| F \|_{L^2_t L^2_x} \| \psi \|_{S[k_2]}
$$
We now split the left-hand side into three contributions.

\divider{Case 1(a).  ($F$ does not dominate $\psi$) The contribution of $P_k (F Q_{\geq j-C} \psi)$.}

In this case we use \eqref{f-l12} to estimate this contribution by
$$
\lesssim 2^{(\frac{n}{2}-1)k} \| P_k (F Q_{\geq j-C} \psi) \|_{L^1_t L^2_x}$$
and then split into two sub-cases.

\divider{Case 1(a).1. ($F$ very close to light cone) $j < 100 k$.}

We discard $P_k$ and estimate the previous by
$$\lesssim 2^{(\frac{n}{2}-1)k} \| F \|_{L^2_t L^2_x} \| Q_{\geq j-C} \psi \|_{L^2_t L^\infty_x}.$$
By dyadic decomposition and \eqref{qbound-strichartz} (noting that the right-hand side of this estimate is decreasing in $j$) we can bound this by
$$
\lesssim 2^{(\frac{n}{2}-1)k} 
\| F\|_{L^2_t L^2_x} \min(2^{-j}, 1) \chi^{(4)}_{j \geq 0} 2^{-j/2} \| \psi \|_{S[k_2]}$$
which is acceptable.  

\divider{Case 1(a).2. ($F$ not too close to light cone) $j \geq 100 k$.}

We use Bernstein's inequality \eqref{bernstein-dual} to estimate this by
$$
\lesssim 2^{(\frac{n}{2}-1)k} 2^{nk/2} \| P_k(F Q_{\geq j-C} \psi) \|_{L^1_t L^1_x}.$$
We discard $P_k$ and estimate the previous by
$$
\lesssim 2^{(\frac{n}{2}-1)k} 2^{nk/2} \| F \|_{L^2_t L^2_x} \| Q_{\geq j-C} \psi \|_{L^2_t L^2_x}.$$
By \eqref{qbound} this is bounded by
$$
\lesssim  2^{(\frac{n}{2}-1)k} 2^{nk/2} 
\| F \|_{L^2_t L^2_x} 2^{-j/2} \frac{1}{1+2^j} 
\| \psi \|_{S[k_2]} $$
which is acceptable.

\divider{Case 1(b).  ($F$ does not dominate output) The contribution of $P_k Q_{\geq j-C} (F Q_{< j-C} \psi)$.}

We use \eqref{f-xsb} to estimate the contribution by
$$
\lesssim 
\| P_k Q_{\geq j-C} (F Q_{<j-C} \psi) \|_{\dot X^{n/2-1, -1/2,\infty}_k}.$$
We estimate this in turn by
$$
\lesssim 2^{(\frac{n}{2}-1)k} 2^{-j/2}
\| P_k(F Q_{<j-C} \psi) \|_{L^2_t L^2_x}.$$

We split into two cases.

\divider{Case 1(b).1. ($F$ very close to light cone) $j < 100k$.}

Discarding the $P_k$, we estimate this contribution by
$$
\lesssim 2^{(\frac{n}{2}-1)k} 2^{-j/2}
\| F \|_{L^2_t L^\infty_x} \| Q_{<j-C} \psi \|_{L^\infty_t L^2_x}.$$
By Lemma \ref{wimpy-strichartz}, \eqref{sk-energy} we can bound this by
$$
\lesssim 2^{(\frac{n}{2}-1)k} 2^{-j/2}
 \chi^{(4)}_{j \geq k}
\| F \|_{L^2_t L^2_x} \| \psi \|_{S[k_2] }$$
which is acceptable.  

\divider{Case 1(b).2.  ($F$ not too close to light cone) $100k \leq j < O(1)$.}

We use Bernstein's inequality \eqref{bernstein-dual} to estimate this by
$$
\lesssim 2^{(\frac{n}{2}-1)k} 2^{-j/2} 2^{nk/2}
\| F Q_{<j-C} \psi \|_{L^2_t L^1_x},$$
which we estimate by
$$
\lesssim  2^{(\frac{n}{2}-1)k} 2^{-j/2} 2^{nk/2}
\| F \|_{L^2_t L^2_x} \| Q_{<j-C} \psi \|_{L^\infty_t L^2_x}$$
which by \eqref{sk-energy} is bounded by
$$
\lesssim 2^{(\frac{n}{2}-1)k} 2^{-j/2} 2^{nk/2}
 \| F \|_{L^2_t L^2_x} \| \psi \|_{S[k_2]}$$
which is acceptable.

\divider{Case 1(c).  ($F$ dominates) The contribution of $P_k Q_{< j-C} (F Q_{< j-C} \psi)$.}

We may assume $j \leq k + O(1)$ since the contribution vanishes otherwise.  Since we are in the unbalanced case, we shall use decompose into caps of the size suggested by Lemma \ref{geometry-lemma}.

We introduce the parameter $l := (k-j)/2 + C/4$, and then split $Q_{<j-C} = Q^+_{<j-C} + Q^-_{<j-C}$.  By conjugation symmetry it suffices to show
$$
\| \sum_{\kappa \in K_l} P_k P_{k,\kappa} Q^+_{<j-C} (
\sum_{\kappa' \in K_l} F P_{k_2,\pm \kappa'} Q^\pm_{<j-C} \psi) \|_{N[k](\R^{1+n})}
\lesssim \chi^{(4)}_{k=0} \chi^{(4)}_{j=k} 2^{-j/2} \| F \|_{L^2_t L^2_x} \| \psi \|_{S[k_2]}.$$
for both choices of sign $\pm$.

Fix $\pm$.  By Lemma \ref{geometry-lemma} we may assume that $\dist(\kappa,\kappa') \sim 2^{-l+C}$.  By \eqref{f-null} we can bound the left-hand side by
$$
\lesssim 2^{(\frac{n}{2}-1)k}
(\sum_{\kappa \in K_l}
\|P_k P_{k,\kappa} Q^+_{<j-C} (
\sum_{\kappa' \in K_l: \dist(\kappa,\kappa') \sim 2^{-l+C}} F  P_{k_2,\pm \kappa'} Q^\pm_{<j-C} \psi) \|_{NFA[\kappa]}^2)^{1/2}
$$
Since the $\kappa'$ summation is only over $O(1)$ elements we can estimate this by
$$
\lesssim 2^{(\frac{n}{2}-1)k}
(\sum_{\kappa, \kappa' \in K_l: \dist(\kappa,\kappa') \sim 2^{-l+C}}
\|P_k P_{k,\kappa} Q^+_{<j-C} (
F P_{k_2,\pm \kappa'} Q^\pm_{<j-C} \psi) \|_{NFA[\kappa]}^2)^{1/2}.
$$
By Lemma \ref{proj-disposable} we may discard $P_k P_{k,\kappa} Q^+_{<j-C}$.  By \eqref{NFAPW-dual} we can thus estimate the previous by
$$
\lesssim 2^{(\frac{n}{2}-1)k} 2^l 2^{-(n-1)l/2}
\| F \|_{L^2_t L^2_x}
(\sum_{\kappa, \kappa' \in K_l: \dist(\kappa,\kappa') \sim 2^{-l+C}}
\| P_{k_2,\pm \kappa'} Q^\pm_{<j-C} \psi \|_{S[k_2,\kappa']}^2)^{1/2}.
$$
The $\kappa$ summation is now trivial and can be discarded.  Note that if we replaced $Q^\pm_{<j-C}$ by $Q^\pm_{<k_2-2l}$ then we could use \eqref{sk-def} to bound the previous by
$$
\lesssim 2^{(\frac{n}{2}-1)k} 2^l 2^{-(n-1)l/2}
\| F \|_{L^2_t L^2_x}
\| \psi \|_{S[k_2]}
$$
which is acceptable by the definition of $l$.  Thus it remains only to control
$$
\lesssim 2^{(\frac{n}{2}-1)k} 2^l 2^{-(n-1)l/2}
\| F \|_{L^2_t L^2_x}
(\sum_{\kappa' \in K_l}
\| P_{k_2,\pm \kappa'} (Q^\pm_{<j-C} - Q^\pm_{<k_2-2l}) \psi \|_{S[k_2,\kappa']}^2)^{1/2}.
$$
By Lemma \ref{f-lemma} we may bound this by
$$
\lesssim 2^{(\frac{n}{2}-1)k} 2^l 2^{-(n-1)l/2}
\| F \|_{L^2_t L^2_x}
\| (Q^\pm_{<j-C} - Q^\pm_{<k_2-2l}) \psi \|_{\dot X^{n/2,1/2,1}_{k_2}},
$$
which we can bound using \eqref{qbound} by
$$
\lesssim 2^{(\frac{n}{2}-1)k} 2^l 2^{-(n-1)l/2}
\| F \|_{L^2_t L^2_x}
(1 + |(j-C) - (k_2-2l)|)
\| \psi \|_{S[k_2]}.
$$
But this is acceptable by our construction of $l$.  This concludes the proof of Case 1.

\divider{Case 2.  (Low-high interactions) $k = k_2 + O(1)$, $j \leq k_1 + O(1)$, and $k_1 \leq k_2 + O(1)$.}

Let $C$ be a large constant.  We may rescale $k=0$, so that $k_2 = O(1)$, $k_1 \leq O(1)$, and $j \leq k_1 + O(1)$.  We may also assume that $k_1 < -C$, since the claim follows from Case 1 otherwise.  Our task is then to show
$$
\| P_0 (F \psi) \|_{N[0](\R^{1+n})} \lesssim 
\chi^{(4)}_{j=k_1} 2^{(n/2-1)k_1} 2^{-j/2} \|F\|_{L^2_t L^2_x} \| \psi \|_{S[k_2]}.$$
We split the left-hand side into three separate contributions.

\divider{Case 2(a). (Output is not too close to the light cone) The contribution of $P_0 Q_{\geq k_1+j-2C} (F \psi)$.}

We apply \eqref{f-xsb} and estimate this contribution by
$$\lesssim \sum_{j' \geq k_1+j-2C} 
2^{-j'/2} \| P_0 Q_{j'}(F \psi) \|_{L^2_t L^2_x}.$$
We use Plancherel to discard $P_0 Q_{j'}$ and estimate this by
$$\lesssim \sum_{j' \geq k_1+j-2C} 
2^{-j'/2} \| F \|_{L^2_t L^\infty_x} \| \psi \|_{L^\infty_t L^2_x}.$$
Applying Lemma \ref{wimpy-strichartz}, and \eqref{sk-energy} we can bound this by
$$\lesssim \sum_{j' \geq k_1+j-2C} 
2^{-j'/2} 2^{nk_1/2} \chi^{(4)}_{j' \geq k_1} \| F \|_{L^2_t L^2_x} \| \psi \|_{S[k_2]}$$
which is acceptable (with a factor of about $2^{k_1/2}$ to spare).

\divider{Case 2(b). (Input $\psi$ is not too close to the light cone) The contribution of $P_0 Q_{<k_1+j-2C} (F Q_{\geq k_1+j-C} \psi)$.}

For this contribution we apply \eqref{f-l12} and discard $P_0 Q_{<k_1+j-2C}$ by Lemma \ref{q-trunc} to control this contribution by
$$
\lesssim \| F \|_{L^2_t L^\infty_x} \| Q_{\geq k_1 + j-C} \psi \|_{L^2_t L^2_x}.$$
Applying Lemma \ref{wimpy-strichartz}, and \eqref{qbound} we can bound this by
$$
\lesssim 2^{nk_1/2} \chi^{(4)}_{j=k_1} \| F \|_{L^2_t L^2_x}
2^{-(k_1+j-C)/2} \| \psi \|_{S[k_2]}$$
which is acceptable (with a factor of about $2^{k_1/2}$ to spare).

\divider{Case 2(c). (Both the input $\psi$ and the output are very close to the light cone; $F$ has the dominant modulation) The contribution of $P_k Q_{<k_1 + j-2C} (F Q_{< k_1+j-C} \psi)$.}

By conjugation symmetry it suffices to show
$$
\| P_0 Q^+_{<k_1+j-2C} (F Q^\pm_{<k_1+j-C} \psi) \|_{N[0](\R^{1+n})} \lesssim \chi^{(4)}_{k_1=j} 2^{(\frac{n}{2}-1)k_1} 2^{-j/2} \| F \|_{L^2_t L^2_x} \| \psi \|_{S[k_2]}$$
for all signs $\pm$.  We can assume that $\pm=+$ since the expression vanishes otherwise.

Since we are again in the imbalanced case, we again use sector decomposition.
We introduce the parameter $l := (k_2 - k_1 - j + C)/2$, and rewrite the left-hand side as
$$
\| \sum_{\kappa \in K_l} \sum_{\kappa' \in K_l}  P_0 P_{0,\kappa} Q^+_{<k_1+j-2C} (F P_{k_2,\kappa'} Q^+_{<k_1+j-C} \psi) \|_{N[0](\R^{1+n})}$$
From Lemma \ref{geometry-lemma} the summand vanishes unless $\dist(\kappa,\kappa') \sim 2^{-l+C}$.  Thus for each $\kappa \in K_l$ there are only $O(1)$ values of $\kappa'$ which have a non-zero contribution.  We now use \eqref{f-null} to estimate the previous by
$$
\lesssim (\sum_{\kappa \in K_l} \| \sum_{\kappa' \in K_l: \dist(\kappa,\kappa') \sim 2^{l+C}} P_0 P_{0,\kappa} Q^+_{<k_1+j-2C} (F P_{k_2,\kappa'} Q^+_{<k_1+j-C} \psi) \|_{NFA[\kappa]}^2)^{1/2}.$$
Since the interior summation is only over $O(1)$ choices, we may estimate this by 
$$
\lesssim (\sum_{\kappa, \kappa' \in K_l: \dist(\kappa,\kappa') \sim 2^{l+C}} 
\| P_0 P_{0,\kappa} Q^+_{<k_1+j-2C} (F P_{k_2,\kappa'} Q^+_{<k_1+j-C} \psi) \|_{NFA[\kappa]}^2)^{1/2}.$$
By Lemma \ref{proj-disposable} we may discard $P_0 P_{0,\kappa} Q^+_{<k_1+j-2C}$.  By \eqref{NFAPW-dual} we can then bound the previous by
$$
\lesssim 2^l 2^{-(n-1)l/2} 
(\sum_{\kappa, \kappa' \in K_l: \dist(\kappa,\kappa') \sim 2^{l+C}}
\| F \|_{L^2_t L^2_x}^2
\| P_{k_2,\kappa'} Q^+_{<k_1 + j-C} \psi \|_{S[0,\kappa']}^2)^{1/2}.$$
The $\kappa$ summation is now trivial and can be discarded. 
By construction of $l$ we have $Q_{<k_1 + j-C} = Q_{<k_2 - 2l}$.  By \eqref{sk-def} we can therefore bound the previous by
$$
\lesssim 2^l 2^{-(n-1)l/2}  \| F \|_{L^2_t L^2_x} 
\| \psi \|_{S[k_2]} \sim 2^{-(n-1)(k_1-j)/2} 2^{(\frac{n}{2}-1)k_1} 2^{-j/2} \| \phi \|_{L^2_t L^2_x} \| \psi \|_{S[k_2]}$$
which is acceptable\footnote{Note that in this case one does not gain a factor of $2^{k_1}$ as in the other cases.  This inability to gain any additional factors will be a substantial headache when it comes to prove Theorem \ref{o-lemma}.}.  This concludes the treatment of Case 2.

\divider{Case 3.  (High-low interactions). $k = k_1 + O(1)$, $j \leq k_2 + O(1)$, and $k_2 \leq k_1 + O(1)$.}

This will be a variant of the argument for Case 2.  There is less room to spare, but on the other hand one does not have to perform as severe an angular decomposition.  Let $C$ be a large constant.  We may rescale $k=0$.  We may take $k_2 < -C$ since the claim follows from Case 1 otherwise.
 
For future reference we observe that
$\| F \|_{\dot X^{n/2-1,-1/2,\infty}_{k_1}} \sim 2^{-j/2} \|F\|_{L^2_t L^2_x}$.
We split into four contributions.

\divider{Case 3(a).  ($F$ does not dominate $\psi$) The contribution of $P_0 (F Q_{\geq j-C} \psi)$.}

For this contribution we apply \eqref{f-l12} and discard $P_0$ to control this contribution by
$$
\lesssim \| F \|_{L^2_t L^2_x} \| Q_{\geq j-C} \psi \|_{L^2_t L^\infty_x}.$$
Applying \eqref{qbound-strichartz} 
we can bound this by
$$
\lesssim \| F \|_{L^2_t L^2_x} \chi^{(4)}_{k_2=j} 2^{-j/2} \| \psi \|_{S[k_2]} 
$$
which is acceptable.

\divider{Case 3(b).  (Output very far away from light cone) The contribution of $P_0 Q_{> k_2 + C} (F Q_{< j-C} \psi)$.}

We apply \eqref{f-xsb} and estimate this contribution by
$$
\lesssim 2^{-k_2/2} \| F Q_{<j-C} \psi \|_{L^2_t L^2_x}
\lesssim 2^{-k_2/2} \| F\|_{L^2_t L^2_x} \| Q_{<j-C} \psi \|_{L^\infty_t L^\infty_x},
$$
which is acceptable by \eqref{sk-infty}.

\divider{Case 3(c).  ($F$ does not dominate output) The contribution of $P_0 Q_{j-2C \leq \cdot \leq k_2 + C} (F Q_{< j-C} \psi)$.}

We apply \eqref{f-xsb} and estimate this contribution by
\be{2b-or-not2b}
\lesssim \sum_{j-2C \leq j' \leq k_2 + C} 
2^{-j'/2} \| P_0 Q_{j'}(F Q_{<j-C} \psi) \|_{L^2_t L^2_x}.
\end{equation}
Fix $j'$.  We introduce the parameter $l := (k_2 - j')/2 + C$, and write 
$$ \| P_0 Q_{j'}(F Q_{<j-C} \psi) \|_{L^2_t L^2_x}
= \| \sum_{\kappa \in K_l} P_0 Q_{j'}( F P_{k_2,\kappa} Q_{<j-C} \psi ) \|_{L^2_t L^2_x}.$$
From the geometry of the cone we see that the summand has Fourier support in the region $\{ \Theta \in C' \kappa \}$ for some $C'$ depending on $C$.  Thus the summands are almost orthogonal, and we may estimate the previous by
$$
\lesssim (\sum_{\kappa \in K_l} \| P_0 Q_{j'}( F P_{k_2,\kappa} Q_{<j-C} \psi ) \|_{L^2_t L^2_x}^2)^{1/2}$$
We use Plancherel to discard $P_0 Q_{j'}$ and then H\"older to estimate this by
$$
\lesssim (\sum_{\kappa \in K_l} \| F \|_{L^2_t L^2_x}^2
\| P_{k_2,\kappa} Q_{<j-C} \psi ) \|_{L^\infty_t L^\infty_x}^2)^{1/2}.$$
By Bernstein's inequality \eqref{bernstein-gen} we may bound this by
$$
\lesssim \| F \|_{L^2_t L^2_x} 2^{nk_2/2} 2^{-(n-1)l/2} (\sum_{\kappa \in K_l} 
\| P_{k_2,\kappa} Q_{<j-C} \psi \|_{L^\infty_t L^2_x}^2)^{1/2}.$$
By Plancherel and orthogonality we may bound this by
$$
\lesssim \| F \|_{L^2_t L^2_x} 2^{nk_2/2} 2^{-(n-1)l/2} \| Q_{<j-C} \psi \|_{L^\infty_t L^2_x}.$$
Applying \eqref{sk-energy} and then inserting this back into \eqref{2b-or-not2b}, we can therefore estimate this contribution by
$$\lesssim \sum_{j-2C \leq j' \leq k_2 + O(1)} 
2^{-j'/2} \| F \|_{L^2_t L^2_x} 2^{-(n-1)l/2} \| \psi \|_{S[k_2]}$$
which is acceptable.

\divider{Case 3(d).  ($F$ dominates) The contribution of $P_0 Q_{< j-2C} (F Q_{<j-C} \psi)$.}

By conjugation symmetry it suffices to show
$$
\| P_0 Q^+_{<j-2C} (F Q^\pm_{<j-C} \psi) \|_{N[0](\R^{1+n})} \lesssim \chi^{(4)}_{k_2=j} 2^{-j/2} \| F \|_{L^2_t L^2_x} \| \psi \|_{S[k_2]} $$
for all $\pm$.

Fix $\pm$.  As usual in the imbalanced case we divide into sectors.  We introduce the parameter $l := -(j-k_2+C)/2$ and rewrite the left-hand side as
$$
\| \sum_{\kappa \in K_l} \sum_{\kappa' \in K_l} P_0 P_{0,\kappa} Q^+_{<j-2C} (F P_{k_2,\pm\kappa'} Q^\pm_{<j-C} \psi) \|_{N[0](\R^{1+n})}$$
By Lemma \ref{geometry-lemma} the summand vanishes unless $\dist(\kappa,\kappa') \sim 2^{-l+C}$.  Thus for each $\kappa \in K_l$ there are only $O(1)$ values of $\kappa'$ which have a non-zero contribution.  We now use \eqref{f-null} to estimate the previous by
$$
\lesssim (\sum_{\kappa \in K_l} \| \sum_{\kappa' \in K_l: \dist(\kappa,\kappa') \sim 2^{-l+C}} P_0 P_{0,\kappa} Q^+_{<j-2C} (F P_{k_2,\pm \kappa'} Q^\pm_{<j-C} \psi) \|_{NFA[\kappa]}^2)^{1/2}.$$
Since the interior summation is only over $O(1)$ choices, we may estimate this by 
$$
\lesssim (\sum_{\kappa \in K_l} 
\sum_{\kappa' \in K_l: \dist(\kappa,\kappa') \sim 2^{-l+C}} \| P_0 P_{0,\kappa} Q^+_{<j-2C} (F P_{k_2,\pm\kappa'} Q^\pm_{<j-C} \psi) \|_{NFA[\kappa]}^2)^{1/2}.$$
We would like to discard $Q^+_{<j-2C}$, but unfortunately the tools we have (Lemmata \ref{disposable}, \ref{proj-disposable}, \ref{q-trunc}) do not allow us to do so automatically. We must therefore split
$$ Q^+_{<j-2C} = 1 - (1-Q^+_{\leq j+C}) - Q^+_{j-2C \leq \cdot \leq j+C}.$$

\divider{Case 3(d).1.  ($\psi$ dominates $F$) The contribution of $1-Q^+_{\leq j+C}$.}

From Lemma \ref{geometry-lemma} we see that this term vanishes.

\divider{Case 3(d).2.  ($\psi$ balances $F$) The contribution of $Q^+_{j-2C \leq \cdot \leq j+C}$.}

By \eqref{into-sklw} we may bound this contribution by
$$
\lesssim 2^{-j/2}
(\sum_{\kappa, \kappa' \in K_l: \dist(\kappa,\kappa') \sim 2^{l+C}} 
\| P_0 P_{0,\kappa} Q^+_{j-2C \leq \cdot \leq j + C} (F P_{k_2,\pm \kappa'} Q^\pm_{<j-C} \psi) \|_{L^2_t L^2_x}^2)^{1/2}.$$
We now use Lemma \ref{proj-disposable} to discard $P_0 P_{0,\kappa} Q^+_{j-2C \leq \cdot \leq j + C}$ and estimate this by
$$
\lesssim 2^{-j/2}
(\sum_{\kappa, \kappa' \in K_l: \dist(\kappa,\kappa') \sim 2^{l+C}} 
\| F \|_{L^2_t L^2_x}^2
\| P_{k_2,\pm \kappa'} Q^\pm_{<j-C} \psi \|_{L^\infty_t L^\infty_x}^2
)^{1/2}.$$
By the improved Bernstein's inequality \eqref{bernstein-gen} we can estimate this by
$$
\lesssim 2^{-j/2} 2^{nk_2/2} 2^{-(n-1)l/2} \| F \|_{L^2_t L^2_x}
(\sum_{\kappa , \kappa' \in K_l: \dist(\kappa,\kappa') \sim 2^{l+C}} \|  P_{k_2,\pm \kappa'} Q^\pm_{<j-C} \psi \|_{L^\infty_t L^2_x}^2
)^{1/2}.$$
The $\kappa$ sum is now trivial and can be discarded.
We can factorize $P_{k_2,\pm \kappa'} Q^\pm_{<j-C}$ as $P_{k_2,\pm \kappa'} Q_{<j-C}$ times a disposable multiplier, which we then discard. By Plancherel we can thus bound the previous by
$$
\lesssim 2^{-j/2} 2^{nk_2/2} 2^{-(n-1)l/2} \| F \|_{L^2_t L^2_x}
\| Q_{<j-C} \psi \|_{L^\infty_t L^2_x},$$
wihich is acceptable by \eqref{sk-energy}.

\divider{Case 3(d).3.  ($\psi$ arbitrary) The contribution of $1$.}

We need to control
$$
(\sum_{\kappa, \kappa' \in K_l: \dist(\kappa,\kappa') \sim 2^{l+C}} 
\| P_0 P_{0,\kappa} (F P_{k_2,\pm \kappa'} Q^\pm_{<j-C} \psi) \|_{NFA[\kappa]}^2)^{1/2}.$$
We discard $P_0 P_{0,\kappa}$ and use \eqref{NFAPW-dual} to bound the previous by
$$
\lesssim 2^l 2^{-(n-1)l/2} 2^{-k_2/2} 
(\sum_{\kappa, \kappa' \in K_l: \dist(\kappa,\kappa') \sim 2^{l+C}} 
\| F \|_{L^2_t L^2_x}^2
\|  P_{k_2,\pm\kappa'} Q^\pm_{<j-C} \psi \|_{S[0,\kappa']}^2
)^{1/2}.$$
The $\kappa$ summation is trivial and can be discarded.
By construction of $l$ we have $Q_{<j-C} = Q_{<k_2 - 2l}$.  By \eqref{sk-def} we can therefore bound the previous by
$$
\lesssim 2^{-l} 2^{(n-1)l/2} 2^{-k_2/2} \| F \|_{L^2_t L^2_x} 
\| \psi \|_{S[k_2]} 
$$
which is acceptable.  This concludes the proof of Case 3, and thus of Lemma \ref{core}.
\endprf
 
In the next three sections we show how the above core product estimates
can be used to prove \eqref{ur-algebra}, \eqref{ur-algebra-sk}, \eqref{algebra}, \eqref{sk-skp}, \eqref{sk-sk}, and \eqref{null}. 

\section{Product estimates: The proof of \eqref{ur-algebra}, \eqref{ur-algebra-sk}}\label{ur-sec}

From Lemma \ref{core} it shall be straightforward to prove \eqref{ur-algebra} and \eqref{ur-algebra-sk}. 

To prove these estimates for $[-T,T] \times \R^n$ it suffices to do so on $\R^{1+n}$. 
By Lemma \ref{minkowski} we may replace $L(\phi,F)$ with $\phi F$.  We may also replace $N_k$, $N_{k'}$ with $N[k]$, $N[k']$.  By \eqref{stil-s} and dyadic decomposition we see that both \eqref{ur-algebra} and \eqref{ur-algebra-sk} will follow from 
\be{ur-tilde}
\| P_{k'} (\phi F) \|_{N[k'](\R^{1+n})} \lesssim \chi^{(4)}_{k' \geq k} \| \phi \|_{S(1)(\R^{1+n})} \| F \|_{N[k](\R^{1+n})}
\end{equation}
for all Schwartz $\phi \in S(1)(\R^{1+n})$ and Schwartz $F \in N[k](\R^{1+n})$.  The estimate \eqref{ur-tilde} shall also be helpful in proving \eqref{o-lemma}.

We may rescale $k'=0$.  We divide into three cases: $k > 10$, $k < -10$, and $-10 \leq k \leq 10$.

\divider{Case 1: $k > 10$ (High-high interactions).}

We may replace $\phi$ by $\phi_{k-5 < \cdot < k+5}$.  By Littlewood-Paley  decomposition it thus suffices to prove the estimate
\be{hh-ur}
\| P_0 (\phi_{k_1} F) \|_{N[0](\R^{1+n})} \lesssim \chi^{(4)}_{k = 0} \| \phi_{k_1} \|_{S[k_1]} \| F \|_{N[k](\R^{1+n})}
\end{equation}
for all $k-5 < k_1 < k+5$.

Fix $k_1$.  We may of course assume $F$ is an $N[k]$ atom.  As usual we split into three cases.

\divider{Case 1(a).  $F$ is an $L^1_t \dot H^{n/2-1}_x$ atom with frequency $2^k$.}

We use \eqref{f-l12} followed by Bernstein's inequality \eqref{bernstein-dual} to estimate the left-hand side of \eqref{hh-ur} by
$$ \lesssim \| \phi_{k_1} F \|_{L^1_t L^1_x}
\lesssim \| \phi_{k_1} \|_{L^\infty_t L^2_x} \| F \|_{L^1_t L^2_x}.$$
By \eqref{sk-energy} and the Case 1(a) assumption this is bounded by
$$ \lesssim 2^{-nk_1/2} \| \phi_{k_1}\|_{S[k_1]} 2^{-(n/2-1)k}$$
which is acceptable.

\divider{Case 1(b).  $F$ is an $\dot X^{n/2-1,-1/2,1}$ atom with frequency $2^k$ and modulation $2^j$.}

We may assume that $j \geq k_1 + 10$ since the claim follows from Lemma \ref{core}
and Lemma \ref{n-xsb} otherwise.

First consider the contribution of $Q_{<j-20} \phi_{k_1}$.  This contribution has Fourier support in $D_0 \sim 1, D_- \sim 2^j$, so by \eqref{f-xsb} we may bound it by
$$
\lesssim 2^{-j/2} \| Q_{<j-20} \phi F \|_{L^2_t L^2_x}.$$
By Bernstein's inequality \eqref{bernstein-dual} and H\"older we may bound this by
$$
\lesssim 2^{-j/2} \| Q_{<j-20} \phi \|_{L^\infty_t L^2_x} \| F \|_{L^2_t L^2_x},$$
which is acceptable by \eqref{sk-infty} and the Case 1(b) hypothesis.

Finally, consider the contribution of $Q_{>j-20} \phi_{k_1}$.  In this case we use \eqref{f-l12} and Bernstein's inequality \eqref{bernstein-dual} to bound this by
$$ \lesssim \| Q_{>j-20} \phi_{k_1} F \|_{L^1_t L^1_x}
\lesssim \| Q_{>j-20} \phi_{k_1} \|_{L^2_t L^2_x} \| F \|_{L^2_t L^2_x}$$
which is acceptable by \eqref{qbound} and the Case 1(b) hypothesis.

\divider{Case 1(c).  $F$ is a null frame atom with frequency $2^{k_1}$.}

In this case we split $F = \sum_{j \leq k + 5} Q_j F$.  From Lemma \ref{n-xsb} and Lemma \ref{core} we have 
$$ 
\| P_0 (\phi_{k_1} Q_j F) \|_{N[0](\R^{1+n})} \lesssim 
\chi^{(4)}_{k=0} \chi^{(4)}_{k=j} 2^{-(k-j)/100} \| \phi_{k_1} \|_{S[k_1]}$$
for all such $j$, and the claim follows by summing in $j$.

\divider{Case 2: $k < -10$ (High-low interactions).}

We may replace $\phi$ by $\phi_{5 < \cdot < 5}$.  By Littlewood-Paley  decomposition it thus suffices to prove the estimate
\be{hh-ur-2}
\| P_0 (\phi_{k_1} F) \|_{N[0](\R^{1+n})} \lesssim \| \phi_{k_1} \|_{S[k_1]} \| F \|_{N[k](\R^{1+n})}
\end{equation}
for all $-5 < k_1 < 5$.

We may of course assume $F$ is an $N[k]$ atom.  As usual we split into three cases.

\divider{Case 2(a).  $F$ is an $L^1_t \dot H^{n/2-1}_x$ atom with frequency $2^k$.}

We use \eqref{f-l12} to estimate the left-hand side of \eqref{hh-ur-2} by
$$ \lesssim \| \phi_{k_1} F \|_{L^1_t L^2_x}
\lesssim \| \phi_{k_1} \|_{L^\infty_t L^2_x} \| F \|_{L^1_t L^\infty_x}.$$
By \eqref{sk-energy}, Bernstein's inequality \eqref{bernstein} and the Case 1(a) assumption this is bounded by
$$ \lesssim \| \phi_{k_1}\|_{S[k_1]} 2^k$$
which is acceptable.

\divider{Case 2(b).  $F$ is an $\dot X^{n/2-1,-1/2,1}$ atom with frequency $2^k$ and modulation $2^j$.}

We may assume that $j \geq k+10$ since the claim follows from Lemma \ref{core} and Lemma \ref{n-xsb} otherwise.

Consider the contribution of $Q_{<j-20} \phi_{k_1}$.  This contribution has Fourier support in $D_0 \sim 1, D_- \sim 2^j$, so by \eqref{f-xsb} we may estimate it by
$$ 2^{-j/2} \| Q_{<j-20} \phi_{k_1} F \|_{L^2_t L^2_x}
\lesssim 2^{-j/2} \| Q_{<j-20} \phi_{k_1} \|_{L^\infty_t L^2_x} \| F \|_{L^2_t L^\infty_x}.$$
But this is acceptable by \eqref{sk-energy}, Bernstein's inequality \eqref{bernstein}, and the Case 2(b) hypothesis.

Finally, consider the contribution of $Q_{\geq j-20} \phi_{k_1}$.  We use \eqref{f-l12} to bound this by
$$ \| Q_{\geq j-20} \phi_{k_1} F \|_{L^1_t L^2_x}
\lesssim \| Q_{\geq j-20} \phi_{k_1} \|_{L^2_t L^2_x} \| F \|_{L^2_t L^\infty_x}.$$
But this is acceptable by \eqref{qbound}, Bernstein's inequality \eqref{bernstein}, and the Case 2(b) hypothesis.

\divider{Case 2(c).  $F$ is a $\pm$-null frame atom with frequency $2^{k_1}$ and angle $2^{-l}$.}

In this case we split $F = \sum_{j \leq k + 5} Q_j F$.  From Lemma \ref{n-xsb} and Lemma \ref{core} we have
$$ 
\| P_0 (\phi_{k_1} Q_j F) \|_{N[0](\R^{1+n})} \lesssim 
\chi^{(4)}_{k=j} \| \phi_{k_1} \|_{S[k_1]}$$
for all such $j$, and the claim follows by summing in $j$.

\divider{Case 3: $-10 \leq k \leq 10$ (Low-high interactions).}

We may then freely apply the Littlewood-Paley  cutoff $P_{<20}$ to $\phi$, so that $\phi$ has Fourier support in the region $D_0 \lesssim 1$.

We may assume that $F$ is an $N[k]$ atom.  We now split into three cases depending on which type of atom $F$ is.

\divider{Case 3(a).  $F$ is an $L^1_t \dot H^{n/2-1}_x$ atom with frequency $2^k$.}

In this case we simply use \eqref{f-l12} to compute
$$ \| P_0 (\phi F) \|_{N[0](\R^{1+n})}
\lesssim \| P_0 (\phi F) \|_{L^1_t L^2_x}
\lesssim \| \phi \|_{L^\infty_t L^\infty_x} \| F \|_{L^1_t L^2_x}$$
which is acceptable by \eqref{stil-def}.

\divider{Case 3(b).  $F$ is an $\dot X^{n/2-1,-1/2,1}$ atom with frequency $2^k$ and modulation $2^j$.}

In this case $F$ has Fourier support in $D_0 \sim 1$, $D_- \sim 2^j$ and
\be{psi-j}
\| F \|_{L^2_t L^2_x} \lesssim 2^{j/2}.
\end{equation}
From Lemma \ref{core} we thus have
$$ \| P_0 (\phi_{k_1} F) \|_{N[0](\R^{1+n})} \lesssim 
\| \phi_k \|_{S[k_1]} 2^{-j/2} \chi^{(4)}_{j=k_1}
\| F \|_{L^2_t L^2_x}$$
whenever $k_1 \geq j - C$.
Summing this and using \eqref{stil-def} we see that
$$ \| P_0 (\phi_{\geq j-C} F) \|_{N[0](\R^{1+n})} \lesssim 
\| \phi \|_{S(1)(\R^{1+n})} 2^{-j/2} 
\| F \|_{L^2_t L^2_x}.$$
It therefore suffices to consider the case when $\phi$ is close to the origin in frequency.  More precisely, we reduce to showing
$$ \| P_0 (\phi_{< j-C} F) \|_{N[0](\R^{1+n})} \lesssim 
\| \phi \|_{S(1)(\R^{1+n})} 2^{-j/2} 
\| F \|_{L^2_t L^2_x}.$$
We split this into two contributions.

\divider{Case 3(b).1.  ($F$ does not dominate $\phi$) The contribution of $P_0 (Q_{>j-2C} \phi_{k_1} F)$.}

By \eqref{f-l12} and the triangle inequality we estimate this contribution by
$$ \lesssim \sum_{k_1 \leq j-C} \| P_0 (Q_{>j-2C} \phi_{k_1} F) \|_{L^1_t L^2_x}.$$
Discarding $P_0$, we can bound this by
$$ \lesssim \sum_{k_1 \leq j-C} \| Q_{>j-2C} \phi_{k_1} \|_{L^2_t L^\infty_x}
\| F \|_{L^2_t L^2_x}.$$
By \eqref{qbound-strichartz}, \eqref{psi-j}, \eqref{stil-def} we can bound this by
$$ \lesssim \sum_{k_1 \leq j-C} 2^{-(j-k_1)} 2^{-j/2} \| \phi \|_{S(1)(\R^{1+n})}
2^{j/2}$$
which is acceptable.

\divider{Case 3(b).2.  ($F$ dominates $\phi$) The contribution of $P_0 (Q_{\leq j-2C} \phi F)$.}

This contribution has Fourier support in the region $D_0 \sim 1$, $D_- \sim 2^j$, and can therefore be estimated using \eqref{f-xsb} by 
$$ \lesssim 2^{j/2} \| P_0 (Q_{\leq j-2C} \phi_{\leq j-C} F) \|_{L^2_t L^2_x}.$$
We then discard $P_0$ and bound this by
$$ 
\lesssim 2^{j/2} \| Q_{\leq j-2C} P_{\leq j-C} \phi \|_{L^\infty_t L^\infty_x}
\| F \|_{L^2_t L^2_x}.$$
But this is acceptable by \eqref{stil-def}, \eqref{psi-j}, and Lemma \ref{disposable}.  This concludes the treatment of the case when $F$ is an $\dot X^{n/2-1,-1/2,1}$ atom.

\divider{Case 3(c).  $F$ is a $\pm$-null frame atom with frequency $2^k$ and angle $2^{-l}$.}

By conjugation symmetry we may take $\pm = +$.
We then have a decomposition $F = \sum_{\kappa \in K_l} F_\kappa$ such that each $F_\kappa$ has Fourier support in the region 
$$
\{ (\tau, \xi): \tau > 0; D_- \lesssim 2^{-2l}; D_+ \sim 1; \Theta \in \frac{1}{2} \kappa \}
$$
and such that
\be{f-nfa}
(\sum_{\kappa \in K_l} \| F_\kappa \|_{NFA[\kappa]}^2)^{1/2} \lesssim 1.
\end{equation}

We now split $P_0 (\phi F)$ into three pieces and treat each contribution separately.

\divider{Case 3(c).1.  ($F$ stays away from light cone, $\phi$ close to origin) The contribution of $P_0(\phi_{<-2l-2C} Q_{\geq -2l-C} F)$.}

The term $Q_{\geq -2l - C} F$ has Fourier support on the region $D_0 \sim 1, D_- \sim 2^{-2l}$.  From Lemma \ref{n-xsb} and \eqref{f-xsb} we thus see that $Q_{\geq -2l-C} F$ is a bounded linear combination of $\dot X^{n/2-1,-1/2,1}$ atoms of frequency $2^k$ and modulation $\sim 2^{-2l}$, in which case the Case 2 argument applies.  

\divider{Case 3(c).2.  ($F$ close to light cone, $\phi$ close to origin) The contribution of $P_0(\phi_{<-2l-2C} Q_{< -2l-C} F)$.}

In this case the idea is to use \eqref{lif} and the fact that $\phi$ is bounded and does not significantly affect frequency support.

We subdivide this contribution as 
$$ \| \sum_{\kappa' \in K_{l+10}} 
P_{0,\kappa'} P_0 (\sum_{\kappa \in K_l} \phi_{<-2l-2C} Q_{<-2l-C} F_\kappa) \|_{N[0](\R^{1+n})}.$$
By \eqref{f-null} we can bound this by
$$ \lesssim (\sum_{\kappa' \in K_{l+10}} 
\| P_{0,\kappa'} P_0 (\sum_{\kappa \in K_l} \phi_{<-2l-2C} Q_{<-2l-C} F_\kappa) \|_{NFA[\kappa']}^2)^{1/2}.$$
We may restrict the $\kappa$ summation to those $\kappa$ for which $\kappa' \subset \kappa$, since the contribution vanishes otherwise after applying $P_{0,\kappa'}$.  Thus for each $\kappa'$ there are only $O(1)$ values of $\kappa$ which contribute, so we can estimate the previous by
$$ \lesssim (\sum_{\kappa' \in K_{l+10}} 
\sum_{\kappa \in K_l: \kappa' \subset \kappa}
\| P_{0,\kappa'} P_0 (\phi_{<-2l-2C} Q_{<-2l-C} F_\kappa) \|_{NFA[\kappa']}^2)^{1/2}.$$
By \eqref{inscribed} we may replace the $NFA[\kappa']$ norm with the $NFA[\kappa]$ norm.  We then discard $P_{0,\kappa'} P_0$ and use \eqref{lif}, \eqref{stil-def} to bound the previous by
$$ \lesssim \| \phi \|_{S(1)(\R^{1+n})} (\sum_{\kappa' \in K_{l+10}} 
\sum_{\kappa \in K_l: \kappa'\ subset \kappa}
\| Q_{<-2l-C} F_\kappa \|_{NFA[\kappa]}^2)^{1/2}.$$
By \eqref{ftpf} we may insert $\tilde P_{k_2,\kappa}$ in front of $F_\kappa$.  By Lemma \ref{proj-disposable} we may discard $\tilde P_{k_2,\kappa} Q_{<-2l-C}$.  The claim now follows from \eqref{f-nfa} and the observation that for each $\kappa$ there are only $O(1)$ values of $\kappa'$ which contribute. 

\divider{Case 3(c).3.  ($\phi$ not too close to origin) The contribution of $P_0(\phi_{\geq -2l-2C} F)$.}

We use the triangle inequality to estimate this contribution by
$$
\lesssim \sum_{k_1 \geq -2l-2C} \sum_j \| P_0 (\phi_{k_1} Q_j F) \|_{N[0](\R^{1+n})}.$$
The case $k_1 > C$ can be ignored since their contribution vanishes.  Similarly for $j > -2l + C$.  For the remaining cases we apply Lemma \ref{core} to estimate the contribution by
$$
\lesssim \sum_{-2l-2C \leq k_1 \leq C} \sum_{j \leq -2l+C} 
\| \phi_{k_1} \|_{S[k_1]} 2^{-j/2} \chi^{(4)}_{j=k_1}
\| Q_j F \|_{L^2_t L^2_x}.$$
By \eqref{stil-def} and Lemma \ref{n-xsb} we can bound this by
$$
\lesssim \sum_{-2l-2C \leq k_1 \leq C} \sum_{j \leq -2l+C} 
\| \phi \|_{S(1)(\R^{1+n})} \chi^{(4)}_{j=k_1} $$
which is acceptable.
\endprf

\section{Algebra estimates: The proof of \eqref{algebra}, \eqref{sk-skp}, \eqref{sk-sk}}\label{algebra-sec}

We now use the core product estimates in Section \ref{prod-sec}, combined with some duality arguments and some additional work, to obtain the algebra estimates \eqref{algebra}, \eqref{sk-skp}, \eqref{sk-sk}.

\divider{Step 1.  Obtain decay near the light cone.}

The purpose of this step is to prove the following lemma, which asserts that if $\phi \in S[k_1]$ and $\psi \in S[k_2]$, then $\phi \psi$ begins to decay once one approaches within $\min(2^{k_1}, 2^{k_2})$ of the light cone in frequency space\footnote{We shall prove this Lemma by duality as this is by far the quickest way to do so.  However, it is an instructive exercise to prove this Lemma directly.  Geometrically, the lemma reflects the fact that the portion of $\phi \psi$ this close to the light cone arises from small angle interactions, which give a better contribution than large angle interactions in dimensions 2 and higher.}.  This lemma shall also be useful in dealing with several sub-cases of the trilinear estimate \eqref{o-lemma}.

\begin{lemma}\label{decay-near-cone-lemma}
Let $j, k_1, k_2, k$ be such that $j \leq \min(k_1,k_2) + O(1)$.  Then
$$
\| P_k Q_j (\phi \psi) \|_{\dot X^{n/2,1/2,1}_k} \lesssim 
\chi^{(4)}_{k = \max(k_1,k_2)} \chi^{(4)}_{j = \min(k_1,k_2)}
\| \phi \|_{S[k_1]} \| \psi \|_{S[k_2]}
$$
for all Schwartz functions $\phi \in S[k_1]$, $\psi \in S[k_2]$.
\end{lemma}

\begin{proof} 
By symmetry we may assume $k_1 \leq k_2$.  By scaling we may take $k=0$. By conjugation symmetry we may assume that $\phi$, $\psi$ are real-valued.

We may assume that either $k_2 = O(1)$ (low-high interaction) or that $k_2 \geq O(1)$ and $k_1 = k_2 + O(1)$ (high-high interaction).  

By scaling we may take $k = 0$. 
We may of course replace $\psi$ by $P_{k_2-5 \leq \cdot \leq k_2+5} \psi$.
The left-hand side is
$$ \sim 2^{j/2} \| P_0 Q_j (\phi \psi) \|_{L^2_t L^2_x},$$
which by duality is 
$$ \sim 2^{j/2} |\langle \psi, P_{k_2-5 \leq \cdot \leq k_2+5}(\phi P_0 Q_j F) \rangle|$$
for some $F$ in the unit ball of $L^2_t L^2_x$.  By \eqref{sf-duality} and \eqref{compat} this is
\be{dnc}
\lesssim 2^{j/2} 2^{-(n-1)k_2} \sum_{k'_2 = k_2 + O(1)}
\| \psi \|_{S[k_2]} \| P_{k'_2}(\phi P_0 Q_j F) \|_{N[k'_2]}.
\end{equation}
Suppose we are in the low-high interaction case $k_2 = O(1)$.  By Lemma \ref{core} we can estimate \eqref{dnc} by
$$\lesssim 2^{j/2} \sum_{k'_2 = k_2 + O(1)}
\| \psi \|_{S[k_2]} \chi^{(4)}_{j = k_1} \| \phi \|_{S[k_1]}
\| P_0 Q_j F \|_{\dot X^{n/2-1,-1/2,1}_0}$$
which is acceptable by the normalization of $F$.

Now suppose we are in the high-high interaction case $k_2 \geq O(1)$, $k_1 = k_2 + O(1)$.  By Lemma \ref{core} we may estimate \eqref{dnc} by
$$
\lesssim 2^{j/2} 2^{-(n-1)k_2}
\sum_{k'_2 = k_2 + O(1)}
2^{j/2} \| \psi \|_{S[k_2]} \chi^{(4)}_{j=0}  \| \phi \|_{S[k_1]}
\| P_0 Q_j F \|_{\dot X^{n/2-1,-1/2,1}_0}
$$
which is acceptable by the normalization of $F$.
\end{proof}

\divider{Step 2.  Control frequency-localized interactions.}

The purpose of this step is to prove the estimate

\be{localized-algebra}
\| P_k (\phi \psi) \|_{S[k]} \lesssim 
\chi^{(4)}_{k=\max(k_1,k_2)} 
\| \phi \|_{S[k_1]} \| \psi \|_{S[k_2]}
\end{equation}
whenever $k, k_1, k_2$ are integers and $\phi \in S[k_1]$, $\psi \in S[k_2]$ are Schwartz functions.  In this step we shall also prove the strengthened estimate
\be{strong-algebra}
\| P_k (\phi \psi) \|_{S[k]} \lesssim 
\| \phi \|_{S(1)} \| \psi \|_{S[k_2] }
\end{equation}
in the low-high interaction case $k = k_2 + O(1)$.

To prove \eqref{localized-algebra} we may use the Littlewood-Paley product trichotomy and symmetry to divide into the high-high interaction case $k \leq k_2 + O(1)$, $k_1 = k_2 + O(1)$ and the low-high interaction case $k = k_2 + O(1)$, $k_1 \leq k_2 + O(1)$.  In the latter case we see that from \eqref{stil-s} that it will suffice to prove \eqref{strong-algebra}.

\divider{Case 2(a).  (High-high interactions) Proof of \eqref{localized-algebra} when $k \leq k_2 + O(1)$, $k_1 = k_2 + O(1)$.}

By scale invariance we may take $k_2 = O(1)$, so $k_1 = O(1)$ and $k \leq O(1)$.

First consider the $L^\infty_t L^2_x$ component of the $S[k]$ norm in \eqref{sk-def}.  Discarding the $P_k$ and applying Leibnitz and H\"older we have
\be{nabla-k}
\| \nabla_{x,t} P_k(\phi \psi) \|_{L^\infty_t L^2_x}
\lesssim \| \nabla_{x,t} \phi \|_{L^\infty_t L^2_x} \| \psi \|_{L^\infty_t L^\infty_x} + \| \phi \|_{L^\infty_t L^\infty_x} \| \nabla_{x,t} \psi \|_{L^\infty_t L^2_x},
\end{equation}
and so by \eqref{sk-def}, \eqref{sk-infty} we see that the $L^\infty_t L^2_x$ component of the $S[k]$ norm is acceptable.

To deal with the other two components we split $P_k(\phi \psi) = P_k Q_{<k_1+C}(\phi \psi) + P_k Q_{\geq k_1 + C}(\phi \psi)$.  To deal with the former term we observe from Lemma \ref{decay-near-cone-lemma} and a summation in $j$ that
$$ \| P_k Q_{<k_1+C}(2^{-k} \nabla_{x,t} \phi \psi) \|_{\dot X^{n/2,1/2,1}_k}
\lesssim
\chi^{(4)}_{k=0}
\| \phi \|_{S[k_1]} \| \psi \|_{S[k_2]},$$
which implies from Lemma \ref{xn11 contain} that $P_k Q_{<k_1+C}(\phi \psi)$ is acceptable.  Thus it only remains to consider $P_k Q_{\geq k_1 + C}(\phi \psi)$.  The third term in \eqref{sk-def} now vanishes, so we reduce to showing
$$ 
2^{nk/2} 2^{j-k}
\| P_k Q_j(\phi \psi) \|_{L^2_t L^2_x}
\lesssim
\chi^{(4)}_{k=0} 2^{-j/2}
\| \phi \|_{S[k_1]} \| \psi \|_{S[k_2]}$$
for all $j \geq k_1 + C$.

Fix $j$.  Consider the contribution of $Q_{>j-C} \phi$.  By Plancherel we may discard $P_k Q_j$, and estimate this contribution by
$$ 
\lesssim 2^{nk/2} 2^{j-k}
\| Q_{>j-C} \phi \|_{L^2_t L^\infty_x} \|  \psi \|_{L^\infty_t L^2_x}$$
which is acceptable by \eqref{qbound-strichartz}, \eqref{sk-energy}.  Similarly one may dispose of the contribution of $Q_{\leq j-C} \phi$ and $Q_{> j-C} \psi$, which leaves only the contribution of $Q_{\leq j-C} \phi$ and $Q_{> j-C} \psi$.  But this vanishes by the assumption on $j$, and so we are done.

\divider{Case 2(b).  (Low-high interactions) Proof of \eqref{strong-algebra} when $k = k_2 + O(1)$, $k_1 \leq k_2 + O(1)$.}

This case will not follow as easily from Lemma \ref{decay-near-cone-lemma} as Case 2(a), because of the possible logarithmic pile-up of low frequencies.  However, we can compensate for this because of the $L^\infty_t L^\infty_x$ control in \eqref{s-def}.

By scale invariance we may take $k = 0$, so that $k_2 = O(1)$.  By a limiting argument we may assume that $\psi$ is Schwartz, and that $\phi$ is a Schwartz function plus a constant.  If $\phi$ is a constant then the claim is trivial from Lemma \ref{tech-lemma}, so we may assume that $\phi$ is actually Schwartz.

Expanding out \eqref{sk-def}, we have to show that
\be{energy-lowhigh}
\| P_0 \nabla_{x,t} (\phi \psi) \|_{L^\infty_t L^2_x} \lesssim 
\| \phi \|_{S(1)(\R^{1+n})} \| \psi \|_{S[k_2]},
\end{equation}
that
\be{qj-lowhigh}
(1 + 2^j) \| P_0 Q_j (\phi \psi) \|_{L^2_t L^2_x} \lesssim 2^{-j/2}  
\| \phi \|_{S(1)(\R^{1+n})} \| \psi \|_{S[k_2]}
\end{equation}
for all $j \in \Z$, and that
\be{f-lowhigh}
(\sum_{\kappa \in K_l} \| P_{0,\kappa} P_0 Q^+_{<-2l} (\phi \psi) \|_{S[0,\kappa]}^2)^{1/2}
\lesssim \| \phi \|_{S(1)(\R^{1+n})} \| \psi \|_{S[k_2]}
\end{equation}
for all $l > 10$ (the corresponding estimate for the $-$ sign follows by conjugation symmetry).

\divider{Step 2(b).1.  Energy estimate: Proof of \eqref{energy-lowhigh}.}

We use \eqref{nabla-k} from Step 2(a).  The claim then follows from \eqref{sk-infty}, \eqref{sk-def}, \eqref{stil-def} and dyadic decomposition. 

\divider{Step 2(b).2.  $\dot X^{n/2,1/2,\infty}$ estimate: Proof of \eqref{qj-lowhigh}.}

Fix $j$.  We first deal with the easy case $j > C$.   Consider the contribution of $Q_{>j-C} \psi$.  By Plancherel we may discard $P_0 Q_j$ and estimate the previous by
$$ \lesssim (1 + 2^j) \| \phi \|_{L^\infty_t L^\infty_x} \| Q_{>j-C} \psi \|_{L^2_t L^2_x}$$
which is acceptable by \eqref{infty-control}, \eqref{qbound}.  The contribution of $Q_{\leq j-C} \phi$ and $Q_{\leq j-C} \psi$ vanish, so we only need consider the contribution of $Q_{>j-C} \phi$ and $Q_{\leq j-C} \psi$.  By Plancherel we may discard $P_0 Q_j$ and estimate this contribution by
$$ \lesssim (1 + 2^j) \| Q_{>j-C} \phi \|_{L^2_t L^2_x} \| Q_{\leq j-C} \psi \|_{L^\infty_t L^\infty_x}.$$
But this is acceptable by \eqref{qbound}, \eqref{sk-infty}, \eqref{stil-def} and dyadic decomposition.

Thus we may assume that $j < O(1)$.  From Lemma \ref{decay-near-cone-lemma}, \eqref{stil-def} and dyadic decomposition we see that
$$
\| P_0 Q_j (\phi_{\geq j-C} \psi) \|_{L^2_t L^2_x} \lesssim 2^{-j/2}  
\| \phi \|_{S(1)(\R^{1+n})} \| \psi \|_{S[k_2]}.$$
Thus we need only show that
$$
\| P_0 Q_j (\phi_{<j-C} \psi) \|_{L^2_t L^2_x} \lesssim 2^{-j/2}  
\| \phi \|_{S(1)(\R^{1+n})} \| \psi \|_{S[k_2]}.$$
We may replace $\psi$ by $Q_{j-10 < \cdot < j+10} \psi$ since the left-hand side vanishes otherwise.  We then use Plancherel to discard $P_0 Q_j$ and estimate the left-hand side by
$$ \| \phi_{<j-C} \|_{L^\infty_t L^\infty_x} \| Q_{j-10 < \cdot < j+10} \psi \|_{L^2_t L^2_x},$$
which is acceptable by \eqref{stil-def}, \eqref{qbound}.  This proves \eqref{qj-lowhigh}.

\divider{Step 2(b).3.  Null frame estimate: Proof of \eqref{f-lowhigh}.}

Fix $l > 10$.  We divide \eqref{f-lowhigh} into two contributions.

\divider{Case 2(b).3(a).  ($\phi$ not too close to origin) The contribution of $\phi_{>-2l-C}$.}

By Lemma \ref{f-lemma} we can bound this contribution by
$$
\| P_0 Q^+_{<-2l} (\phi_{>-2l-C} \psi) \|_{\dot X^{n/2,1/2,1}_0}
\lesssim \sum_{j < -2l+C} 2^{j/2} 
\| P_0 Q^+_j (\phi_{>-2l-C} \psi) \|_{L^2_t L^2_x}.$$
But the right-hand side is acceptable by Lemma \ref{decay-near-cone-lemma} (estimating $Q^+_j$ by $Q_j$), \eqref{stil-def} and dyadic decomposition.

\divider{Case 2(b).3(b).  ($\phi$ close to origin) The contribution of $\phi_{\leq -2l-C}$.}

In this case the idea is to use \eqref{liff} and the fact that $\phi$ is bounded and does not significantly affect frequency support.

We may assume that $\psi$ has Fourier support in $2^{-2.5} \leq D_0 \leq 2^{2.5}$ since this contribution vanishes otherwise.  We first take advantage of Lemma \ref{commutator} to write
$$ P_0(\phi_{\leq -2l-C} \psi) = \phi_{\leq -2l-C} \psi_0 + L(\nabla \phi_{\leq -2l-C}, \psi).$$
We consider the contribution of each term separately.

\divider{Case 2(b).3(b).1.  (Commutator term) The contribution of $L(\nabla \phi_{\leq -2l-C}, \psi)$.}

By Lemma \ref{f-lemma} we may control this by
$$
\lesssim \| Q^+_{<-2l} L(\nabla \phi_{\leq -2l-C}, \psi) \|_{\dot X^{n/2,1/2,1}_0}$$
which we can bound using \eqref{l-lp-d} by
$$
\lesssim \sum_{j < -2l+C} \sum_{k_1 \leq -2l-C} 2^{j/2} 2^{k_1}
\| Q_j L(\phi_{k_1}, \psi) \|_{L^2_t L^2_x}.$$
However, from \eqref{qj-lowhigh}, \eqref{stil-def}, and Lemma \ref{minkowski} we have
$$
\| Q_j L(\phi_{k_1}, \psi) \|_{L^2_t L^2_x} \lesssim 2^{-j/2} \| \phi \|_{S(1)(\R^{1+n})} \| \psi \|_{S[k_2]}$$
while from Lemma \ref{decay-near-cone-lemma}, \eqref{stil-def}, and Lemma \ref{minkowski} we have
$$
\| Q_j L(\phi_{k_1}, \psi) \|_{L^2_t L^2_x} \lesssim 2^{-j/2} \chi^{(4)}_{k_1=j} \| \phi \|_{S(1)(\R^{1+n})} \| \psi \|_{S[k_2]}$$
when $j \leq k_1+O(1)$.  Combining these estimates we see that the contribution of $L(\nabla \phi_{\leq -2l-C}, \psi)$ is acceptable (in fact we have a sizeable gain in $l$).

\divider{Case 2(b).3(b).2.  (Main term) The contribution of $\phi_{\leq -2l-C} \psi_0$.}

We need to show
\be{intern}
(\sum_{\kappa \in K_l} \| P_{0,\kappa} Q^+_{<-2l} (\phi_{\leq -2l-C} \psi_0) \|_{S[0,\kappa]}^2)^{1/2}
\lesssim
\| \phi \|_{S(1)(\R^{1+n})} \| \psi \|_{S[k_2]}.
\end{equation}
We may freely replace $\psi_0$ by 
$$ \sum_{\kappa' \in K_{l+10}: \kappa' \subset \kappa } P_{0,\kappa'} Q^+_{<-2l+C} \psi_0.$$
By Lemma \ref{proj-disposable} we may discard 
$P_{0,\kappa} Q^+_{<-2l}$ and estimate this contribution by
$$
\lesssim (\sum_{\kappa \in K_l} \sum_{\kappa' \in K_{l+10}: \kappa' \subset \kappa}
\| \phi_{\leq -2l-C} 
P_{0,\kappa'} Q^+_{<-2l+C} \psi_0 \|_{S[0,\kappa]}^2)^{1/2}.
$$

We now discard a technical contribution.

\divider{Case 2(b).3(b).2(a).  ($\psi$ stays away from light cone) The contribution of $Q_{-2l-20 \leq \cdot <-2l+C} \psi_0$.}

This contribution can be estimated by \eqref{into-sklw} by
$$ \lesssim (\sum_{\kappa \in K_l} 
\sum_{\kappa' \in K_{l+10}: \kappa' \subset \kappa}
\| \phi_{\leq -2l-C} P_{0,\kappa'} Q^+_{-2l-2C \leq \cdot \leq -2l+C} \psi_0 \|_{\dot X^{n/2,1/2,1}_0}^2)^{1/2}.
$$
We may estimate the $\dot X^{n/2,1/2,1}_0$ norm by $\lesssim 2^{-l}$ times the $L^2_t L^2_x$ norm.  By \eqref{stil-def} we can then estimate the previous by
$$ \lesssim 2^{-l} \| \phi \|_{S(1)(\R^{1+n})} (\sum_{\kappa \in K_l} 
\sum_{\kappa' \in K_{l+10}: \kappa' \subset \kappa}
\| P_{0,\kappa'} Q^+_{-2l-2C \leq \cdot \leq -2l+C} \psi_0 \|_{L^2_t L^2_x}^2)^{1/2},
$$
which by Plancherel is bounded by
$$ \lesssim 2^{-l} \| \phi \|_{S(1)(\R^{1+n})}
\|  Q^+_{-2l-2C \leq \cdot \leq -2l+C} \psi_0 \|_{L^2_t L^2_x},
$$
which is acceptable by \eqref{qbound}.

\divider{Case 2(b).3(b).2(b).  ($\psi$ close to light cone) The contribution of $Q_{<-2l-20} \psi_0$.}

We need to control
\be{side}
\lesssim (\sum_{\kappa \in K_l} 
\sum_{\kappa' \in K_{l+10}: \kappa' \subset \kappa} 
\| \phi_{\leq -2l-C} P_{0,\kappa'} Q^+_{<-2l-20} \psi_0 \|_{S[0,\kappa]}^2)^{1/2}.
\end{equation}
By \eqref{liff} and \eqref{stil-def} we can estimate this
$$ \lesssim \| \phi \|_{S(1)(\R^{1+n})}
(\sum_{\kappa \in K_l} 
\sum_{\kappa' \in K_{l+10}: \kappa' \subset \kappa} 
\| P_{0,\kappa'} Q^+_{<-2l-20} \psi_0 \|_{S[0,\kappa]}^2)^{1/2}.
$$
By \eqref{inscribed-sklw} we may replace the $S[0,\kappa]$ norm by the $S[0,\kappa']$ norm. By \eqref{sk-def} the previous is therefore bounded by
$$ \lesssim \| \phi \|_{S(1)(\R^{1+n})} \| \psi_0 \|_{S[0]},$$
which is acceptable by Lemma \ref{tech-lemma}.  This concludes the proof of \eqref{localized-algebra}.

\divider{Step 3.  Prove \eqref{algebra}.}

Since the $L^\infty_t L^\infty_x$ control on $\phi \psi$ is immediate from \eqref{infty-control}, it thus suffices by \eqref{s-def} to show that
$$
\| P_k(\phi \psi) \|_{S[k]} \lesssim c_k \| \phi \|_{S(c)(\R^{1+n})} \| \psi \|_{S(c)(\R^{1+n})}
$$
for all $k$.  By scale invariance we may take $k=0$.  By the Littlewood-Paley trichotomy we can decompose
\begin{eqnarray*}
P_0 (\phi \psi) &= P_0 ( (\phi_{\leq -5}) \psi ) + P_0( (\phi_{-5 < \cdot < 5}) \psi) + P_0 ( \phi_{\geq 5}) \psi)\\
&= P_0( (\phi_{\leq -5}) \psi_{-5 < \cdot < 5}) + P_0( (\phi_{-5 < \cdot < 5}) \psi_{<10}) + \sum_{k_1 \geq 5} \sum_{k_1-5 \leq k_2 \leq k_1+5} P_0 ( \phi_{k_1} \psi_{k_2} ).
\end{eqnarray*}
Applying the triangle inequality and \eqref{localized-algebra}, \eqref{strong-algebra} we obtain
\begin{eqnarray*}
 \| P_0 (\phi \psi) \|_{S[0]}
&\lesssim \sum_{k_2 = O(1)} \| \phi \|_{S(1)} \| \psi_{k_2} \|_{S[k_2]} \\
&+ 
\sum_{k_1 = O(1)} \| \phi_{k_1}\|_{S[k_1]} \| \psi \|_{S(1)}\\
&+
\sum_{k_1 \geq O(1)} \sum_{k_2 = k_1 + O(1)} \chi^{(4)}_{k_1 = 0}
\| \phi_{k_1} \|_{S[k_1]} \| \psi_{k_2} \|_{S[k_2]}.
\end{eqnarray*}
But this is acceptable by \eqref{s-def}, \eqref{stil-s}.  This completes the proof of \eqref{algebra}.

\divider{Step 4.  Prove \eqref{sk-skp}.}

By Lemma \ref{minkowski} we may replace $L(\phi,\psi)$ with $\phi \psi$.  By \eqref{S_k-def} and dyadic decomposition it suffices to show
$$
2^{\delta_1 |k'-k|} \sum_{k_1} \| P_{k'} (\phi_{k_1} \psi) \|_{S[k']} \lesssim 
\| \phi \|_{S_k} \| \psi \|_{S(c)}$$
for all $k'$.  Recall from hypothesis that $\psi$ has Fourier support in the region $D_0 \lesssim 2^k$.

Fix $k'$.  We divide the $k_1$ summation into three pieces.

\divider{Case 4(a).   (Low-high interactions) The contribution when $k_1 \leq k'-5$.}

In this case we may replace $\psi$ by $\psi_{k'-5 < \cdot < k'+5}$.  We may thus assume $k' \leq k + O(1)$ since this contribution vanishes otherwise.  By \eqref{localized-algebra} we may therefore estimate this contribution by
$$
2^{\delta_1 (k-k')} \sum_{k_1 \leq k'-5} \sum_{k_2 = k' + O(1)}
\| \phi_{k_1}\|_{S[k_1]} \|\psi_{k_2} \|_{S[k_2]}$$
which is acceptable by \eqref{S_k-def}, \eqref{s-def}.

\divider{Case 4(b).   (High-high interactions) The contribution when $k_1 \geq k'+5$.}

In this case we may replace $\psi$ by $\psi_{k_1-5 < \cdot < k_1+5}$. We may thus assume that $k' \geq k+O(1)$ since this contribution vanishes otherwise. By \eqref{localized-algebra} we may therefore estimate this contribution by
$$
2^{\delta_1 (k'-k)} \sum_{k_1 \geq k'-5} \sum_{k_2 = k_1 + O(1)} \chi^{(4)}_{k'=k_1} \| \phi_{k_1}\|_{S[k_1]} \|\psi_{k_2} \|_{S[k_2]}$$
which is acceptable by \eqref{S_k-def}, \eqref{s-def} (note that the $\chi^{(4)}$ gain dominates any $\delta_1$ losses).

\divider{Case 4(c).   (High-low interactions) The contribution when $k'-5 < k_1 < k'+5$.}

In this case we apply \eqref{strong-algebra} to dominate this by
$$
2^{\delta |k-k'|} \sum_{k_1 = k' + O(1)} \| \phi_{k_1}\|_{S[k_1]} \|\psi \|_{S(1)(\R^{1+n})}$$
which is acceptable by \eqref{S_k-def}, \eqref{stil-s}.  This completes the proof of \eqref{sk-skp}.

\divider{Step 5.  Prove \eqref{sk-sk}.}

As in Step 4 we may reduce to showing
$$
2^{\delta_1 |k'-k|} \sum_{k_1} \| P_{k'} (\phi_{k_1} \psi) \|_{S[k']} \lesssim 
\| \phi \|_{S_k} \| \psi \|_{S_k}$$
for all $k'$.  Again, we divide into three cases.

\divider{Case 5(a).   (Low-high interactions) The contribution when $k_1 \leq k'-5$.}

In this case we may replace $\psi$ by $\psi_{k'-5 < \cdot < k'+5}$. By \eqref{localized-algebra} we may therefore estimate this contribution by
$$
2^{\delta_1 |k-k'|} \sum_{k_1 \leq k'-5} \sum_{k_2 = k' + O(1)}
\| \phi_{k_1}\|_{S[k_1]} \|\psi_{k_2} \|_{S[k_2]}$$
which is acceptable by \eqref{S_k-def}.

\divider{Case 5(b).   (High-high interactions) The contribution when $k_1 \geq k'+5$.}

In this case we may replace $\psi$ by $\psi_{k_1-5 < \cdot < k_1+5}$. By \eqref{localized-algebra} we may therefore estimate this contribution by
$$
2^{\delta_1 |k-k'|} \sum_{k_1 \geq k'-5} \sum_{k_2 = k_1 + O(1)} \chi^{(4)}_{k'=k_1} \| \phi_{k_1}\|_{S[k_1]} \|\psi_{k_2} \|_{S[k_2]}$$
which is acceptable by \eqref{S_k-def} (note that the $\chi^{(4)}$ gain dominates any $\delta_1$ losses).

\divider{Case 5(c).   (High-low interactions) The contribution when $k'-5 < k_1 < k'+5$.}

In this case we may replace $\psi$ by $\psi_{<k_1+10}$.  We apply \eqref{localized-algebra} to dominate this by
$$
2^{\delta_1 |k-k'|} \sum_{k_1 = k' + O(1)} \sum_{k_2 \leq k_1 + O(1)}
\| \phi_{k_1}\|_{S[k_1]} \|\psi_{k_2} \|_{S[k_2]}$$
which is acceptable by \eqref{S_k-def}.  This completes the proof of \eqref{sk-sk}.

\section{Null forms: The proof of \eqref{null}}\label{null-sec}

We now use the estimates of the previous three sections, combined with \eqref{null-form}, to prove \eqref{null}.

The identity \eqref{null-form} is especially good for small angle interactions in which $\phi$, $\psi$, and $\phi \psi$ are all close to the light cone in frequency space\footnote{Indeed, it allows one - heuristically at least - to consider large angle interactions as the dominant term, although in practice the small angle interactions do introduce several technical nuisances, notably the angular decomposition into sectors which appears in the definition of $S[k]$ and $N_k$.}.
However if $\phi$ is much smaller frequency than $\psi$, and $\psi$ is far away from the light cone, then \eqref{null-form} is no longer efficient, however one can obtain a very satisfactory estimate in this case just by ignoring the null structure and writing $\phi_{,\alpha} \psi^{,\alpha}$ just as $L(\nabla_{x,t} \phi, \nabla_{x,t} \psi)$.
This improvement shall be important in the proof of Lemma \ref{o-lemma} in the next section.

We shall prove \eqref{null-form} in several stages.

\divider{Step 1.  Obtain decay away from light cone}

The purpose of this step is to obtain the bound
\be{j-improv}
\| P_k (\phi_{,\alpha} Q_j \psi^{,alpha}) \|_{N[k]} \lesssim 
\chi^{(4)}_{k=k_2} \chi^{(4)}_{j \leq k_1} \| \phi \|_{S[k_1]} \| \psi \|_{S[k_2]}
\end{equation}
whenever $k_1 \leq k_2 + O(1)$ and $j > k_1 + O(1)$.  In other words, we obtain a decay in the low-high interaction case if $\psi$ is more than $2^{k_1}$ away from the light cone.

This bound \eqref{j-improv} is quite easy to prove and does not exploit the null structure.  We may of course assume that $k \leq k_2 + O(1)$ since the contribution vanishes otherwise.  We first dispose of the contribution of $Q_{>k_1-C} \phi_{,\alpha}$.  In this case we use \eqref{f-l12} and H\"older to estimate this contribution by
$$ 2^{(n/2-1)k} \| Q_{>k_1-C} \nabla_{x,t} \phi \|_{L^2_t L^\infty_x} \|Q_j \nabla_{x,t} \psi \|_{L^2_t L^2_x},$$
which by \eqref{qbound-strichartz}, \eqref{qbound} is bounded by
$$ 2^{k_1/2} 2^{-j/2} 2^{(n/2-1) k_2} \| \phi \|_{S[k_1]} \| \psi \|_{S[k_2]}$$
which is acceptable.

It remains to control the contribution of $Q_{<k_1-C} \phi_{,\alpha}$.  This function has Fourier support in the region $D_+ \sim 2^{k_1}$,while $Q_j \psi_{,\alpha}$ has Fourier support in the region $D_+ \sim 2^{k_2} + 2^j$.  Thus by standard multiplier arguments (cf. Lemma \ref{minkowski}) it suffices to show that
$$ 2^{k_1} (2^{k_2} + 2^j) \| P_k (Q_{<k_1-C} \phi Q_j \psi) \|_{N[k]}
\lesssim 
\chi^{(4)}_{k=k_2} \chi^{(4)}_{j \leq k_1} \| \phi \|_{S[k_1]} \| \psi \|_{S[k_2]}
$$
for all $\phi \in S[k_1]$, $\psi \in S[k_2]$. 

First suppose $j = k_1 + O(1)$.  In this case we use \eqref{qbound} and \eqref{f-xsb} to obtain
$$ \| Q_j \psi \|_{N[k_2]} \lesssim 2^{-k_2} 2^{-j/2} \| \psi \|_{S[k_2]},$$
while from Lemma \ref{disposable} we have
$$ \| Q_{<k_1-C} \phi \|_{S[k_1]} \lesssim \| \phi \|_{S[k_1]}.$$
The claim then follows from Lemma \ref{core}.  

Now suppose that $j > k_1 + C$.  Then the expression inside the norm has Fourier support in $D_0 \sim 2^k$, $D_- \sim 2^j$, so we may use \eqref{f-xsb} to obtain
$$ 2^{k_1} (2^{k_2} + 2^j) \| P_k (Q_{<k_1-C} \phi Q_j \psi) \|_{N[k]}
\lesssim 2^{k_1} (2^{k_2} + 2^j) 2^{(n/2-1)k} 2^{-j/2} \| Q_{<k_1-C} \phi Q_j \psi \|_{L^2_t L^2_x}.$$
If $k = k_2 + O(1)$ the claim then follows by \eqref{sk-infty} for $\phi$ and \eqref{qbound} for $\psi$.  If $k < k_2 - C$, then $k_1 = k_2 + O(1)$.  By Bernstein's inequality  \eqref{bernstein-dual} and H\"older we may estimate the previous by
$$ \lesssim 2^{k_1} (2^{k_2} + 2^j) 2^{(n/2-1)k} 2^{-j/2} 2^{nk/2} \| Q_{<k_1-C} \phi \|_{L^\infty_t L^2_x} \| Q_j \psi \|_{L^2_t L^2_x}$$
which is acceptable by \eqref{sk-energy}, \eqref{qbound}.

\divider{Step 2.  Control null forms from $S[k_1] \times S[k_2] \to N[k]$}

The purpose of this step is to prove

\begin{lemma}\label{core-null}  Let $k_1$, $k_2$, $k$ be integers.  Then we have
\be{main-lhn}
\| P_k (\phi_{,\alpha} \psi^{,\alpha}) \|_{N[k]} \lesssim 
\chi^{(4)}_{k = \max(k_1,k_2)} 
\| \phi \|_{S[k_1]} \| \psi \|_{S[k_2]}
\end{equation}
for all $\phi \in S[k_1]$ and $\psi \in S[k_2]$. 
\end{lemma}

This lemma shall also be useful in the proof of Lemma \ref{o-lemma}.  Our main tools shall be \eqref{j-improv}, Lemma \ref{core}, \eqref{null-form}, and the estimate (from Plancherel and the identity $|\tau^2 - |\xi|^2| = D_+ D_-$)
\be{dminus}
\| \Box \psi \|_{\dot X^{n/2-1,-1/2,\infty}_k} \sim \| \nabla_{x,t} \psi \|_{\dot X^{n/2-1,1/2,\infty}_k}
\end{equation} 
whenever $\psi$ has Fourier support in $2^{k-5} \leq D_0 \leq 2^{k+5}$.  

By the Littlewood-Paley  product trichotomy and symmetry we may split into the high-high interaction case $k \leq k_2 + O(1)$, $k_1 = k_2 + O(1)$ and the low-high interaction case $k_1 \leq k_2 + O(1)$, $k = k_2 + O(1)$.

\divider{Case 2(a).  (High-high interactions) $k_1 = k_2 + O(1)$, $k \leq k_2 + O(1)$.}

By scale invariance we may take $k_1 = 0$, hence $k_2 = O(1)$ and $k \leq O(1)$. 

We begin by discarding those contributions where $\phi$ or $\psi$ is far away from the light cone.
From \eqref{j-improv} and the triangle inequality, the contribution of $Q_{\geq k_2} \psi$ is acceptable, thus we need only consider the contribution of $Q_{<k_2} \psi$.  Note that $\| Q_{<k_2} \psi \|_{S[k_2]} \lesssim \| \psi \|_{S[k_2]}$ by Lemma \ref{disposable}.  By similar reasoning we may now dispose of the contribution of $Q_{\geq 0} \phi$.  We are thus left with showing
$$\| P_k (Q_{<0} \phi_{,\alpha} Q_{<k_2} \psi^{,\alpha}) \|_{N[k]} \lesssim \chi^{(4)}_{k = 0} 
\| \phi \|_{S[0]} \| \psi \|_{S[k_2]}.$$
From \eqref{j-improv} and dyadic decomposition we obtain the bounds
$$\| P_k (Q_{<0} \phi \Box Q_{<k_2} \psi) \|_{N[k]} \lesssim \chi^{(4)}_{k = 0} 
\| \phi \|_{S[0]} \| \psi \|_{S[k_2]}.$$
$$\| P_k (\Box Q_{<0} \phi Q_{<k_2} \psi) \|_{N[k]} \lesssim \chi^{(4)}_{k = 0} 
\| \phi \|_{S[0]} \| \psi \|_{S[k_2]}.$$
By \eqref{null-form} it thus remains to show
$$\| P_k \Box (Q_{<0} \phi Q_{<k_2} \psi) \|_{N[k]} \lesssim \chi^{(4)}_{k = 0} \| \phi \|_{S[0]} \| \psi \|_{S[k_2]}.$$
Split $P_k = P_k Q_{<k+2C} + P_k Q_{\geq k+2C}$.  To deal with the contribution of $P_k Q_{<k+2C}$ we use \eqref{f-xsb} and \eqref{dminus} to estimate this by
$$ \| P_k Q_{<k+2C} (Q_{<0} \phi Q_{<k_2} \psi) \|_{\dot X^{n/2,1/2,1}_k}.$$
By Lemma \ref{decay-near-cone-lemma} we may bound this by
$$ \chi^{(4)}_{k=0} \| Q_{<0} \phi \|_{S[0]} \| Q_{<k_2} \psi \|_{S[k_2]}$$
which is acceptable by Lemma \ref{disposable}.

It remains to control interactions where the output stays away from the light cone, or more precisely that
$$\| P_k Q_{\geq k+2C} \Box (Q_{<0} \phi Q_{<k_2} \psi) \|_{N[k]} \lesssim \chi^{(4)}_{k = 0} \| \phi \|_{S[0]} \| \psi \|_{S[k_2]}.$$
By \eqref{f-xsb} we can estimate the left hand side by
$$\sum_{k+2C \leq j \leq C} 2^{(n/2-1)k} 2^{3j/2} \| P_k Q_j (Q_{<0} \phi Q_{<k_2} \psi) \|_{L^2_t L^2_x}.$$
since the $j>C$ contributions vanish.  We split this into three contributions.

\divider{Case 2(a).1.  ($\phi$ not too close to light cone) The contribution of $Q_{k+j-C < \cdot < 0} \phi Q_{<k_2} \psi$.}

In this case we discard $Q_j$, then use Bernstein's inequality \eqref{bernstein-dual} and H\"older to estimate the previous by
$$\sum_{k+2C \leq j \leq C} 2^{nk/2} 2^{(n/2-1)k} 2^{3j/2} \| Q_{k+j-C < \cdot <0} \phi \|_{L^2_t L^2_x} \| Q_{<k_2} \psi \|_{L^\infty_t L^2_x},$$
which by \eqref{qbound} and \eqref{sk-infty} is bounded by
$$\sum_{k+2C \leq j \leq C} 2^{nk/2} 2^{(n/2-1)k} 2^{3j/2} 2^{-nk/2} 2^{-(k+j-C)/2} 2^{k-j} \| \phi \|_{S[0]} \| \psi \|_{S[k_2]},$$
which is acceptable.

\divider{Case 2(a).2.  ($\phi$ is very close to light cone, but $\psi$ is not) The contribution of $Q_{\leq k+j-C} \phi Q_{k_2+k+j-C < \cdot < k_2} \psi$.}

In this case we discard $Q_j$, then use Bernstein's inequality \eqref{bernstein-dual} and H\"older to estimate the previous by
$$\sum_{k+2C \leq j \leq C} 2^{nk/2} 2^{(n/2-1)k} 2^{3j/2} \| Q_{\leq k+j-C} \phi \|_{L^\infty_t L^2_x} \| Q_{k_2+j-C < \cdot < k_2} \psi \|_{L^2_t L^2_x},$$
which by \eqref{qbound} and \eqref{sk-infty} is bounded by
$$\sum_{k+2C \leq j \leq C} 2^{nk/2} 2^{(n/2-1)k} 2^{3j/2}  \| \phi \|_{S[0]} 
2^{-nk/2} 2^{-(k_2+j-C)/2} 2^{k-j} \| \psi \|_{S[k_2]},$$
which is acceptable.

\divider{Case 2(a).3.  ($(++)$ case: both $\phi$, $\psi$ are very close to light cone) The contribution of $Q_{\leq k+j-C} \phi Q_{\leq k_2+k+j-C} \psi$.}

This case is only non-zero when $j = O(1)$.  We may assume that $k < -C$ since the sum is vacuous otherwise.  In this case the output dominates, and so we shall use sector decomposition.

Fix $j$, and $l := (C-j-k)/2$.  From the triangle inequality we have
$$\| P_k Q_j (Q_{\leq k+j-2} \phi Q_{<k_2+k+j-C} \psi) \|_{L^2_t L^2_x}
\lesssim 
\sum_{\pm,\pm'} \sum_{\kappa,\kappa' \in K_l}
\| P_k Q_j (P_{0,\pm\kappa} Q^\pm_{\leq -2l} \phi P_{k_2,\pm' \kappa} Q^{\pm'}_{<k_2-2l} \psi) \|_{L^2_t L^2_x}.
$$
From the geometry of the cone we see that the expression inside the norm vanishes unless $\pm = \pm'$ and $\dist(\kappa,-\kappa') \lesssim 2^{-l}$.  By \eqref{NFAPW-dual} we can thus estimate the previous by
$$ 2^{-(n-1)l/2}
\sum_\pm \sum_{\kappa,\kappa' \in K_l: \dist(\kappa,-\kappa') \lesssim 2^{-l}} 
\| P_{0,\pm \kappa} Q^\pm_{\leq -2l} \phi \|_{S[0,\kappa]}
\| P_{k_2,\pm \kappa} Q^\pm_{<k_2-2l} \psi) \|_{S[k_2,\kappa]}.$$
By Cauchy-Schwartz and \eqref{sk-def} we may bound this by
$$ 2^{-(n-1)l/2} \| \phi \|_{S[0]} \| \psi \|_{S[k_2]}$$
which is acceptable.

\divider{Case 2(b).  (Low-high interactions) $k_1 \leq k_2 + O(1)$, $k = k_2 + O(1)$.}

By scale invariance we may take $k = 0$, thus $k_2 = O(1)$ and $k_1 \leq O(1)$.

From \eqref{j-improv} and dyadic decomposition it suffices to show
$$\| P_0 (\phi_{,\alpha} Q_{\leq k_1+O(1)} \psi^{,\alpha}) \|_{N[0]} \lesssim \| \phi \|_{S[k_1]} \| \psi \|_{S[k_2]}.$$
First consider the contribution of $Q_j \phi_{,\alpha}$ for some $j > k_1 + C$, where is a large constant.  This term has Fourier support in $\{ D_0 \sim 1, D_- \sim 2^j\}$, so by \eqref{f-xsb} and H\"older we may estimate this contribution by
$$ 2^{-j/2} \| Q_j \nabla_{x,t} \phi \|_{L^2_t L^\infty_x} \| Q_{\leq k_1 + O(1)} \nabla_{x,t} \psi \|_{L^\infty_t L^2_x}.$$
By \eqref{qbound-strichartz}, \eqref{sk-energy} we may bound this by
by
$$ 2^{-j/2} 2^{k_1-j/2} \| \phi \|_{S[k_1]}
\| \psi \|_{S[k_2]}.$$
Summing over $j > k_1 + C$ we see that this contribution is acceptable, thus we reduce to
$$\| P_0 (Q_{\leq k_1+C} \phi_{,\alpha} Q_{\leq k_1+O(1)} 
\psi^{,\alpha}) \|_{N[0]} \lesssim \| \phi \|_{S[k_1]} \| \psi \|_{S[k_2]}.$$

By \eqref{null-form} it suffices to show that
\be{box-1}
\| P_0 (Q_{\leq k_1+C} \phi \Box Q_{\leq k_1+O(1)} \psi) \|_{N[0]} \lesssim \| \phi \|_{S[k_1]} \| \psi \|_{S[k_2]}.
\end{equation}
\be{box-2}
\| P_0 (\Box Q_{\leq k_1+C}\phi Q_{\leq k_1+O(1)} \psi) \|_{N[0]} \lesssim \| \phi \|_{S[k_1]} \| \psi \|_{S[k_2]}.
\end{equation}
\be{box-3}
\| P_0 \Box (Q_{\leq k_1+C} \phi Q_{\leq k_1+O(1)} \psi) \|_{N[0]} \lesssim \| \phi \|_{S[k_1]} \| \psi \|_{S[k_2]}.
\end{equation}

\divider{Step 2(b).1.  Proof of \eqref{box-1}.}

By Lemma \ref{core} we may bound the left-hand side by
$$ \lesssim \sum_{j \leq k_1 + O(1)} \chi^{(4)}_{j = k_1}
\| Q_{\leq k_1+C} \phi \|_{S[k_1]} \| \Box Q_j \psi \|_{\dot X^{n/2-1,-1/2,\infty}_{k_2}}.$$
By Lemma \ref{disposable} we may discard $Q_{\leq k_1+C}$.  By \eqref{dminus}, \eqref{s-def} this is bounded by
$$ \lesssim \sum_{j \leq k_1 + O(1)} \chi^{(4)}_{j = k_1}
\| \phi \|_{S[k_1]} \| \psi \|_{S[k_2]}$$
which is acceptable.

\divider{Step 2(b).2.  Proof of \eqref{box-2}.}

By repeating the argument in Step 2(a).1 (but with the roles of $\phi$, $\psi$ reversed) we can see that
$$ \| P_0 (\Box Q_{\leq k_1+C}\phi Q_{\leq k_2 + O(1)} \psi) \|_{N[0]} \lesssim \| \phi \|_{S[k_1]} \| \psi \|_{S[k_2]}.$$
Thus it suffices to show
$$ \sum_{k_1 + O(1) < j \leq k_2 + O(1)}
\| P_0 (\Box Q_{\leq k_1+C}\phi Q_j \psi) \|_{N[0]} \lesssim \| \phi \|_{S[k_1]} \| \psi \|_{S[k_2]}.$$
First consider the case when $j < k_1 + 2C$.  In this case we use Lemma \ref{core}, \eqref{dminus}, \eqref{qbound} to estimate
$$
\| P_0 (\Box Q_{\leq k_1+C}\phi Q_j \psi) \|_{N[0]} \lesssim \| \phi \|_{S[k_1]}
\| \psi \|_{S[k_2]}$$
as desired.  Now suppose $j > k_1 + 2C$.  Then the expression inside the norm has Fourier support in $D_0 \sim 0, D_- \sim 2^j$, so by \eqref{f-xsb} and H\"older we may bound this by
$$ 2^{-j/2} \| \Box Q_{\leq k_1+C} \phi \|_{L^2_t L^\infty_x} \| Q_j \psi \|_{L^\infty_t L^2_x}.$$
By Bernstein's inequality \eqref{bernstein} and \eqref{sk-energy} this is bounded by
$$ 2^{-j/2} 2^{nk_1/2} \| \Box Q_{\leq k_1+C} \phi \|_{L^2_t L^2_x}
\| \psi \|_{S[k_2]}.$$
Applying \eqref{dminus}, \eqref{qbound} we may estimate this by
$$ 2^{-j/2} 2^{k_1/2} \| \phi \|_{S[k_1]}
\| \psi \|_{S[k_2]},$$
which is acceptable after summing in $j$.

\divider{Step 2(b).3.  Proof of \eqref{box-3}.}

We may freely insert $Q_{\leq k_1 + 2C}$ in front of $P_0$.  We first consider the contribution of
$$
\| P_0 Q_{\leq k_1 + 2C} \Box (Q_{\leq k_1+C} \phi \psi) \|_{N[0]}.$$
By \eqref{f-xsb}, \eqref{dminus} we may estimate this by
$$
\lesssim \| Q_{\leq k_1 + 2C} (Q_{\leq k_1+C} \phi \psi) \|_{\dot X^{n/2,1/2,1}_1}.$$
By dyadic decomposition and Lemma \ref{decay-near-cone-lemma} we may bound this by
$$
\lesssim \sum_{j \leq k_1 + 2C} \chi^{(4)}_{j=k_1}
\| Q_{\leq k_1+C} \phi \|_{S[k_1]}  \| \psi \|_{S[k_2]}.$$
By Lemma \ref{disposable} we may discard $Q_{\leq k_1+C}$, and so we see this contribution is acceptable.

It remains to control
$$
\| P_0 Q_{\leq k_1 + 2C} \Box (Q_{\leq k_1+C} \phi Q_{\geq k_1 + O(1)} \psi) \|_{N[0]}.$$
We may freely replace $Q_{\geq k_1 + O(1)}$ with $Q_{k_1 + O(1) \leq \cdot \leq k_1 + 3C}$.  We then use \eqref{f-xsb} to estimate this by
$$
\lesssim 2^{k_1/2} \| Q_{\leq k_1+C} \phi Q_{k_1 + O(1) \leq \cdot \leq k_1 + 3C} \psi \|_{L^2_t L^2_x}.$$
But this is acceptable by \eqref{sk-infty} for $\phi$ and \eqref{qbound} for $\psi$.

\divider{Step 3.  Proof of \eqref{null}.}

By Lemma \ref{minkowski} we may replace $L(\phi_{,\alpha}, \psi^{,\alpha})$ with $\phi_{,\alpha} \psi^{,\alpha}$.
It suffices to prove this estimate on $\R^{1+n}$, since the $[-T,T] \times \R^n$ estimate then follows by truncation.  Since $N_k$ contains $N[k]$, it suffices to show that
$$
\| P_k (\phi_{,\alpha} \psi^{,\alpha}) \|_{N[k]} \lesssim 
\chi^{(1)}_{k = \max(k_1,k_2)} 
\| \phi \|_{S_{k_1}} \| \psi \|_{S_{k_2}}.$$
By the triangle inequality and Littlewood-Paley  decomposition the left-hand side is bounded by
$$\lesssim \sum_{k'_1, k'_2} \| P_k (\phi_{k'_1,\alpha} \psi_{k'_2}^{,\alpha}) \|_{N[k]}.$$
By \eqref{core-null} and \eqref{S_k-def} we may bound this by
$$
\sum_{k'_1, k'_2}
\chi^{(4)}_{k = \max(k'_1,k'_2)} \chi^{(1)}_{k_1=k'_1} \chi^{(1)}_{k_2=k'_2}
\| \phi \|_{S{k_1}} \| \psi \|_{S_{k_2}}$$
which is acceptable.

\section{Trilinear estimates: The proof of \eqref{o-lemma}}\label{o-lemma-sec}

We now prove the final major estimate in Theorem \ref{spaces}, namely \eqref{o-lemma}.  The major difficulty here is in obtaining the crucial $\chi^{(1)}_{k_1 \leq \min(k_2,k_3)}$ gain; without this gain, the lemma follows fairly easily from \eqref{ur-tilde} and \eqref{null}.  In some cases one can use the refined estimates \eqref{j-improv} and Lemma \ref{decay-near-cone-lemma} to eke out this gain, but there is one difficult case (Case 4, when $k_1$ is much larger than $k_2$, $k_3$) in which this does not suffice, in which case we must use more difficult arguments.

We apologize in advance for the convoluted nature of the arguments here, which are an amalgam of five or six separate attacks on this estimate by the author.  Each attack was only able to deal with a portion of the cases, but together they eventually manage to cover all the possible trilinear interactions between $\phi^{(1)}$, $\phi^{(2)}$, and $\phi^{(3)}$.  The crude piecing together of these distinct attacks is responsible for the small constant $\delta_1$ in the exponential gain, which is of course not optimal.

We turn to the details.  As usual we shall work on $\R^{1+n}$ rather than $[-T,T] \times \R^n$. We may replace $N_k$ by $N[k]$.  By Lemma \ref{minkowski} we may replace $L(\phi^{(1)}, \phi_{,\alpha}^{(2)}, \phi^{(3),\alpha})$ with $\phi^{(1)} \phi_{,\alpha}^{(2)} \phi^{(3),\alpha}$.
By the argument in Step 3 of the previous section, it suffices to prove the estimate 
$$ \| P_k [\phi^{(1)} \phi_{,\alpha}^{(2)} \phi^{(3),\alpha}] \|_{N[k]}
\lesssim  
\chi^{(2)}_{k = \max(k_1,k_2,k_3)}
\chi^{(2)}_{k_1 \leq \min(k_2,k_3)}
\| \phi^{(1)} \|_{S[k_1]} 
\| \phi^{(2)} \|_{S[k_2]}
\| \phi^{(3)} \|_{S[k_3]}
$$
for all integers $k, k_1, k_2, k_3$ and Schwartz $\phi^{(i)} \in S[k_i]$ for $i=1,2,3$. 

By symmetry we may assume $k_2 \leq k_3$; by scaling we may take $k_3 = 0$. We thus need to show
$$ \| P_k [\phi^{(1)} \phi_{,\alpha}^{(2)} \phi^{(3),\alpha}] \|_{N[k]}
\lesssim  
\chi^{(2)}_{k=\max(k_1,0)}
\chi^{(2)}_{k_1 \leq k_2}
\| \phi^{(1)} \|_{S[k_1]} 
\| \phi^{(2)} \|_{S[k_2]}
\| \phi^{(3)} \|_{S[0]}.
$$

To do this we split into four (slightly overlapping) cases.

\divider{Case 1: ($\phi^{(1)}$ has the lowest frequency) $k_1 < k_2 + O(1)$.}

This is the easiest case as one does not need to obtain the decay $\chi^{(2)}_{k_1 \leq k_2}$.

From Lemma \ref{core-null} we have
$$
\| P_{k'} [\phi_{,\alpha}^{(2)} \phi^{(3),\alpha}] \|_{N[k']}
\lesssim \chi^{(4)}_{k'=0}
\| \phi^{(2)} \|_{S[k_2]}
\| \phi^{(3)} \|_{S[0]}
$$
for all $k'$.  From \eqref{ur-tilde}, \eqref{stil-s} and Littlewood-Paley  decomposition we thus have
$$ \| P_k [\phi^{(1)} \phi_{,\alpha}^{(2)} \phi^{(3),\alpha}] \|_{N[k]}
\lesssim  
\chi^{(4)}_{k=0}
\| \phi^{(1)} \|_{S[k_1]} 
\| \phi^{(2)} \|_{S[k_2]}
\| \phi^{(3)} \|_{S[0]}
$$
which is acceptable.

\divider{Case 2: ($\phi^{(1)}$ has the intermediate frequency) $k_2 \leq k_1+O(1) \leq O(1)$.}

The idea is to play off Lemma \ref{core} (which gives a gain if $\phi_{,\alpha}^{(2)} \phi^{(3),\alpha}$ is far from the cone) against a variant of \eqref{j-improv} (which will gain when $\phi_{,\alpha}^{(2)} \phi^{(3),\alpha}$ is close to the cone).

We may assume that $k = O(1)$ since the left-hand side vanishes otherwise.  We may also assume that $k_2 < -C$ for some large constant $C$ we can use Case 1 otherwise.  By the triangle inequality it thus suffices to show
\be{good-spread}
\| P_k [\phi^{(1)} G_{k'}] \|_{N[k]}
\lesssim  
\chi^{(2)}_{k_1 \leq k_2}
\| \phi^{(1)} \|_{S[k_1]} 
\| \phi^{(2)} \|_{S[k_2]}
\| \phi^{(3)} \|_{S[0]}
\end{equation}
for all $k,k' = O(1)$, where 
\be{gkp-def}
G_{k'} := P_{k'} (\phi_{,\alpha}^{(2)} \phi^{(3),\alpha}).
\end{equation}

Fix $k$, $k'$.  We apply

\begin{lemma}\label{gnk-lemma}  For any $k'$, $k_2$, $k_3$ we have
\be{gnk}
\| P_{k'} (\phi_{,\alpha}^{(2)} \phi^{(3),\alpha}) \|_{N[k']} \lesssim 
\| \phi^{(2)} \|_{S[k_2]}
\| \phi^{(3)} \|_{S[k_3]}
\end{equation}
and
\be{gnk-j}
\| Q_j P_{k'} (\phi_{,\alpha}^{(2)} \phi^{(3),\alpha}) \|_{N[k']} \lesssim \chi^{(4)}_{j = \min(k_2,k_3)} 
\| \phi^{(2)} \|_{S[k_2]}
\| \phi^{(3)} \|_{S[k_3]}
\end{equation}
whenever $j > \min(k_2,k_3)+O(1)$.  In particular we have
\end{lemma}

\begin{proof}  By symmetry we may assume $k_2 \leq k_3$.  We may then assume that $k' \leq k_3 + O(1)$ since the expression vanishes otherwise.

The estimate \eqref{gnk} follows from Lemma \ref{core-null}, so we turn to \eqref{gnk-j}.  This shall be a variant of \eqref{j-improv}.  We may assume that $j > k_2 + 10$ since the claim follows from \eqref{gnk}, Lemma \ref{n-xsb}, and \eqref{f-xsb} otherwise.

Consider the contribution of $Q_{> j-5} \phi^{(3),\alpha}$. By \eqref{f-xsb} and H\"older we can bound this contribution by
$$
2^{(n/2-1)k'} 2^{-j/2} \| \nabla_{x,t} \phi^{(2)} \|_{L^\infty_t L^\infty_x} \| Q_{>j-5} \phi^{(3)} \|_{L^2_t L^2_x}$$
which is acceptable by \eqref{sk-infty}, \eqref{qbound}.  

Now consider the contribution of $Q_{\leq j-5} \phi^{(3),\alpha}$.  We may then freely insert $Q_{>j-10}$ in front of $\phi^{(2)}_{,\alpha}$.  By \eqref{f-xsb} and H\"older we can then bound the previous by
$$
2^{(n/2-1)k'} 2^{-j/2} \| \nabla_{x,t} Q_{>j-10} \phi^{(2)} \|_{L^2_t L^\infty_x} \| Q_{>j-5} \phi^{(3)} \|_{L^2_t L^\infty_x}$$
which is acceptable by \eqref{qbound-strichartz}, \eqref{sk-infty}.  
\end{proof}

From the above lemma we have
$$ \| Q_{>(k_1+k_2)/2} G_{k'} \|_{N[k']} \lesssim \chi^{(4)}_{k_1=k_2}
\| \phi^{(2)} \|_{S[k_2]}
\| \phi^{(3)} \|_{S[0]}
$$
which when combined with \eqref{ur-tilde}, \eqref{stil-s} shows that the contribution of $Q_{>(k_1+k_2)/2} G_{k'}$ is acceptable.  It thus suffices to show that
$$ \| P_k[ \phi^{(1)} Q_{\leq (k_1 + k_2)/2} G_{k'} \|_{N[k]} 
\lesssim \chi^{(2)}_{k_1=k_2}
\| \phi^{(1)} \|_{S[k_1]} 
\| \phi^{(2)} \|_{S[k_2]}
\| \phi^{(3)} \|_{S[0]}.$$
We first consider the contribution of $Q_{>k_1-C} \phi^{(1)}$.  By Lemma \ref{l12} and H\"older we may bound this by
$$
\lesssim \| Q_{>k_1-C} \phi^{(1)} \|_{L^2_t L^\infty_x} \| Q_{\leq (k_1 + k_2)/2} G_{k'} \|_{L^2_t L^2_x}.$$
From \eqref{qbound-strichartz} and Lemma \ref{n-xsb} we can estimate this by
$$
\lesssim 2^{-k_1/2} \| \phi^{(1)} \|_{S[k_1]} 
2^{(k_1+k_2)/4} \| G_{k'} \|_{N[k']}$$
which is acceptable by \eqref{gnk}.

Thus we are left with controlling
$$ \| P_k[ Q_{\leq k_1-C} \phi^{(1)} Q_{\leq (k_1 + k_2)/2} G_{k'} \|_{N[k]},$$
We may freely insert $Q_{\leq k_1 + 10}$ in front of $P_k$, and insert $P_{k_1-10 < \cdot < k_1+10}$
in front of $Q_{\leq k_1-C}$.  By the triangle inequality we can thus estimate the previous by
$$\lesssim \sum_{j \leq (k_1 + k_2)/2} \| P_k Q_{\leq k_1+10} [ P_{k_1-10 < \cdot < k_1+10} Q_{\leq k_1-C} \phi^{(1)} Q_j G_{k'} \|_{N[k]},$$
which by Lemma \ref{core} is bounded by
$$\lesssim \sum_{j \leq (k_1 + k_2)/2} \| P_{k_1-10 < \cdot < k_1+10} Q_{\leq k_1 - C} \phi^{(1)} \|_{S[k_1]}
2^{-j/2} \chi^{(4)}_{j-k_1}
\| Q_j G_{k'} \|_{L^2_t L^2_x}.$$
By Lemma \ref{disposable} we may discard $P_{k_1-10 < \cdot < k_1+10} Q_{\leq k_1 - C}$.  By Lemma \ref{n-xsb} we can then estimate the previous by
$$ \lesssim \sum_{j \leq (k_1 + k_2)/2} \| \phi^{(1)} \|_{S[k_1]}
\chi^{(4)}_{j=k_1} \| G_{k'} \|_{N[k]}$$
which is acceptable by \eqref{gnk}.

\divider{Case 3: ($\phi^{(1)}$ has the highest frequency; $\phi^{(2)}$, $\phi^{(3)}$ have different frequencies) $k_1 \gg 1$ and $k_2 \leq -\delta_2 k_1$.}

This will be similar to Case 2, but uses the separation between $k_2$ and $k_3=0$ as a proxy for the separation between $k_2$ and $k_1$ for the purposes of obtaining the gain $\chi^{(2)}_{k_1 \leq k_2}$.

In this case we may assume $k = k_1 + O(1)$.  It suffices by the triangle inequality to show that
$$ \| P_k [\phi^{(1)} G_{k'}] \|_{N[k]}
\lesssim  
\chi^{(2)}_{k_1 \leq k_2}
\| \phi^{(1)} \|_{S[k_1]} 
\| \phi^{(2)} \|_{S[k_2]}
\| \phi^{(3)} \|_{S[0]}
$$
for all $k'=O(1)$, where $G_{k'}$ is defined in \eqref{gkp-def}; note that $G_{k'}$ vanishes for $k' \neq O(1)$.  

Fix $k'$.  We shall repeat the argument in \eqref{good-spread}.  From Lemma \ref{gnk-lemma} we have 
$$ \| Q_{>k_2/2} G_{k'} \|_{N[k']} \lesssim \chi^{(4)}_{k_2=0}
\| \phi^{(2)} \|_{S[k_2]}
\| \phi^{(3)} \|_{S[0]}.$$
From \eqref{ur-tilde} and the Case 3 hypothesis we thus see that the contribution of $Q_{>k_2/2} G_{k'}$ is acceptable.  It thus remains to estimate 
$$ \sum_{j \leq k_2/2} \| P_k [\phi^{(1)} Q_j G_{k'}] \|_{N[k]}.$$
By Lemma \ref{core} and Lemma \ref{n-xsb} we can estimate this by
$$ \lesssim \sum_{j \leq k_2/2}  \chi^{(4)}_{j=0}
\| G_{k'} \|_{N[k']}
\| \phi^{(1)} \|_{S[k_1]},$$
which is acceptable by \eqref{gnk} and the Case 3 hypothesis.

\divider{Case 4: ($\phi^{(1)}$ has the highest frequency; $\phi^{(2)}$, $\phi^{(3)}$ have comparable frequencies) $k_1 \gg 1$ and $-\delta_2 k_1 \leq k_2 \leq 0$.}

In this case we may assume $k = k_1 + O(1)$, and we reduce to showing
$$
\| P_k [\phi^{(1)} \phi_{,\alpha}^{(2)} \phi^{(3),\alpha}] \|_{N[k]}
\lesssim  
\chi^{(2)}_{k_1 = 0}
\| \phi^{(1)} \|_{S[k_1]} 
\| \phi^{(2)} \|_{S[k_2]}
\| \phi^{(3)} \|_{S[0]}.
$$

In principle this case is easy because both derivatives are on low frequency terms, however by the same token the null structure is much less advantageous in this case.  In fact in two dimensions (where the null structure is crucial) the gain $\chi^{(2)}_{k_1=0}$ is not easy to obtain even when the $\phi^{(i)}$ are all free solutions\footnote{In particular, it cannot be proven solely using $\dot X^{s,b}$ spaces even in the free case, because these spaces just barely fail to exclude the possibility of the low frequencies concentrating in the region Minkowski-orthogonal to the high frequency, which has the effect of keeping the above trilinear expression close to the light cone in frequency space.  A similar difficulty is also encountered in equations of Maxwell-Klein-Gordon or Yang-Mills type near the critical regularity.  This possible concentration can be shown to be impossible by using $L^p$ bilinear estimates for $p>2$, but in low dimensions these estimates are quite difficult (cf. \cite{wolff:cone}, \cite{tao:cone}).  Another approach would be to refine the frequency-sector control of $S(c)$ further, perhaps to small balls, in order to isolate the portions of $\phi^{(i)}$ which are contributing to the bad Minkowski-orthogonal case.  In both cases one would have to re-engineer the space $S(c)$ substantially.  We have chosen a different approach based on Bernstein's inequality.}.  Our strategy will be to use some multilinear multiplier calculus to try to move the derivatives (and the null structure) back onto the high frequency term, hopefully gaining in the process a factor proportional to the ratio between the high and low frequencies.  This turns out to be achievable except when the high frequency is Minkowski-orthogonal to the two low frequencies, but in this case one can repeat the derivation of Lemma \ref{ur-algebra}, improving upon Bernstein's inequality for the low frequencies at one crucial juncture to obtain the desired gain.

We now consider the contribution of
$$ \| P_k [\phi^{(1)} P_{\leq -\delta_2 k_1} (\phi_{,\alpha}^{(2)} \phi^{(3),\alpha})] \|_{N[k]}.$$
This contribution is only non-zero when $k_2 = O(1)$.
By \ref{ur-tilde} and Lemma \ref{core-null} we can bound this by
$$ \lesssim \sum_{k' \leq - \delta_2 k_1} \| \phi^{(1)} \|_{S[k_1]}
\chi^{(2)}_{k' = 0} \| \phi^{(2)} \|_{S[k_2]} \| \phi^{(3)} \|_{S[k_3]}$$
which is acceptable.

Next, consider the contribution of
$$ \| P_k [\phi^{(1)} P_{-\delta_2 k_1 < \cdot < 5} Q_{\leq -2\delta_2 k_1} (\phi_{,\alpha}^{(2)} \phi^{(3),\alpha})] \|_{N[k]}.$$
We estimate this by
$$\lesssim \sum_{-\delta_2 k_1 < k' < 10} \sum_{j \leq -2\delta_2 k_1}
\| P_k [\phi^{(1)} P_{k'} Q_j (\phi_{,\alpha}^{(2)} \phi^{(3),\alpha})] \|_{N[k]}$$
which by Lemma \ref{core}, Lemma \ref{n-xsb} is bounded by
$$\lesssim \sum_{-\delta_2 k_1 < k' < 5} \sum_{j \leq -2\delta_2 k_1}
\chi^{(4)}_{j \geq k'}
\| \phi^{(1)}\|_{S[k_1]}
\| P_{k'} (\phi_{,\alpha}^{(2)} \phi^{(3),\alpha})] \|_{N[k']}
$$
which by Lemma \ref{core-null} is bounded by
$$ \lesssim \sum_{-\delta_2 k_1 < k' < 5} \sum_{j \leq -2\delta_2 k_1}
\chi^{(4)}_{j \geq k'}
\| \phi^{(1)}\|_{S[k_1]}
\| \phi^{(2)}\|_{S[k_2]}
\| \phi^{(3)}\|_{S[k_3]},
$$
which is acceptable.

It thus remains to show
$$
\| P_k [\phi^{(1)} \Pi (\phi_{,\alpha}^{(2)} \phi^{(3),\alpha})] \|_{N[k]}
\lesssim  
\chi^{(2)}_{k_1=k_2}
\| \phi^{(1)} \|_{S[k_1]} 
\| \phi^{(2)} \|_{S[k_2]}
\| \phi^{(3)} \|_{S[0]}
$$
where
$$ \Pi := P_{-\delta_2 k_1 < \cdot < 5} Q_{> -2\delta_2 k_1}.$$
By adapting the proof of Lemma \ref{proj-disposable} we observe that
$\chi^{(2)}_{k_1=0} \Pi$ is disposable for some $C$ depending only on $n$.  Thus we can usually discard $\Pi$ whenever desired by paying a price of $\chi^{-(2)}_{k_1=0}$.  The projection $\Pi$ will be a minor nuisance for the next few steps, but is crucially needed much later in the argument (in Cases 4(e).3(b)-(c)), and is difficult to insert later on.

We now use Littlewood-Paley projections in time to split into five sub-cases, of which the first four are quite easy.  The fifth case is the main one, in which $D_+ \sim D_0$ for all three functions $\phi^{(i)}$.

\divider{Case 4(a): $\phi^{(2)}$, $\phi^{(3)}$ both have Fourier support in the region $|\tau| \gtrsim 1$.}

In this case we use \eqref{f-l12} and H\"older and discard $\Pi$ (paying $\chi^{-(2)}_{k_1=0}$) to estimate this contribution by
\begin{eqnarray*}
&\lesssim \| \phi^{(1)} \Pi (\phi_{,\alpha}^{(2)} \phi^{(3),\alpha}) \|_{L^1_t L^2_x}\\
&\lesssim \| \phi^{(1)} \|_{L^\infty_t L^2_x}
\chi^{-(2)}_{k_1=0}
\| \nabla_{x,t} \phi^{(2)} \|_{L^2_t L^\infty_x}
\| \nabla_{x,t} \phi^{(3)} \|_{L^2_t L^\infty_x}.
\end{eqnarray*}
By \eqref{sk-energy}, \eqref{qbound-strichartz} we may bound this by
$$\lesssim 
2^{-nk_1/2} \| \phi^{(1)} \|_{S[k_1]} \chi^{-(2)}_{k_1=0}
\| \phi^{(2)} \|_{S[k_2]} \| \phi^{(3)} \|_{S[k_3]},$$
which is acceptable.

\divider{Case 4(b): $\phi^{(2)}$, $\phi^{(3)}$ have Fourier support in the regions $|\tau| \lesssim 1$ and $|\tau| \gg 1$ respectively.}

Split $\phi^{(1)} = Q_{<-C} \phi^{(1)} + Q_{\geq -C} \phi^{(1)}$.  

To control the contribution of $Q_{\geq -C}$ we use \eqref{f-l12} and H\"older,
and discard $\Pi$ (paying $\chi^{-(2)}_{k_1=0}$) to estimate this contribution by
$$ \lesssim \| Q_{\geq -C} \phi^{(1)} \Pi (\phi_{,\alpha}^{(2)} \phi^{(3),\alpha}) \|_{L^1_t L^2_x}
\lesssim 
\| Q_{\geq -C} \phi^{(1)} \|_{L^2_t L^2_x}
\chi^{-(2)}_{k_1=0}
\| \nabla_{x,t} \phi^{(2)} \|_{L^\infty_t L^\infty_x}
\| \nabla_{x,t} \phi^{(3)} \|_{L^2_t L^\infty_x}.$$
by \eqref{qbound}, \eqref{sk-infty}, \eqref{qbound-strichartz} this is bounded by
$$
\lesssim 
2^{-nk_1/2} \| \phi^{(1)} \|_{S[k_1]} \chi^{-(2)}_{k_1=0}
\| \phi^{(2)} \|_{S[k_2]} \| \phi^{(3)} \|_{S[k_3]}$$
which is acceptable as before.  

To control the contribution of $Q_{<-C}$, we observe that this term has Fourier support
in the region $D_0 \sim 2^k$, $D_- \gtrsim 1$.  Thus we may use \eqref{f-xsb} to estimate this contribution by
$$ \lesssim 2^{(n/2-1)k} \| Q_{< -C} \phi^{(1)} \Pi (\phi_{,\alpha}^{(2)} \phi^{(3),\alpha}) \|_{L^2_t L^2_x}.
$$
By H\"older and discard $\Pi$ (paying $\chi^{-(2)}_{k_1=0}$) we may estimate this contribution by
$$
\lesssim 2^{(n/2-1)k}
\| Q_{< -C} \phi^{(1)} \|_{L^\infty_t L^2_x}
\chi^{-(2)}_{k_1=0}
\| \nabla_{x,t} \phi^{(2)} \|_{L^\infty_t L^\infty_x}
\| \nabla_{x,t} \phi^{(3)} \|_{L^2_t L^\infty_x}.$$
But this is acceptable by \eqref{sk-energy}, \eqref{sk-infty}, \eqref{qbound-strichartz}.

\divider{Case 4(c): $\phi^{(2)}$, $\phi^{(3)}$ have Fourier support in the regions $|\tau| \gg 1$ and $|\tau| \lesssim 1$ respectively.}

This is similar to Case 4(b) and is left to the reader.  (Some additional factors of $2^{k_2}$ appear but the reader can easily verify that they are negligible since we are in Case 4.)

\divider{Case 4(d): $\phi^{(2)}$, $\phi^{(3)}$ both have Fourier support in the region $|\tau| \lesssim 1$, and $\phi^{(1)}$ has Fourier support on $|\tau| \gg 2^{k_1}$.}

This contribution has Fourier support on $D_0 \sim 2^{k}$, $D_- \gtrsim 2^{k}$, so by \eqref{f-xsb} we may estimate it by
$$ \lesssim 2^{(n/2-1)k} 2^{k} \| \phi^{(1)} 
\Pi (\phi_{,\alpha}^{(2)} \phi^{(3),\alpha}) \|_{L^2_t L^2_x}.$$
Using H\"older and discarding $\Pi$ (paying $\chi^{-(2)}_{k_1=0}$) we may bound this by
$$
\lesssim 
2^{(n/2-1)k} 2^{k}
\| \phi^{(1)} \|_{L^2_t L^2_x}
\chi^{-(2)}_{k_1=0}
\| \nabla_{x,t} \phi^{(2)} \|_{L^\infty_t L^\infty_x}
\| \nabla_{x,t} \phi^{(3)} \|_{L^\infty_t L^\infty_x}.$$
But this is acceptable by \eqref{qbound} and two applications of \eqref{sk-infty}.

\divider{Case 4(e): (No function is too far away from the light cone) $\phi^{(2)}$, $\phi^{(3)}$ both have Fourier support in the region $|\tau| \lesssim 1$, and $\phi^{(1)}$ has Fourier support on $|\tau| \lesssim 2^{k_1}$.}

We now divide the $\phi^{(i)}$ into spacetime sectors of angular width $2^{-\delta_4 k_1}$.  More precisely, we let $\Omega_*$ be a maximal $2^{-\delta_4 k_1}$-separated subset of the set $\{ D_0 = 1; D_+ \lesssim 1\}$
and for each $\Xi = (\tau,\xi)$ in $\Omega_*$, let $\pi_{\Xi}$ be a spacetime Fourier multiplier whose symbol $m_\Xi$ is homogeneous of degree zero, is adapted to the sector
$$ \{ (\tau', \xi'): \angle (\tau',\xi'), \Xi \lesssim 2^{-\delta_4 k_1} \}$$
(where $\angle$ denotes the standard Euclidean angle) and we have the partition of unity
$$ 1 = \sum_{\Xi \in \Omega_*} \pi_\Xi$$
for functions with Fourier support on the region $D_+ \sim D_0$.
We can then bound the contribution of this case by
\be{split-sum}
\lesssim \|
\sum_{\Xi_1, \Xi_2, \Xi_3 \in \Omega_*}
P_k [\pi_{\Xi_1} \phi^{(1)} \Pi( \pi_{\Xi_2} \phi_{,\alpha}^{(2)} \pi_{\Xi_3} \phi^{(3),\alpha})] \|_{N[k]}.
\end{equation}
Observe that the $\pi_{\Xi_i}$ are disposable.

If $\Xi = (\tau,\xi), \Xi' = (\tau',\xi')$ are spacetime vectors, we define the Minkowski inner product $\eta(\Xi,\Xi')$ by
$$ \eta(\Xi,\Xi') := \xi \cdot \xi' - \tau \tau'.$$

We now split the sum into three pieces.

\divider{Case 4(e).1: ($\phi^{(1)}_{,\alpha} \phi^{(2),\alpha}$ is non-degenerate) The portion of the summation where $|\eta(\Xi_1, \Xi_2)| \gg 2^{-\delta_4 k_1}$}

In this case the idea is to move the null form from $\phi^{(2)}_{,\alpha} \phi^{(3),\alpha}$ to the null form $\phi^{(1)}_{,\alpha} \phi^{(2),\alpha}$ (which is now non-degenerate), gaining a factor of about $2^{-k_1}$ in the process.

We turn to the details.  Consider a single summand of \eqref{split-sum}.  The Fourier transform of this at $(\tau,\xi)$ is
\begin{eqnarray*}
C \int
\delta(\tau_1 + \tau_2 + \tau_3)  \delta(\xi_1 + \xi_2 + \xi_3)&
(\xi_2 \cdot \xi_3 - \tau_2 \tau_3)
(1-m_0(2^{\delta_2 k_1} |\xi_2 + \xi_3|))\\
(1-m_0(2^{2\delta_2 k_1} ||\tau_2 + \tau_3| - |\xi_2 + \xi_3||))&
\prod_{i=1}^3 m_{\Xi_i}(\tau_i,\xi_i) {\cal F} \phi^{(i)}(\tau_i,\xi_i)\ d\tau_i d\xi_i.
\end{eqnarray*}
We rewrite this as
$$
C 2^{-k_1} 2^{\delta_4 k_1} \int
\delta(\tau_1 + \tau_2 + \tau_3) \delta(\xi_1 + \xi_2 + \xi_3)
 (\xi_1 \cdot \xi_2 - \tau_1 \tau_2) M
\prod_{i=1}^3 m_{\Xi_i}(\tau_i,\xi_i) {\cal F} \phi^{(i)}(\tau_i,\xi_i)\ d\tau_i d\xi_i$$
where the symbol $M = M(\tau_1,\tau_2,\tau_3,\xi_1,\xi_2,\xi_3)$ is defined by
$$
M := 
(1-m_0(2^{\delta_2 k_1} |\xi_2 + \xi_3|))
(1-m_0(2^{2\delta_2 k_1} ||\tau_2 + \tau_3| - |\xi_2 + \xi_3||))
\frac{2^{k_1} (\xi_2 \cdot \xi_3 - \tau_2 \tau_3)}{2^{\delta_4 k_1} 
(\xi_1 \cdot \xi_2 - \tau_1 \tau_2)}.
$$
The function
$$ \prod_{i=1}^3 m_{\Xi_i}(\tau_i,\xi_i) {\cal F}  \phi^{(i)}(\tau_i,\xi_i)$$
is supported on a $3n+3$-dimensional box, which is the product of three $n+1$-dimensional boxes of dimension $2^{k_i} \times 2^{-\delta_4 k_1} 2^{k_i} \times \ldots 2^{-\delta_4 k_1}$ for $i=1,2,3$.  On this box the symbol $M$ is bounded, and by the hypothesis of Case 4(e).1 can in fact be replaced by a bump function adapted to a dilate of this $3n+3$-dimensional box, except that each normalized derivative could pull down a factor of $\chi^{-(2)}_{k_1=0}$.  In particular, we can then write $M$ on this box as a Fourier series (cf. \cite{tao:boch-rest})
$$ M(\tau_1,\tau_2,\tau_3,\xi_1,\xi_2,\xi_3) := \sum_{X \in \Gamma} c_X \prod_{i=1}^3 e^{2\pi i (t_i \tau_i + x_i \cdot \xi_i)}$$
where $\Gamma$ is some discrete subset of $(\R \times \R^n)^3$ (basically the dual lattice of the dilate of the $3n+3$-dimensional box), $X = (t_1,t_2,t_3,x_1,x_2,x_3)$, and the $c_X$ are a set of co-efficients with an $l^1$ norm of $O(\chi^{-(2)}_{k_1=0})$.  Undoing the Fourier transform, we thus see that
$$
\pi_{\Xi_1} \phi^{(1)}(t,x) \pi_{\Xi_2} \phi_{,\alpha}^{(2)}(t,x) \pi_{\Xi_3} \phi^{(3),\alpha}(t,x)
$$
can be written as
$$
C 2^{-k_1} \chi^{-(4)}_{k_1=0}
\sum_{X \in \Gamma} c_X
\pi_{\Xi_1} \phi^{(1)}_{,\alpha}(t-t_1,x-x_1) \pi_{\Xi_2} \phi^{(2),\alpha}(t-t_2,x-x_2) \pi_{\Xi_3} \phi^{(3)}(t-t_3,x-x_3).$$
The summand in \eqref{split-sum} can thus be estimated by
$$
\lesssim
2^{-k_1} \chi^{-(4)}_{k_1=0}
\sup_{X \in \Gamma} 
\| P_k[\pi_{\Xi_1} \phi^{(1)}_{,\alpha}(t-t_1,x-x_1) \pi_{\Xi_2} \phi^{(2),\alpha}(t-t_2,x-x_2) \pi_{\Xi_3} \phi^{(3)}(t-t_3,x-x_3)] \|_{N[k]}.$$
By Lemma \ref{core-null} and \eqref{ur-tilde} we can bound this by\footnote{It is possible to adapt the proofs of Lemma \ref{core-null} and \eqref{ur-tilde} to apply directly to the summands of \eqref{split-sum} in this case, and thus forego the need to perform the Fourier series decomposition of $M$.  But this would require a tedious repetition the proofs of these two Lemmata - in addition to the repetition which will occur in Case 4(e).3 - so we opted to take this shortcut instead.} 
$$
\lesssim
2^{-k_1} \chi^{-(4)}_{k_1=0}
\sup_{X \in \Gamma}
\prod_{i=1}^3
\| \pi_{\Xi_i} \phi^{(i)}(\cdot - t_i, \cdot - x_i) \|_{S[k_i]}
$$
On the Fourier support of $\phi^{(i)}$, $\pi_{\Xi_i}$ is a multiplier with integrable kernel, so we can estimate this by
$$
\lesssim
2^{-k_1} \chi^{-(4)}_{k_1=0}
\| \phi^{(1)} \|_{S[k_1]}
\| \phi^{(2)} \|_{S[k_2]}
\| \phi^{(3)} \|_{S[k_3]}.
$$
Since there are at most $O(\chi^{-(4)}_{k_1=0})$ summands, we see that the total contribution of these terms are acceptable (since $n=2,3,4$).  

\divider{Case 4(e).2: ($\phi^{(1)}_{,\alpha} \phi^{(3),\alpha}$ is non-degenerate) The portion of the summation not covered by 4(e).1 where $|\eta(\Xi_1, \Xi_3)| \gg 2^{-\delta_4 k_1}$.}

This is done similarly to Case 4(e).1. A few additional powers of $2^{k_2} = O(\chi^{-(2)}_{k_1=0})$ appear in the estimates, but they make no difference to the final argument. 

\divider{Case 4(e).3: (Low frequencies are Minkowski-orthogonal to the high frequency) The portion of the summation where
$\eta(\Xi_1, \Xi_2), \eta(\Xi_1,\Xi_3) = O(2^{-\delta_4 k_1})$.}

We are thus left with showing
\be{mid-split}
\|
\sum_*
P_k [\pi_{\Xi_1} \phi^{(1)} \Pi( \pi_{\Xi_2} \phi_{,\alpha}^{(2)} \pi_{\Xi_3} \phi^{(3),\alpha})] \|_{N[k]}
\lesssim  
\chi^{(2)}_{k_1=0}
\prod_{i=1}^3
\| \phi^{(i)} \|_{S[k_i]},
\end{equation}
where the summation $\sum_*$ is over all $\Xi_1, \Xi_2, \Xi_3 \in \Omega_*$ for which
$\eta(\Xi_1, \Xi_2), \eta(\Xi_1,\Xi_3) = O(2^{-\delta_4 k_1})$.  Note that there are at most $O(\chi^{-(3)}_{k_1=0})$ terms in this summation.

Define the quantity $l := k_1 + \delta_2 k_1$.
We split $\phi^{(1)} = Q_{\geq k_1 -2l} \phi^{(1)} + Q_{< k_1 - 2l} \phi^{(1)}$.  To deal with the latter term we shall also split $P_k = P_k Q_{<k_1-2l} + P_k Q_{\geq k_1-2l}$.  This gives us three contributions to consider:

\divider{Case 4(e).3(a).  (The input $\phi^{(1)}$ is not too close to light cone) The contribution of $Q_{\geq k_1-2l} \phi^{(1)}$.}

By the triangle inequality and \eqref{f-l12} we can estimate this contribution by
$$
\lesssim
\sum_*
2^{(\frac{n}{2}-1)k_1}
\|
\pi_{\Xi_1} Q_{\geq k_1-2l} \phi^{(1)} \Pi( \pi_{\Xi_2} \phi_{,\alpha}^{(2)} \pi_{\Xi_3} \phi^{(3),\alpha}) \|_{L^1_t L^2_x}
$$
We use H\"older and discard $\Pi$ (paying $\chi^{-(2)}_{k_1=0}$) to estimate this by
$$
\lesssim
\sum_*
2^{(\frac{n}{2}-1)k_1}
\|
Q_{\geq k_1-2l} \phi^{(1)} \|_{L^2_t L^2_x}
\chi^{-(2)}_{k_1=0} \| \nabla_{x,t} \phi^{(2)} \|_{L^\infty_t L^\infty_x}
\| \nabla_{x,t} \phi^{(3)}) \|_{L^\infty_t L^\infty_x}
$$
which by \eqref{qbound}, \eqref{sk-infty} can be bounded by
$$
\lesssim 
\chi^{-(4)}_{k_1=0}
2^{(\frac{n}{2}-1)k_1}
2^{-nk_1/2} 2^{-(k_1-2l)/2} \| \phi^{(1)} \|_{S[k_1]}
\chi^{-(2)}_{k_1=0} \| \phi^{(2)} \|_{S[k_2]} \| \phi^{(3)} \|_{S[k_3]}.$$
But this simplifies to
$$
\lesssim 
\chi^{-(4)}_{k_1=0}
2^{-k_1/2}
\| \phi^{(1)} \|_{S[k_1]} \| \phi^{(2)} \|_{S[k_2]} \| \phi^{(3)} \|_{S[k_3]},$$
which is acceptable.  

\divider{Case 4(e).3(b).  (The output is not too close to light cone) The contribution of $Q_{< k_1-2l} \phi^{(1)}$ and $P_k Q_{\geq k_1 - 2l}$.}

We split $Q_{< k_1-2l} = Q^+_{<k_1-2l} + Q^-_{<k_1-2l}$
and $Q_{\geq k_1-2l} = Q^+_{\geq k_1-2l} + Q^-_{\geq k_1-2l}$, and observe that the contribution from opposing signs is zero.  From conjugation symmetry it thus suffices to estimate
$$
\lesssim \|
\sum_*
Q^+_{\geq k_1 - 2l} P_k [\pi_{\Xi_1} Q^+_{< k_1 - 2l} \phi^{(1)} \Pi( \pi_{\Xi_2} \phi_{,\alpha}^{(2)} \pi_{\Xi_3} \phi^{(3),\alpha})] \|_{N[k]}.$$
By the triangle inequality and \eqref{f-xsb} we may estimate this by
$$
\lesssim \sum_*
2^{(\frac{n}{2}-1)k_1} 2^{-(k_1-2l)/2}
\|
P_k [\pi_{\Xi_1} Q^+_{< k_1 - 2l} \phi^{(1)} \Pi( \pi_{\Xi_2} \phi_{,\alpha}^{(2)} \pi_{\Xi_3} \phi^{(3),\alpha})] \|_{L^2_t L^2_x}.$$
We expand out $\Pi$ and estimate this by
\begin{eqnarray*}
\lesssim \sum_* &
2^{(\frac{n}{2}-1)k_1} 2^{-(k_1-2l)/2}
\|
\phi^{(1)} \|_{L^\infty_t L^2_x}\\
&\sum_{-2\delta_2 k_1 < k' < 10}
\| P_{k'} Q_{> -2\delta_2 k_1}( \pi_{\Xi_2} \phi_{,\alpha}^{(2)} \pi_{\Xi_3} \phi^{(3),\alpha}) \|_{L^2_t L^\infty_x}.
\end{eqnarray*}
By \eqref{sk-energy} and Bernstein's inequality \eqref{bernstein} we can bound this by
\begin{eqnarray*}
\lesssim \sum_* &
2^{(\frac{n}{2}-1)k_1} 2^{-(k_1-2l)/2} 2^{-nk_1/2}
\| \phi^{(1)} \|_{S[k_1]}\\
&\sum_{-2\delta_2 k_1 < k' < 10} 2^{nk'/2}
\| P_{k'} Q_{> -2\delta_2 k_1}( \pi_{\Xi_2} \phi_{,\alpha}^{(2)} \pi_{\Xi_3} \phi^{(3),\alpha}) \|_{L^2_t L^2_x}.
\end{eqnarray*}
By Lemma \ref{n-xsb} we can bound this by
\begin{eqnarray*}
\lesssim 
\chi^{-(2)}_{k_1=0}
\sum_* &
2^{(\frac{n}{2}-1)k_1} 2^{-(k_1-2l)/2} 2^{-nk_1/2} 
\| \phi^{(1)} \|_{S[k_1]} \\
& \sum_{-2\delta_2 k_1 < k' < 10} 2^{nk'/2} 2^{-(\frac{n}{2}-1)k'} 
\| P_{k'} ( \pi_{\Xi_2} \phi_{,\alpha}^{(2)} \pi_{\Xi_3} \phi^{(3),\alpha}) \|_{N[k']}
\end{eqnarray*}
which by Lemma \ref{core-null} is bounded by
\begin{eqnarray*}
\lesssim
\chi^{-(2)}_{k_1=0}
\sum_*
&
2^{(\frac{n}{2}-1)k_1} 2^{-(k_1-2l)/2} 2^{-nk_1/2}
\| \phi^{(1)} \|_{S[k_1]}
\\
&\sum_{-2\delta_2 k_1 < k' < 10} 2^{nk'/2} 2^{-(\frac{n}{2}-1)k'} \| \phi^{(2)}\|_{S[k_2]} \| \phi^{(3)} \|_{S[k_3]}.
\end{eqnarray*}
But this simplifies to 
$$
\lesssim 
2^{-k_1/2} \chi^{-(4)}_{k_1=0}
\| \phi^{(1)} \|_{S[k_1]}
\| \phi^{(2)}\|_{S[k_2]} \| \phi^{(3)} \|_{S[k_3]}$$
which is acceptable.

\divider{Case 4(e).3(c).  (The input $\phi^{(1)}$ and the output are both very close to light cone) The contribution of $Q_{< k_1-2l} \phi^{(1)}$ and $P_k Q_{< k_1 - 2l}$.}

We repeat the argument at the start of Case 4(e).3(b), and reduce to showing
$$
\|
\sum_*
Q^+_{< k_1 - 2l} P_k [\pi_{\Xi_1} Q^+_{< k_1 - 2l} \phi^{(1)} \Pi( \pi_{\Xi_2} \phi_{,\alpha}^{(2)} \pi_{\Xi_3} \phi^{(3),\alpha})] \|_{N[k]}
\lesssim  
\chi^{(2)}_{k_1=0}
\prod_{i=1}^3
\| \phi^{(i)} \|_{S[k_i]}.
$$
We estimate the left-hand side by
$$
\lesssim \|
\sum_{\kappa \in K_l}
Q^+_{< k_1 - 2l} P_k P_{k,\kappa} [
\sum_* \sum_{\kappa' \in K_l}
\pi_{\Xi_1} Q^+_{< k_1 - 2l} P_{k_1,\kappa'} \phi^{(1)} \Pi( \pi_{\Xi_2} \phi_{,\alpha}^{(2)} \pi_{\Xi_3} \phi^{(3),\alpha})] \|_{N[k]}.
$$
Observe from the geometry of the cone that we may assume that 
$$ 2^{-l} \ll \dist(\kappa,\kappa') \lesssim 2^{-k_1}$$
since the other terms have a zero contribution.  By \eqref{f-null} we may thus estimate the previous by
\begin{eqnarray*}
\lesssim 2^{(\frac{n}{2}-1)k_1}
(\sum_{\kappa \in K_l}
\|
Q^+_{< k_1 - 2l} & P_k P_{k,\kappa} [ 
\sum_* \sum_{\kappa' \in K_l: 2^{-l} \ll \dist(\kappa,\kappa') \lesssim 2^{-k_1}}
\\
&\pi_{\Xi_1} Q^+_{< k_1 - 2l} P_{k_1,\kappa'} \phi^{(1)} \Pi( \pi_{\Xi_2} \phi_{,\alpha}^{(2)} \pi_{\Xi_3} \phi^{(3),\alpha})] \|_{NFA[\kappa]}^2)^{1/2}.
\end{eqnarray*}
By Lemma \ref{proj-disposable} we can discard $Q^+_{< k_1 - 2l} P_k P_{k,\kappa}$.

Let $\theta_\kappa$ denote the null direction 
$$ \theta_\kappa := (1,\omega_\kappa).$$
Note that we may assume that $\angle(\Xi_1, \theta_\kappa) \lesssim 2^{-\delta_4 k_1}$ since the other terms have a zero contribution.  We therefore see that for each $\kappa$ there are only $O(1)$ values of $\Xi_1$ which contribute.  Also for each $\kappa$ there are only $O(\chi^{-(2)}_{k_1=0})$ values of $\kappa'$ which contribute.  We may therefore estimate the previous by
\begin{eqnarray*}
\lesssim 2^{(\frac{n}{2}-1)k_1} \chi^{-(2)}_{k_1=0}
(\sum_{\kappa \in K_l} &
\sum_{\kappa' \in K_l: 2^{-l} \ll \dist(\kappa,\kappa') \lesssim 2^{-k_1}}
\sum_{\Xi_1 \in \Omega_*: \angle(\Xi_1, \theta_\kappa) \lesssim 2^{-\delta_4 k_1}}\\
&\|
\pi_{\Xi_1} Q^+_{< k_1 - 2l} P_{k_1,\kappa'} \phi^{(1)} \Pi( 
\sum_{**} 
\pi_{\Xi_2} \phi_{,\alpha}^{(2)} \pi_{\Xi_3} \phi^{(3),\alpha}) \|_{NFA[\kappa]}^2)^{1/2}
\end{eqnarray*}
where for each $\Xi_1$, $\sum_{**}$ ranges over all the pairs $(\Xi_2, \Xi_3)$ which would have contributed to $\sum_*$.
By \eqref{NFAPW-dual}, then discarding $\pi_{\Xi_1}$, we can estimate the previous by
\begin{eqnarray*}
\lesssim 2^{(\frac{n}{2}-1)k_1} &\chi^{-(2)}_{k_1=0} 2^l 2^{-(n-1)l/2} 2^{-k_1/2}
(\sum_{\kappa \in K_l}
\sum_{\kappa' \in K_l: 2^{-l} \ll \dist(\kappa,\kappa') \lesssim 2^{-k_1}}
\sum_{\Xi_1 \in \Omega_*: \angle \Xi_1, \theta \lesssim 2^{-\delta_4 k_1}}\\
& \|
Q^+_{< k_1 - 2l} P_{k_1,\kappa'} \phi^{(1)} 
\|_{S[k_1,\kappa']}
\| \Pi( 
\sum_{**} 
\pi_{\Xi_2} \phi_{,\alpha}^{(2)} \pi_{\Xi_3} \phi^{(3),\alpha}) \|_{L^2_t L^2_x}^2)^{1/2}.
\end{eqnarray*}
Assume for the moment that we can show
\be{sandwich}
\| \Pi( 
\sum_{**} 
\pi_{\Xi_2} \phi_{,\alpha}^{(2)} \pi_{\Xi_3} \phi^{(3),\alpha}) \|_{L^2_t L^2_x}
\lesssim \chi^{(4)}_{k_1=0} \| \phi^{(2)} \|_{S[k_2]} \| \phi^{(3)} \|_{S[k_3]}.
\end{equation}
for all $\kappa$, $\kappa'$, $\Xi_1$ in the above summation.
Substituting this itno the previous and evaluating the $\Xi_1$ and $\kappa$ summations using \eqref{sk-def} we obtain the bound
$$
\lesssim 2^{(\frac{n}{2}-1)k_1} \chi^{-(2)}_{k_1=0} 2^l 2^{-(n-1)l/2} 2^{-k_1/2}
\chi^{-(2)}_{k_1=0}
\| \phi^{(1)} \|_{S[k_1]}
\chi^{(4)}_{k_1=0} \| \phi^{(2)} \|_{S[k_2]} \| \phi^{(3)} \|_{S[k_3]}.
$$
But this simplifies to
$$
\lesssim \chi^{(4)}_{k_1=0}
\| \phi^{(1)} \|_{S[k_1]}
\| \phi^{(2)} \|_{S[k_2]} \| \phi^{(3)} \|_{S[k_3]}
$$
which is acceptable.

It remains to prove \eqref{sandwich}.  This is an estimate similar to those in Lemma \ref{core-null}, but the conditions of Case 4(e).3 will allow us to squeeze out the additional $\chi^{(3)}_{k_1=0}$ gain.

Fix $\kappa$, $\kappa'$, $\Xi_1$.  We can factorize the left-hand side of \eqref{sandwich} as
\be{sand-2}
\| \Pi( 
[\sum_{\Xi_2 \in \Omega_{**}} 
\pi_{\Xi_2} \phi_{,\alpha}^{(2)} ]
[\sum_{\Xi_3 \in \Omega_{**}}
\pi_{\Xi_3} \phi^{(3),\alpha})] \|_{L^2_t L^2_x}
\end{equation}
where
$$
\Omega_{**} := \{ \Xi \in \Omega_*: \eta(\Xi_1,\Xi) = O(2^{-\delta_4 k_1}) \}.$$
Write $\Xi_i = (\tau_i, \xi_i)$ for $i=1,2,3$.  From the hypotheses
$\eta(\Xi_1,\Xi_2), \eta(\Xi_1,\Xi_3), \angle \Xi_1, \theta_\kappa = O(2^{-\delta_4 k_1})$ we see that
\be{2-plane}
\tau_2 = \xi_2 \cdot \omega_\kappa + O(2^{-\delta_4 k_1});
\quad \tau_3 = \xi_3 \cdot \omega_\kappa + O(2^{-\delta_4 k_1}).
\end{equation}

Define $\xi'_i := \xi_i - (\xi_i \cdot \omega_\kappa) \omega_\kappa$ for $i=2,3$, thus $\xi'_i$ is the portion of $\xi_i$ orthogonal to $\omega_\kappa$.
We now split the summation into four cases.

\divider{Case 4(e).3(c).1: ($\phi^{(2)}$, $\phi^{(3)}$ aligned with $\phi^{(1)}$) $|\xi'_2| < 2^{-\delta_3 k_1}$ and $|\xi'_3| < 2^{-\delta_3 k_1}$.}

In this case we see from \eqref{2-plane} that
$$ |\tau_2 + \tau_3| = |\xi_2 + \xi_3| + O(2^{-\delta_3 k_1})$$
and therefore that the contribution of this case vanishes once $\Pi$ is applied.

\divider{Case 4(e).3(c).2: ($\phi^{(2)}$ is aligned with $\phi^{(1)}$, but $\phi^{(3)}$ is not) $|\xi'_2| < 2^{-\delta_3 k_1}$ and $|\xi'_3| \geq 2^{-\delta_3 k_1}$.}

In this case we may assume that $|\xi'_3| \geq 2^{-C \delta_2 k_1}$, since otherwise this contribution vanishes by the arguments in Case 1.  We then discard $\Pi$  (paying $\chi^{-(2)}_{k_1=0}$) and estimate \eqref{sand-2} by
\be{jusenkyo}
\lesssim \chi^{-(2)}_{k_1=0}
\| \sum_{\Xi_2 \in \Omega_{**}: |\xi'_2| < 2^{-\delta_3 k_1} } \pi_{\Xi_2} \nabla_{x,t} \phi^{(2)} \|_{L^\infty_t L^\infty_x}
\|
\sum_{\Xi_3 \in \Omega_{**}: |\xi'_3| \geq 2^{-C\delta_2 k_1} } \pi_{\Xi_3} \nabla_{x,t} \phi^{(3)} \|_{L^2_t L^2_x}.
\end{equation}
Consider the first factor of \eqref{jusenkyo}.  By Bernstein's inequality \eqref{bernstein}, taking advantage of the restriction $|\xi'_2| < 2^{-\delta_3 k_10}$, we may bound this by
$$ \lesssim \chi^{(4)}_{k_1=0}
\| \sum_{\Xi_2 \in \Omega_{**}: |\xi'_2| < 2^{-\delta_3 k_1} } \pi_{\Xi_2} \nabla_{x,t} \phi^{(2)} \|_{L^\infty_t L^2_x}.
$$
By \eqref{2-plane} we see that the Fourier supports of $\pi_{\Xi_2}$ have finitely overlapping $\xi$-projections as $\Xi_2$ varies across $\Omega_{**}$.  Thus we can use orthogonality to estimate the previous by
$$ \lesssim \chi^{(4)}_{k_1=0}
\| \nabla_{x,t} \phi^{(2)} \|_{L^\infty_t L^2_x} \lesssim
\chi^{(4)}_{k_1=0} 2^{nk_2/2}
\| \phi^{(2)} \|_{S[k_2]}
$$
by \eqref{sk-energy}.

Now consider the second factor of \eqref{jusenkyo}.  From \eqref{2-plane} and the restriction $|\xi'_3| \geq 2^{-C\delta_2 k_1}$ we see that we may freely insert a factor of $Q_{>- C\delta_2 k_1}$ inside the norm.  We then use the almost orthogonality of the $\pi_{\Xi_3}$ as $\Xi_3$ varies over $\Omega_{**}$ to estimate this factor by
$$
\lesssim \|
Q_{>-C\delta_2 k_1} \nabla_{x,t} \phi^{(3)} \|_{L^2_t L^2_x}
$$
which is $O(\chi^{-(2)}_{k_1=0} \| \phi^{(3)} \|_{S[0]})$
by \eqref{qbound}.  Combining this with the previous we see that the contribution of this case is acceptable.

\divider{Case 4(e).3(c).3: ($\phi^{(3)}$ is aligned with $\phi^{(1)}$, but $\phi^{(2)}$ is not) $|\xi'_2| \geq 2^{-\delta_3 k_1}$ and $|\xi'_3| < 2^{-\delta_3 k_1}$.}

This is similar to Case 2 and is left to the reader (the factor $2^{nk_2/2} = O(\chi^{(2)}_{k_1=0})$ moves around but is irrelevant). 

\divider{Case 4(e).3(c).4: (Neither $\phi^{(2)}$ nor $\phi^{(3)}$ is aligned with $\phi^{(1)}$) $|\xi'_2| \geq 2^{-\delta_3 k_1}$ and $|\xi'_3| \geq 2^{-\delta_3 k_1}$.}

In this case we discard
$ \Pi$ (paying $\chi^{-(2)}_{k_1=0}$) and estimate \eqref{sand-2} by
$$
\lesssim \chi^{-(2)}_{k_1=0}
\| \sum_{\Xi_2 \in \Omega_{**}: |\xi'_2| \geq 2^{-\delta_3 k_1} } \pi_{\Xi_2} \nabla_{x,t} \phi^{(2)} \|_{L^2_t L^2_x}
\|
\sum_{\Xi_3 \in \Omega_{**}: |\xi'_3| \geq 2^{-\delta_3 k_1}} \pi_{\Xi_3} \nabla_{x,t} \phi^{(3)} \|_{L^\infty_t L^\infty_x}.
$$
By \eqref{2-plane}, the second factor has Fourier support in a $O(2^{-\delta_4 k_1}) \times O(1) \times \ldots \times O(1)$ slab, so by Bernstein's inequality \eqref{bernstein} we can bound the previous by
$$
\lesssim \chi^{-(2)}_{k_1=0}
\| \sum_{\Xi_2 \in \Omega_{**}: |\xi'_2| \geq 2^{-\delta_3 k_1} } \pi_{\Xi_2} \nabla_{x,t} \phi^{(2)} \|_{L^2_t L^2_x}
\chi^{(4)}_{k_1=0}
\|
\sum_{\Xi_3 \in \Omega_{**}: |\xi'_3| \geq 2^{-\delta_3 k_1}} \pi_{\Xi_3} \nabla_{x,t} \phi^{(3)} \|_{L^2_t L^2_x}.
$$
By \eqref{2-plane} and the constraints $|\xi'_2|, |\xi'_3| \geq 2^{-\delta_3 k_10}$ we may freely insert $Q_{>-C \delta_3 k_1}$ into both factors.  By orthogonality and \eqref{qbound} we can estimate the previous by
\begin{eqnarray*}
\lesssim &\chi^{-(2)}_{k_1=0}
\| \nabla_{x,t} Q_{>-C \delta_3 k_1} \phi^{(2)} \|_{L^2_t L^2_x}
\chi^{(4)}_{k_1=0}
\|
\nabla_{x,t} Q_{>-C \delta_3 k_1} \phi^{(3)} \|_{L^2_t L^2_x}\\
\lesssim &
\chi^{-(2)}_{k_1=0}
\chi^{-(3)}_{k_1=0} \| \phi^{(2)} \|_{S[k_2]} \chi^{(4)}_{k_1=0} \chi^{(4)}_{k_1=0} \chi^{-(3)}_{k_1=0} \| \phi^{(3)} \|_{S[0]}
\end{eqnarray*}
which is acceptable.  This concludes the proof of \eqref{sandwich} and therefore of 
\eqref{o-lemma}.
\endprf

\end{document}